\documentclass[12pt,a4paper]{article}
\usepackage{latexsym,fullpage,amsfonts,amssymb,amsmath,amscd,graphics,epic,theorem}
\usepackage[all]{xy}
\baselineskip 18pt \textwidth 15cm  \theoremstyle{plain}
\newtheorem{lemma}{Lemma}[section]
\newtheorem{proposition}[lemma]{Proposition}
\newtheorem{remark}[lemma]{Remark}
\newtheorem{example}[lemma]{Example}
\newtheorem{theorem}{Theorem}
\newtheorem{definition}[lemma]{Definition}

{\theorembodyfont{\rmfamily} \font\newf=cmr10
\begin{document}
\newcommand{\pperp}{\hbox{$\perp\hskip-6pt\perp$}}
\newcommand{\N}{{\mathbb N}}
\newcommand{\PP}{{\mathbb P}}
\newcommand{\Z}{{\mathbb Z}}
\newcommand{\Q}{{\mathbb Q}}
\newcommand{\R}{{\mathbb R}}
\newcommand{\C}{{\mathbb C}}
\newcommand{\K}{{\mathbb K}}
\newcommand{\F}{{\mathbb F}}
\newcommand{\proofend}{\hfill$\Box$\bigskip}
\newcommand{\eps}{{\varepsilon}}
\newcommand{\ko}{{\mathcal O}}
\newcommand{\wx}{{\widetilde x}}
\newcommand{\wz}{{\widetilde z}}
\newcommand{\wa}{{\widetilde a}}
\newcommand{\bz}{{\boldsymbol z}}
\newcommand{\bp}{{\boldsymbol p}}
\newcommand{\wy}{{\widetilde y}}
\newcommand{\wc}{{\widetilde c}}
\newcommand{\bi}{{\omega}}
\newcommand{\bx}{{\boldsymbol x}}
\newcommand{\Log}{{\operatorname{Log}}}
\newcommand{\pr}{{\operatorname{pr}}}
\newcommand{\Graph}{{\operatorname{Graph}}}
\newcommand{\jet}{{\operatorname{jet}}}
\newcommand{\Tor}{{\operatorname{Tor}}}
\newcommand{\sqh}{{\operatorname{sqh}}}
\newcommand{\const}{{\operatorname{const}}}
\newcommand{\Arc}{{\operatorname{Arc}}}
\newcommand{\Sing}{{\operatorname{Sing}}}
\newcommand{\Span}{{\operatorname{Span}}}
\newcommand{\Aut}{{\operatorname{Aut}}}
\newcommand{\Ker}{{\operatorname{Ker}}}
\newcommand{\Int}{{\operatorname{Int}}}
\newcommand{\Aff}{{\operatorname{Aff}}}
\newcommand{\Area}{{\operatorname{Area}}}
\newcommand{\val}{{\operatorname{Val}}}
\newcommand{\conv}{{\operatorname{conv}}}
\newcommand{\rk}{{\operatorname{rk}}}
\newcommand{\ow}{{\overline w}}
\newcommand{\ov}{{\overline v}}
\newcommand{\ks}{{\cal S}}
\newcommand{\red}{{\operatorname{red}}}
\newcommand{\kc}{{\cal C}}
\newcommand{\ki}{{\cal I}}
\newcommand{\kj}{{\cal J}}
\newcommand{\ke}{{\cal E}}
\newcommand{\kz}{{\cal Z}}
\newcommand{\tet}{{\theta}}
\newcommand{\Del}{{\Delta}}
\newcommand{\bet}{{\beta}}
\newcommand{\mm}{{\mathfrak m}}
\newcommand{\kap}{{\kappa}}
\newcommand{\del}{{\delta}}
\newcommand{\sig}{{\sigma}}
\newcommand{\alp}{{\alpha}}
\newcommand{\Sig}{{\Sigma}}
\newcommand{\Gam}{{\Gamma}}
\newcommand{\gam}{{\gamma}}
\newcommand{\Lam}{{\Lambda}}
\newcommand{\lam}{{\lambda}}
\title{Patchworking singular algebraic
curves, non-Archimedean amoebas and enumerative geometry}
\author{E. Shustin\thanks{{\it AMS Subject Classification}:
Primary 14H15. Secondary 12J25, 14H20, 14M25, 14N10}
\thanks{The
author was supported by Grant No. G-616-15.6/99 from the
German-Israeli Foundation for Research and Development, by the
Hermann-Minkowski Minerva Center for Geometry at Tel Aviv
University and by the Bessel research award from the Alexander von
Humboldt Foundation.}}
\date{}
\maketitle
\begin{abstract} We demonstrate a tropical approach to enumeration of
singular curves on toric surfaces, which consists of reducing the
enumeration of algebraic curves to enumeration of non-Archimedean
amoebas, the images of algebraic curves by a real-valued
non-Archimedean valuation. This idea was proposed by Kontsevich
and recently realized by Mikhalkin, who enumerated nodal curves on
toric surfaces \cite{M2}. We give a detailed algebraic-geometric
explanation for the correspondence between nodal curves and their
amoebas. Our main technical tool is a new patchworking theorem for
singular algebraic curves. We also treat the case of curves with a
cusp and the case of real nodal curves.
\end{abstract}

\vskip10pt

{\newf \hskip2.23in ''... the naive approach with indeterminate coefficients\\[-20pt]

\hskip2.53in and the implicit function theorem."

\hskip3in {\it From a referee report.}}

\vskip10pt

\section{Introduction}

The rapid development of tropical algebraic geometry over the
recent two years has led to interesting applications in
enumerative geometry of singular algebraic curves, proposed by
Kontsevich (see \cite{M}). The first result in this direction has
been obtained by Mikhalkin \cite{M2,M3}, who counted curves with a
given number of nodes on toric surfaces via lattice paths in
convex lattice polygons. The main goal of the present paper is to
explain this breakthrough result and to give a detailed proof for
the link between nodal curves and non-Archimedean amoebas, which
is the core of the tropical approach to enumerative geometry. Our
point of view is purely algebraic-geometric and differs from
Mikhalkin's one, which is based on symplectic geometry techniques.

{\bf Tropical approach to enumerative geometry.} Let
$\Del\subset\R^2$ be a convex lattice polygon, $\Tor_\K(\Del)$ the
toric surface associated with the polygon $\Del$ and defined over
an algebraically closed field $\K$ of characteristic zero. Denote
by $\Lam_\K(\Del)$ the linear system on $\Tor_\K(\Del)$, generated
by the monomials $x^iy^j$, $(i,j)\in\Del\cap\Z^2$. We would like
to count $n$-nodal curves belonging to $\Lam_\K(\Del)$ and passing
through $r=\dim\Lam_\K(\Del)-n=|\Del\cap\Z^2|-1-n$ generic points
in $\Tor_\K(\Del)$. The required number is just the degree of the
so-called Severi variety $\Sig_\Del(nA_1)$. We choose $\K$ to be
the field of convergent Puiseux series over $\C$, i.e., power
series of the form $b(t)=\sum_{\tau\in R}c_\tau t^\tau$, where
$R\subset\R$ is contained in the sum of finitely many bounded from
below arithmetic progressions, and $\sum_{\tau\in
R}|c_\tau|t^\tau<\infty$ for sufficiently small positive $t$. The
latter field is equipped with a non-Archimedean valuation
$\val(b)=-\min\{\tau\in R\ :\ c_\tau\ne 0\}$, which takes $\K^*$
onto $\R$ and satisfies
$$\val(ab)=\val(a)+\val(b),\quad \val(a+b)\le\max\{\val(a),\ \val(b)\},\quad
a,b\in\K^*\ .$$

A curve $C\in\Lam_\K(\Del)$ with $n$ nodes is given by a
polynomial
\begin{equation}f(x,y)=\sum_{(i,j)\in\Del\cap\Z^2}a_{ij}(t)x^iy^j,\quad
a_{ij}(t)\in\K\ .\label{e50}\end{equation} Choosing points
$(x_i,y_i)\in(\K^*)^2$, $i=1,...,r$, so that the exponents of $t$
in $x_i,y_i$ are integral, and imposing conditions $f(x_i,y_i)=0$,
$i=1,...,r$, we shall necessarily have only integral exponents for
$t$ in $a_{ij}(t)$, $(i,j)\in\Del$. Thus, polynomial (\ref{e50})
defines an analytic surface $X$ in $Y=\Tor(\Del)\times
(D\backslash\{0\})$\footnote{From now on the symbol $\Tor(*)$
always means a toric variety over $\C$.}, $D$ being a small disc
in $\C$ centered at $0$, such that the fibres $X_t$ are complex
algebraic curves, which belong to the linear system $\Lam(\Del)$
on the surface $\Tor(\Del)$, and have $n$ nodes (cf. Lemma
\ref{l9}, section \ref{sec7}).

To the pair $(\Tor_\K(\Del),C)$ we assign a certain limit of the
family $(Y_t,X_t)$ as $t\to 0$, where
$Y_t=\Tor(\Del)\times\{t\}\subset Y$. The result $(Y_0,X_0)$ of
this operation we call the {\it tropicalization} (or {\it
dequantization}) of the pair $(\Tor_\K(\Del),C)$. Namely, the
surface $Y_0$ splits into irreducible components
$Y_{0,1},...,Y_{0,N}$, corresponding to a subdivision of $\Del$
into convex lattice polygons, and this subdivision is dual to the
{\it non-Archimedean amoeba} $A_f\subset\R^2$ of the polynomial
$f$, which passes through the points
$(\val(x_i),\val(y_i))\in\R^2$. Next we define a {\it refinement
of the tropicalization} as the tropicalization of the polynomial
$f$ after a certain change of coordinates. The refinement
corresponds to (weighted) blow-ups of the threefold $\overline
Y=Y\cup Y_0$ at some singular points of $X_0$ or along multiple
components of $X_0$, and extends $Y_0$ by adding exceptional
divisors and extends the curve $X_0$ by adding new components,
which we call deformation patterns.

We show that the refined tropicalizations $(Y_0,X_0)$ of $n$-nodal
curves $C\in\Tor_\K(\Del)$ passing through $(x_i,y_i)\in(\K^*)^2$,
$i=1,...,r$, belong to a certain finite set $T$. Using our {\it
patchworking} theorem we decide how many $n$-nodal curves
$C\in\Tor_\K(\Del)$ passing through $(x_i,y_i)\in(\K^*)^2$,
$i=1,...,r$, arise from an element $(Y_0,X_0)$ of $T$, and thus,
we obtain $\deg\Sig_\Del(nA_1)$ as the sum of weights of elements
of $T$. In fact, we look for the family $X_t$ in the form
(\ref{e50}), in which the tropicalization provides some initial
terms in the coefficients $a_{ij}(t)$.

Here we do not touch the merely combinatorial (and, in fact,
elementary) problem to count the elements of $T$. Mikhalkin
\cite{M2} has found a nice way to do this, tracing non-Archimedean
amoebas through points on a straight line and attaching the dual
subdivisions of the amoebas to lattice paths in $\Del$.

We would also like to point out that the tropical approach can be
applied to counting curves with other singularities, and here we
demonstrate this for a relatively simple case of curves with an
ordinary cusp. The main difficulty in the general case is to
describe possible tropicalizations, whereas the patchworking
Theorem \ref{t1} applies to curves with arbitrary singularities.

Furthermore, if the given points in $(\K^*)^2$ are invariant with
respect to the complex conjugation, one can count real
tropicalizations and thus the real singular curves passing through
the given points. We discuss this in section \ref{sec301} in
connection with the Welschinger invariant \cite{Wel}.

{\bf Patchworking construction.} In 1979-80, O. Viro
\cite{Vi1,Vi2,Vi3,Vi4} invented a patchworking construction for
real non-singular algebraic hypersurfaces. We would like to
mention that almost all known topological types of real
non-singular algebraic curves are realized in this way.

In general, the initial data of the construction consist of
\begin{itemize}\item a one-parametric flat family $F\to(\F,0)$ of algebraic
varieties $Y_t$ of dimension $\ge 2$, with $\F=\C$ or $\R$, where
$Y_0$ is assumed to be reduced reducible, and $Y_t$, $t\ne 0$,
irreducible, \item a line bundle ${\cal L}$ on $Y$, \item the zero
locus $X_0\subset Y_0$ of some section $S$ of ${\cal
L}\big|_{Y_0}$, which is assumed to be a hypersurface in
$Y_0$.\end{itemize} The construction extends $S$ up to a section
of ${\cal L}$, whose zero locus $X\subset Y$ defines a family of
hypersurfaces $X_t\subset Y_t$, which inherit some properties of
$X_0$. In \cite{Vi1,Vi2,Vi3,Vi4}, $Y$ is a toric variety
associated with a convex lattice polytope and fibred into toric
hypersurfaces $Y_t$, $t>0$, which degenerate into the union of
some divisors on $Y$, corresponding to facets of the polytope, and
$X_0$ is the union of {\it non-singular} real algebraic
hypersurfaces. The real non-singular hypersurfaces $X_t\subset
Y_t$, $t\ne 0$, appear as a result of a topological gluing
(patchworking) of the components of $X_0$.

In the early 1990's the author suggested to use the patchworking
construction for tracing other properties of objects defined by
polynomials, for example, prescribed singularities of algebraic
hypersurfaces \cite{Sh1,Sh2,Sh3}, critical points of polynomials
\cite{Sh4,Sh2}, singular points and limit cycles of planar
polynomial vector fields \cite{IS1}, resultants of bivariate
polynomials \cite{Sh5}. Considering the patchworking of singular
algebraic curves (i.e., $\dim Y=3$, $\dim Y_t=2$, $\dim X_0=1$) in
\cite{Sh1,Sh2,Sh3}, we always supposed that the components of the
curve $X_0$ are reduced and meet the intersection lines of the
components of the surface $Y_0$ transversally at their
non-singular points. The novelty of the patchworking theorem
presented in this paper (Theorem \ref{t1}, section \ref{sec4}) is
that we allow $X_0$ to be non-reduced and to have singularities
along $\Sing(Y_0)$.

In this connection we would like to point out that, in \cite{Ch}
(see also \cite{CC}), a deformation $Y\to(\C,0)$ of surfaces in
$\PP^3$ with reducible $Y_0$ was considered, where the components
of $X_0$ are nodal curves tangent to the intersection lines of the
components of $Y_0$. For example, Theorem 2.1 in \cite{Ch}, claims
that a point on the intersection line of two components of $Y_0$,
at which non-singular germs of the corresponding components of
$X_0$ have contact of order $m$, gives rise to $m-1$ nodes of
$X_t\subset Y_t$, $t\ne 0$, and the proof is based on a
technically tricky result by Caporaso and Harris \cite{CH}, Lemma
4.1. Our approach is to interpret this as a patchworking, i.e., a
replacement of a neighborhood of a singular point by some
algebraic curve, or more precisely, by an affine curve with Newton
triangle $\{(0,0),(0,2),(m,1)\}$ which can have any number $0\le
k\le m-1$ of nodes (cf. \cite{Sh1}, Proposition 2.5). Extensive
development of this idea covering a broad class of possible
singularities, is done in \cite{ST}. However, the result of
\cite{ST} is not sufficient, for example, for patchworking nodal
curves as required in the enumerative problem.

{\bf Organization of the material.} In the first section we
provide preliminary information on non-Archimedean amoebas and
tropicalizations of polynomials. The second section contains
Theorem \ref{t5}, which reduces the enumeration of nodal curves in
toric surfaces, associated with convex lattice polygons, to the
count of nodal non-Archimedean amoebas passing through the
respective number of generic points in the real plane. The third
section contains Theorem \ref{t6} reducing the enumeration of
curves with one cusp to the count of appropriate cuspidal
non-Archimedean amoebas. In the proof of Theorems \ref{t5} and
\ref{t6} we formulate explicit patchworking statements, which
invert the tropicalization procedure, and which follow from the
main patchworking Theorem \ref{t1}, presented in section
\ref{sec8}. Finally, in section \ref{sec301} we demonstrate an
application of our technique in the computation of the Welschinger
number for real nodal curves in toric surfaces.

{\bf Acknowledgment}. I am very grateful to G. Mikhalkin and I.
Itenberg for useful discussions. I also wish to thank
Universit\"at Kaiserslautern for its hospitality and excellent
working conditions.

\section{Non-Archimedean amoebas}

\subsection{Preliminaries} Amoebas of complex algebraic hypersurfaces have been
introduced in \cite{GKZ} and studied further in
\cite{FPT,GKZ,H,M1,M,MR,PR,R}. We are interested in
``non-Archimedean amoebas", i.e., defined over fields with a
non-Archimedean valuation (see \cite{KT,M}). The field $\K$ of
convergent Puiseux series over $\C$ serves as an example. For a
non-empty finite set $I\subset\Z^k$, denote by $F_\K(I)$ the set
of Laurent polynomials
$$f(\bz)=\sum_{\bi\in I}c_\bi\bz^\bi,\quad \bz=(z_1,...,z_k),\quad
c_\bi\in\K^*,\ \bi\in I\ .$$ Put $Z_f=\{f=0\}\subset(\K^*)^k$ and
define the amoeba of $F$ as
$$A_f=\val(Z_f)\subset\R^k,\quad\text{where}\quad \val(z_1,...,z_k)=(\val(z_1),...,\val(z_n))\
.$$ We denote the set of amoebas $A_f$, $f\in F_\K(I)$ by ${\cal
A}(I)$. If $I$ is the set of all integral points in a convex
lattice polygon $\Del$, we write ${\cal A}(\Del)$.

The following simple observation, which we supply with a proof, is
due to Kapranov \cite{K}.

\begin{theorem}\label{t3} Amoeba $A_f$ coincides with the
corner locus of the piece-wise linear convex function
$$N_f(\bx)=\max_{\bi\in I}(\bi\bx+\val(c_\bi)),\quad\bx\in\R^k\ .$$
(Here and further on, product of vectors means the standard scalar
product.)
\end{theorem}

{\bf Proof}. Let $\bz\in Z_f$, i.e., $\sum_{\bi\in
I}c_\bi\bz^\bi=0$, and hence
$$c_{\bi_0}\bz^{\bi_0}=-\sum_{\bi\in
I\backslash\{\bi_0\}}c_\bi\bz^\bi\ ,$$ with $\bi_0\in I$ chosen so
that
$$\val(c_{\bi_0}\bz^{\bi_0})=\val(c_{\bi_0})+\bi_0\val(\bz)=\min_{\bi\in
I}(\val(c_\bi)+\bi \val(\bz))\ .$$ Thus, there are some
$\bi_1,...,\bi_r\in I\backslash\{\bi_0\}$, $r\ge 1$, such that
$$\val(c_{\bi_0}\bz^{\bi_0})=...=\val(c_{\bi_r}\bz^{\bi_r})>\val(c_\bi\bz^\bi),
\quad\bi\in I\backslash\{\bi_0,...,\bi_r\}$$ or, equivalently,
\begin{eqnarray}&\val(c_{\bi_0})+\bi_0\val(\bz)=...=\val(c_{\bi_r})+\bi_r\val(\bz)
>\val(c_\bi)+\bi \val(\bz)\ ,\nonumber\\ &\quad\bi\in I\backslash\{\bi_0,...,\bi_r\}\
,\nonumber\end{eqnarray} which means that $\bx=\val(\bz)$ belongs
to the corner locus of the graph of $N_f$.

Suppose now that $\bx=(s_1,...,s_k)\in\R^k$ satisfies
\begin{equation}\bi_1\bx+\val(c_{\bi_1})=...=\bi_r\bx+\val(c_{\bi_r})\stackrel{\text{def}}{=}\rho>\bi\bx+\val(c_\bi)
\label{e22}\end{equation} for some $r\ge 2$ and all $\bi\in
I\backslash\{\bi_1,...,\bi_r\}$. Assume that the first coordinate
of $\bi_1,...,\bi_r$ takes values $m_1,...,m_p$, $p\ge 2$. Choose
$z^0_2,...,z^0_k\in\K^*$ so that \mbox{$\val(z^0_i)=s_i$},
$i=2,...,k$, and the coefficients $b_{m_1}$, ..., $b_{m_p}$ of
$z_1^{m_1}$, ..., $z_1^{m_p}$, respectively, in the polynomial
$\varphi(z_1)=f(z_1,z^0_2,...,z^0_k)$ satisfy
$$\val(b_{m_1})=\rho-m_1s_1,\quad ...\quad,\ \val(b_{m_p})=\rho-m_ps_1\
.$$ Notice that in view of (\ref{e22}), for any other coefficient
$c_m$ of $z_1^m$ in $\varphi(z_1)$, $m\ne m_1,...,m_p$, it holds
$\val(b_m)<\rho-ms_1$, i.e., in the Newton diagram of $\varphi$
spanned by points $(i,\val(b_i))$, the monomials of degrees
$m_1,...,m_p$ form an edge, and thus, there is a root $z^0_1$ of
$\varphi$ with $\val(z^0_1)=s_1$. \proofend

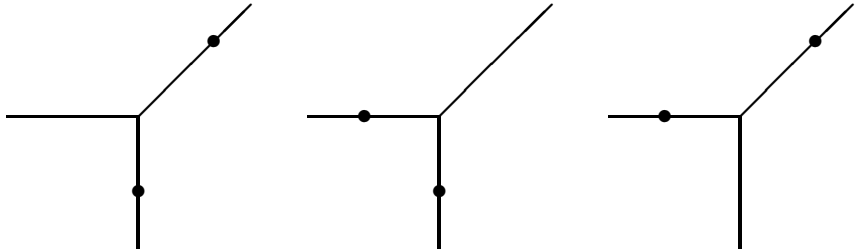
\begin{figure}
\setlength{\unitlength}{1cm}
\begin{picture}(13,6)(0,0)
\thicklines \put(2.5,2.5){\line(-1,0){1.75}}
\put(2.5,2.5){\line(0,-1){1.75}} \put(2.5,2.5){\line(1,1){1.5}}
\put(6.5,2.5){\line(-1,0){1.75}} \put(6.5,2.5){\line(0,-1){1.75}}
\put(6.5,2.5){\line(1,1){1.5}} \put(10.5,2.5){\line(-1,0){1.75}}
\put(10.5,2.5){\line(0,-1){1.75}} \put(10.5,2.5){\line(1,1){1.5}}
\put(2.4,1.4){$\bullet$} \put(6.4,1.4){$\bullet$}
\put(3.4,3.4){$\bullet$} \put(5.4,2.4){$\bullet$}
\put(9.4,2.4){$\bullet$} \put(11.4,3.4){$\bullet$}
\end{picture}
\caption{Plane amoebas of the first order}\label{f1}
\end{figure}

Non-Archimedean amoebas unexpectedly reveal many common properties
with algebraic varieties. For example (see \cite{M}), there is one
and only one amoeba of a straight line through two generic points
in the plane (see Figure \ref{f1}). Similarly there exists one and
only one amoeba of a conic curve through five generic points in
the plane. To introduce the reader to the subject, we extend this
existence and uniqueness result to amoebas of polynomials with
arbitrary support and in any number of variables.

\begin{theorem}\label{t4}
Given arbitrary integers $k\ge 2$, $n\ge 1$ and a finite set
$I\subset\Z^k$ consisting of $n+1$ points, for a generic $n$-tuple
$(\bx_1,...,\bx_n)\in(\R^k)^n$, there exists one and only one
amoeba $A\in{\cal A}(I)$ passing through $\bx_1,...,\bx_n$.
\end{theorem}

{\bf Proof}. The existence part is trivial: just take the amoeba
of a hypersurface $Z_f$, $f\in F_\K(I)$, passing through any
$n$-tuple $(w_1,...,w_n)\in
\val^{-1}(\bx_1,...,\bx_n)\subset((\K^*)^k)^n$.

To prove the uniqueness, we impose the following condition on
$\bx_1,...,\bx_n$. Consider the $n\times(n+1)$ matrix $M$, whose
$i$-th row entries are $\bx_i\bi$, $\bi\in I$. Assume that all the
sums of $n$ entries of $M$, taken one from each row and one from
each but one column, are distinct. This, clearly, excludes a
finite number of hyperplanes in $(\R^k)^n$. Given $n$ points
$w_1,...,w_n\in(\K^*)^k$, the coefficients $c_\bi$, $\bi\in I$ of
the polynomial $f\in F_\K(I)$ vanishing at these points can be
found as the $n\times n$ minors (with signs) of the $n\times(n+1)$
matrix $N$ whose $i$-th row entries are $w_i^\bi$, $\bi\in I$. If
$\val(w_1,...,w_n)=(\bx_1,...,\bx_n)$ then $\val(c_\bi)$ will be
the maximal sum of $n$ entries of $M$, taken one from each row and
one from each but the $\bi$-th column. Thus, $\val(c_\bi)$ does
not depend on the choice of $(w_1,...,w_n)$ in
$\val^{-1}(\bx_1,...,\bx_n)$, and hence by Theorem \ref{t3} all
such polynomials produce the same amoeba. \proofend

\subsection{Amoebas and subdivisions of Newton polytope}\label{sec7}
For a polynomial $f\in F_\K(I)$, one can define a subdivision of
the Newton polytope \mbox{$\Del=\conv(I)$} into convex polytopes
with vertices from $I$. Namely, take the convex hull $\Del_v(F)$
of the set $\{(\bi,-\val(c_\bi))\in\R^{k+1}\ :\ \bi\in I\}$ and
define the function
$$\nu_f:\Del\to\R,\quad\nu_f(\bi)=\min\{x\ :\
(\bi,x)\in\Del_v(f)\}\ .$$ This is a convex piece-wise linear
function, whose linearity domains are convex polytopes with
vertices in $I$, which form a subdivision $S_f$ of $\Del$. It is
easy to see (for example, from the fact that the functions $N_f$
and $\nu_f$ are dual by the Legendre transform) that

\begin{lemma}\label{l8}
The subdivision $S_f$ of $\Del$ is combinatorially dual to the
pair $(\R^k,A_f)$.
\end{lemma}

Notice that, in general, the geometry of an amoeba $A\in{\cal
A}(\Del)$ determines a dual subdivision $S$ of $\Del$ not
uniquely, but up to a combinatorial isotopy, in which all edges
remain orthogonal to the corresponding edges of $A$, and vice
versa. Combinatorially isotopic amoebas form a subset\footnote{It
is, in fact, the interior of a convex polyhedron in ${\cal
A}(I)$.} in ${\cal A}(I)$, whose dimension we call the {\it rank
of amoeba} (or the {\it rank of subdivision}) and denote
$\rk(A_f)=\rk(S_f)$.

\begin{lemma}\label{l12} For the case $k=2$, and $S_f:\ \Del=\Del_1\cup...\cup\Del_N$,
\begin{equation}\rk(S_f)\ge\rk_{\text{\rm exp}}(S_f)\stackrel{\text{\rm def}}{=}|V(S_f)|-1-\sum_{i=1}^N(|V(\Del_i)|-3)\
,\label{e32}\end{equation} where $V(S_f)$ is the set of vertices
of $S_f$, $V(\Del_i)$ is the set of vertices of the polygon
$\Del_i$, $i=1,...,N$. More precisely,
\begin{equation}\rk(S_f)=\rk_{\text{\rm exp}}(S_f)+d(S_f)\ ,\label{e40}
\end{equation}
where \begin{itemize}\item $d(S_f)=0$ if all the polygons
$\Del_1,...,\Del_N$ are triangles or parallelograms,
\item otherwise,
\begin{equation}0\le 2d(S_f)\le\sum_{m\ge
2}((2m-3)N_{2m}-N'_{2m})+\sum_{m\ge 2}(2m-2)N_{2m+1}-1\ ,
\label{e41}\end{equation} where $N_m$, $m\ge 3$, is the number of
$m$-gons in $S_f$, $N'_{2m}$ is the number of $2m$-gons in $S_f$,
whose opposite edges are parallel, $m\ge 2$.
\end{itemize}
\end{lemma}

{\bf Proof}. Inequality (\ref{e32}) is obvious, since an
$m$-valent vertex of $A_f$ imposes $m-3$ linear conditions on the
planes forming the graph of $N_f$.

Assume that all $\Del_1,...,\Del_N$ are triangles or
parallelograms, and show that the conditions imposed by the
$4$-valent vertices of $A_f$ are independent. Take a vector
$\overline a\in\R^2$ with an irrational slope and coorient each
edge of any parallelogram so that the normal vector forms an acute
angle with $\overline a$. This coorientation defies a partial
ordering on the set of parallelograms, which we complete somehow
up to a linear ordering. Notice that each parallelogram has two
neighboring edges cooriented outside. Altogether this means that
the coefficients of the linear conditions imposed by the
$4$-valent vertices of $A_f$ can be arranged into a triangular
matrix, and hence are independent, i.e., $d(S_f)=0$.

If $S_f$ contains polygons, different from triangles and
parallelograms, we define a linear ordering on the set of all
non-triangles in the same manner as above. Denote by $e_-(\Del_i)$
(resp., $e_+(\Del_i)$) the number of edges of a polygon $\Del_i$
cooriented outside (resp., inside) $\Del_i$. Passing inductively
over non-triangular polygons $\Del_i$, each time we add at least
\mbox{$\min\{e_-(\Del_i)-1,\ |V(\Del_i)|-3\}$} new linear
conditions independent of all the preceding ones. Thus,
$$d(S_f)\le\sum_{i=2}^N\left(|V(\Del_i)|-3-\min\{e_-(\Del_i)-1,\ |V(\Del_i)|-3\}\right)$$
$$=\sum_{i=2}^N\max\{|V(\Del_i)|-e_-(\Del_i)-2,\ 0\}\ ,$$ since, for the initial
polygon $\Del_1$, all $|V(\Del_1)|-3$ imposed conditions are
independent. Replacing $\overline a$ by $-\overline a$, we obtain
$$d(S_f)\le\sum_{i=1}^{N-1}\max\{|V(\Del_i)|-e_+(\Del_i)-2,\ 0\}\
.$$ Since \begin{itemize}\item $1\le e_-(\Del_i)\le|V(\Del_i)|-1$
and $e_-(\Del_i)+e_+(\Del_i)=|V(\Del_i)|$ yield
$$\max\{|V(\Del_i)|-e_-(\Del_i)-2,\ 0\}+\max\{|V(\Del_i)|-e_+(\Del_i)-2,\
0\}\le|V(S_f)|-3\ ,$$ \item for a $2m$-gon with parallel opposite
edges, $$e_-=e_+=m\quad\Longrightarrow\quad\max\{2m-e_--2,\
0\}+\max\{2m-e_+-2,\ 0\}=2m-4\ ,$$
\end{itemize}
we get \begin{equation}2d(S_f)\le\sum_{m\ge
2}((2m-3)N_{2m}-N'_{2m})+\sum_{m\ge 2}(2m-2)N_{2m+1}\
.\label{e43}\end{equation} If among $\Del_1,...,\Del_N$ there is a
polygon $\Del_i$ with an odd $\ge 5$ number of edges, or a polygon
with an even number of edges and a pair of non-parallel opposite
sides, then $\overline a$ can be chosen so that
$\min\{e_-(\Del_i),e_+(\Del_i)\}\ge 2$, and thus, the contribution
of $\Del_i$ in the latter bound for $2d(S_f)$ will be
$|V(\Del_i)|-4$, which allows us to gain $-1$ on the right-hand
side of (\ref{e43}) and obtain (\ref{e41}).

Lastly, assume that all non-triangular polygons in $S_f$ have an
even number of edges and their opposite sides are parallel, and
furthermore, that there is $\Del_i$ with $|V(\Del_i)|=2m\ge 6$.
The union of the finite length edges of $A_f$ is the adjacency
graph of $\Del_1,...,\Del_N$. Take the vertex corresponding to
$\Del_i$, pick a generic point $O$ in a small neighborhood of this
vertex, and orient each finite length edge of $A_f$ so that it
forms an acute angle with the radius-vector from $O$ to the middle
point of the chosen edge. The adjacency graph, equipped with such
an orientation, has no oriented cycles, since the terminal point
of any edge is further from $O$ than the initial one. Thus, we
obtain a partial ordering on $\Del_1,...,\Del_N$ such that, for
any $\Del_k$ with an even number of edges, at least half of them
is cooriented outside. Then we apply the preceding argument to
estimate $d(S_f)$ and notice that the contribution of the initial
polygon $\Del_i$ to such a bound is zero, whereas on the
right-hand side of (\ref{e43}) it is at least two, and this
completes the proof of (\ref{e41}). \proofend

\subsection{Algebraic curves over $\K$ and $\C$: general fibre and
tropicalization}\label{sec111} Let $\Del\subset\R^2$ be a
non-degenerate convex lattice polygon, $C\in\Lam_\K(\Del)$ a curve
with only isolated singularities, which is defined by a polynomial
$f(x,y)$ as in (\ref{e50}). This curve gives rise to some complex
algebraic curves.

First, evaluating the coefficients of $f(x,y)$ at small positive
$t$ (or at complex non-zero $t$ close to zero, if the exponents of
$t$ in the coefficients $a_{ij}(t)$ of $f(x,y)$ are integral), we
obtain a family of curves $C^{(t)}\in\Lam(\Del)$. The relation
between $C$ and $C^{(t)}$ is formulated in the following
statement, in which we understand a topological type of isolated
singular point (over any algebraically closed field of
characteristic zero) as a minimal resolution tree with given
multiplicities of the point itself and of its infinitely near
points, or, equivalently, the number of local branches, their
characteristic Puiseux exponents and pair-wise intersection
multiplicities.

\begin{lemma}\label{l9} The collection of topological types of
singular points of a reduced curve $C\in\Lam_\K(\Del)$ coincides
with the collection of topological types of singular points of a
generic curve $C^{(t)}\in\Lam(\Del)$. A curve $C\in\Lam_\K(\Del)$
is reducible if and only if a generic curve $C^{(t)}\in\Lam(\Del)$
is reducible.
\end{lemma}

This immediately follows from the fact that the set of curves in a
given linear system, having singularities of prescribed
topological types, over any algebraically closed field of
characteristic zero is defined by the same system of polynomial
equalities and inequalities with integer coefficients. The same
argument confirms the simultaneous reducibility of $C$ and
$C^{(t)}$.

Notice that, shrinking the range of $t$ if necessary, we obtain
that the curve $C$ bears a one-parametric equisingular deformation
of complex curves.

\medskip

We shall also define certain limits of $C^{(t)}$ as $t\to 0$.
Namely, let $\nu_f:\Del\to\R$ be a convex function, $S_f$ the
corresponding subdivision $\Del=\Del_1\cup...\cup\Del_N$, as
defined in the preceding section. The restriction
$\nu_f\big|_{\Del_i}$ coincides with a linear (affine) function
$\lam_i:\Del\to\R$, $\lam_i(\bx)=\bi_i\bx+\gam_i$,
$\bi_i=(\alp_1,\alp_2)\in\R^2$, $\gam_i\in\R$, $i=1,...,N$. Then
the polynomial
\begin{equation}t^{-\gam_i}f(z_1t^{-\alp_1},z_2t^{-\alp_2})=\sum_{\bi\in
\Del\cap\Z^2}\wc_\bi\bz^\bi\label{e25}\end{equation} satisfies the
following condition:
$$\val(\wc_\bi)\begin{cases}=0,\quad&\text{if}\ \bi\ \text{is a
vertex of}\ \Del_i,\\ \le 0,\quad&\text{if}\ \bi\in\Del_i,\\
<0,\quad&\text{if}\ \bi\not\in\Del_i\ .\end{cases}$$ In other
words, letting $t=0$ on the right-hand side of (\ref{e25}), we
obtain a {\it complex} polynomial $f_i$ with Newton polygon
$\Del_i$, which in turn define complex curves
$C_i\in\Lam(\Del_i)$, $i=1,...,N$. Notice that multiplying
$f(x,y)$ by a constant from $\K^*$ does not change $S_f$ and
$C_1,...,C_N$, but adds a linear function to $\nu_f$. The
collection $(\nu_f,S_f;C_1,...,C_N)$ is called the {\it
tropicalization} (or {\it dequantization}) of the curve $C$, and
denoted by ${\cal T}(C)$. We also call $f_i$ and $C_i$ the
tropicalizations of the polynomial $f$ and the curve $C$ on the
polygon $\Del_i$, $1\le i\le N$.

Assume that the exponents of $t$ in the coefficients $a_{ij}(t)$
of $f(x,y)$ are rational. By a change of parameter $t\mapsto t^m$,
we can make all these exponents integral and the function $\nu_f$
integral-valued at integral points. Introduce the polyhedron
$$\widetilde\Del=\{(\alp,\bet,\gam)\in\R^3\ :\
(\alp,\bet)\in\Del,\ \gam\ge\nu_f(\alp,\bet)\}\ .$$ It defines a
toric variety $Y=\Tor(\widetilde\Del)$, which naturally fibers
over $\C$ so that the fibres $Y_t$ over $t\ne 0$ are isomorphic to
$\Tor(\Del)$, and $Y_0$ is the union of toric surfaces
$\Tor(\widetilde\Del_i)$, $i=1,...,N$, with
$\widetilde\Del_1,...,\widetilde\Del_N$ being the faces of the
graph of $\nu_f$. By the choice of $\nu_f$,
$\Tor(\widetilde\Del_i)\simeq\Tor(\Del_i)$, and we shall simply
write that $Y_0=\bigcup_i\Tor(\Del_i)$. Then the curve $C$ can be
interpreted as an analytic surface in a neighborhood of $Y_0$,
which fibers into the complex curves $C^{(t)}\subset
Y_t\simeq\Tor(\Del)$, and whose closure intersects $Y_0$ along the
curve $C^{(0)}$ that can be identified with
$\bigcup_iC_i\subset\bigcup_i\Tor(\Del_i)$. Passing if necessary
to a finite cyclic covering ramified along $Y_0$, we can make
$\Tor(\widetilde\Del)$ non-singular everywhere but may be at
finitely many points, corresponding to the vertices of
$\widetilde\Del$, and, in addition, make the surfaces
$\Tor(\Del_k)\backslash\Sing(\Tor(\widetilde\Del))$, $k=1,...,N$,
smooth and intersecting transversally in
$\Tor(\widetilde\Del)\backslash\Sing(\Tor(\widetilde\Del))$.

The singular points of the curves $C^{(t)}$ define sections
$s:D\backslash\{0\}\to\Tor(\widetilde\Del)$, $D\subset\C$ being a
small disc centered at $0$. The limit points $z=\lim_{t\to 0}s(t)$
are singular points of $C^{(0)}$. We say that such a point $z\in
C^{(0)}$ bears the corresponding singular points of $C^{(t)}$. If
$z\in C^{(0)}$ does not belong to the intersection lines
$\bigcup_{i\ne j}\Tor(\Del_i\cap\Del_j)$ and bears just one
singular point of $C^{(t)}$, which is topologically equivalent to
$z$, we call $z$ a regular singular point, otherwise it is
irregular. If $C^{(0)}$ has irregular singular points, we can
define a {\it refinement of the tropicalization} in the following
way: transform the polynomial $f(x,y)$ into $f(x+a,y+b)$ with
$a,b\in\K$ such that the irregular singular point of $C^{(0)}$
goes to the origin, and consider the tropicalization of the curve
defined by the new polynomial $f(x+a,y+b)$. This provides
additional information on the behavior of singular points of
$C^{(t)}$ tending to irregular singular points of $C^{(0)}$, and
corresponds, in a sense, to blowing-up the threefold $Y$ at the
irregular singular points of $C^{(0)}$ (cf. \cite{ST}).

\section{Counting nodal curves}\label{sec302}

\subsection{Formulation of the result}\label{sec5}
Let $\Del\subset\R^2$ be a non-degenerate lattice polygon, which
has integral points in its interior. It is well known that the
number of nodes of an irreducible curve in $\Lam_\K(\Del)$ does
not exceed $|\Int(\Del)\cap\Z^2|$. For any positive integer
$n\le|\Int(\Del)\cap\Z^2|$, denote by $\Sig_\Del(nA_1)$ the set of
reduced curves in $\Lam_\K(\Del)$ having exactly $n$ nodes as
their only singularities and defined by polynomials with Newton
polygon $\Del$. This is a smooth quasiprojective subvariety of
$\Lam_\K(\Del)$ (so-called {\it Severi variety}) of codimension
$n$, i.e., $\dim\Sig_\Del(nA_1)=r=|\Del\cap\Z^2|-1-n$ in view of
$\dim\Lam_\K(\Del)=|\Del\cap\Z^2|-1$. Imposing the condition to
pass through $r$ generic points in $(\K^*)^2\subset\Tor_\K(\Del)$,
we obtain a finite set of curves in $\Sig_\Del(nA_1)$, whose
cardinality is just $\deg\Sig_\Del(nA_1)$.

Now we describe amoebas which are projections of nodal curves,
passing through generic points in $(\K^*)^2$ with distinct
valuation projections to $\R^2$. An amoeba $A\in{\cal A}(\Del)$ is
called {\it nodal}, if its dual subdivision $S$ of $\Del$ is as
follows: \begin{itemize}
\item all the points in $\partial\Del\cap\Z^2$ are vertices of
$S$, \item $S$ consists of triangles and parallelograms.
\end{itemize}
Define the weight of a nodal amoeba $A$ by
$$W(A)=\prod_{\renewcommand{\arraystretch}{0.6}
\begin{array}{c}
\scriptstyle{\Del'\in P(S)}\\
\scriptstyle{|V(\Del')|=3}
\end{array}}|\Del'|\ ,$$ where $P(S)$ denotes the set of polygons
of $S$, $|\Del'|$ stands for the double Euclidean area of $\Del'$.

\begin{theorem}\label{t5} In the previous notation,
$$\deg\Sig_\Del(nA_1)=\sum_AW(A)\ ,$$ where the sum ranges over all
nodal amoebas of rank $\ r$, passing through $\ r$ fixed generic
points in $\R^2$.
\end{theorem}

\begin{remark}\label{r1}
{\rm Our formula coincides with that given by Mikhalkin \cite{M2}.
Namely, the multiplicity of a lattice path defined in \cite{M2} is
just the sum of multiplicities of nodal amoebas which correspond
to subdivisions of $\Del$, arising from the given path along the
construction of \cite{M2}. We notice also that the generality
requirement for the position of $r$ points in $\R^2$ will be
specified in the proof, and one can easily check that the
configurations considered by Mikhalkin \cite{M3}, i.e., generic
points on a generic straight line, satisfy these generality
conditions.}
\end{remark}

The proof comprises three main steps. First, we determine amoebas
and tropicalizations of nodal curves in the count, in particular,
that the amoebas are nodal of rank $r$ (section \ref{sec9}). Then
we refine tropicalizations in a suitable way (sections
\ref{sec10}, \ref{sec12}). Finally, using the patchworking
theorem, we show that the refined tropicalization gives rise to an
explicit number of nodal curves passing through given points
(section \ref{sec11}).

\subsection{Deformation of reducible surfaces and curves}
We start with the following auxiliary statement.

\begin{lemma}\label{l6} Let a complex threefold $Y$ be smooth at a
point $z$, $U\subset Y$ a small ball centered at $z$. Assume that
$\pi:U\to(\C,0)$ is a flat family of reduced surfaces such
$U_0=\pi^{-1}(0)$ consists of two smooth components $U'_0,U''_0$
which intersect transversally along a line $L\supset\{z\}$, and
$U_t=\pi^{-1}(t)$ are nonsingular as $t\ne 0$. Let $C'_0\subset
U'_0$, $C''_0\subset U''_0$ be reduced algebraic curves, which
cross $L$ only at $z$ and with the same multiplicity $m\ge 2$.
Assume also that $U'_0$, $U''_0$ are regular neighborhoods for the
(possibly singular) point $z$ of $C'_0$ and $C''_0$, respectively.
Let $\del'=\del(C'_0,z)$, $\del''=\del(C''_0,z)$ be the
$\del$-invariants, $r'$, $r''$ the numbers of local branches of
$C'_0$, $C''_0$ at $z$, respectively. Then in any flat deformation
$C_t$, $t\in(\C,0)$, of $C_0=C'_0\cup C''_0$ such that $C_t\subset
U_t$, the total $\del$-invariant of $C_t$, $t\ne 0$, in $U_t$ does
not exceed
$$\del'+\del''+m-\max\{r',\ r''\}\ .$$
\end{lemma}

{\bf Proof}. Topologically, the curves $C'_0$ and $C''_0$ (in $U$)
are bouquets of $r'$ and $r''$ discs, respectively. Notice that
the circles of $C'_0\cap\partial U$ and $C''_0\cap\partial U$ move
slightly when $t$ changes, and they are not contractible in $U_t$
for $t\ne 0$. For instance, a circle of $C'_0\cap\partial U$ is
(positively) linked with the line $L$ in $U'_0$, and hence remains
(positively) linked with the surface $U''_0$ in $U$; thus, it
cannot be contracted in $U_t$, $t\ne 0$, which does not intersect
$U''_0$. This means that the curve $C_t\subset U_t$, $t\ne 0$, is
the union of a few immersed surfaces with a total of $r'+r''$
holes and at least $\max\{r',r''\}$ handles.

Now the asserted upper bound can be derived either from a local
count of intersections and self-intersections of the components of
$C_t$, or by a ``global" reasoning. For the latter, we consider
the following model situation, which is quite relevant to our
consideration and is explored in more detail below in the proof of
Theorem \ref{t5}. Namely, assume that
\begin{itemize}\item $p\gg m$, $p\in\N$, \item $\Del\subset\R^2$ is the triangle with
vertices $(0,0)$, $(2p,0)$, $(0,2p)$,
\item $\nu:\Del\to\R$ is the function such that $\nu(\alp,\bet)=0$ as $\alp+\bet\le p$, and
$\nu(\alp,\bet)=\alp+\bet-p$ as $p\le\alp+\bet$, \item
$\widetilde\Del=\{(\alp,\bet,\gam)\in\R^3\ :\ (\alp,\bet)\in\Del,\
\nu(\alp,\bet)\le\gam\le p+1\}$.\end{itemize} Then
$Y=\Tor(\widetilde\Del)$ is a non-singular threefold,
$Y'_0=\Tor(\Del')$, $Y''_0=\Tor(\Del'')$ are surfaces isomorphic
respectively to $\PP^2$ and $\PP^2$ with a blown up point, where
\mbox{$\Del'=\conv\{(0,0,0),(p,0,0),(0,p,0)\}$},
\mbox{$\Del''=\conv\{(p,0,0),(0,p,0),(2p,0,p),(0,2p,p)\}$} are the
faces of $\Graph(\nu)$. Furthermore, $Y'_0$ and $Y''_0$ intersect
transversally along the line $L=\Tor(\sig)$,
$\sig=[(p,0,0),(0,p,0)]$. A neighborhood $V$ of $Y_0=Y'_0\cup
Y''_0$ in $Y$ admits a fibration $V\to(\C,0)$ with the zero fibre
$Y_0$ and other fibres $Y_t$ being the closures of the images of
the hyperplanes $\{x_3=t\}\subset(\C^*)^3$ by the standard
embedding of $(\C^*)^3$ into $Y$ with the coordinate
correspondence $(\alp,\bet,\gam)\leftrightarrow(x_1,x_2,x_3)$ of
$\R^3$ and $(\C^*)^3$. Clearly, $Y_t\simeq\PP^2$, $t\ne 0$. Assume
that the curves $C'_0\subset Y'_0$, $C''_0\subset Y''_0$ are given
by polynomials with Newton polygons $\Del',\Del''$, respectively,
with a common truncation to $\sig$, and such that they have a
common point $z\in L$ as in the statement of the lemma, are
non-singular outside $z$, and intersect $L$ transversally outside
$z$ (at common points). The flatness of a deformation $C_t\subset
Y_t$, $t\in(\C,0)$, of the curve $C_0=C'_0\cup C''_0$ means that
$C_t$, $t\ne 0$, tends to a curve of degree $2p$ by the
isomorphism $Y_t\simeq\PP^2$. Denoting by $U$ a neighborhood of
$C'_0\cap C''_0$ in $Y$, we obtain for $\check\chi(C_t)$, the
Euler characteristic of the normalization of $C_t$, the following
bound
$$\check\chi(C_t)=\chi(C_t\backslash U)+\check\chi(C_t\cap
U)=\chi(C'_0\backslash U)+\chi(C''_0\backslash
U)+\check\chi(C_t\cap U)$$
$$\le(-p^2+2p+m-r'+2\del')+(-3p^2+4p+m-r''+2\del'')+(r'+r''-2\max\{r',r''\})$$
$$=-4p^2+6p+2m+2\del'+2\del''-\max\{r',r''\}\ ,$$ and hence for
the total $\del$-invariant of $C_t$
$$\del(C_t)=\frac{(2p-1)(2p-2)}{2}-g(C_t)=\frac{(2p-1)(2p-2)}{2}-1+\frac{\check\chi(C_t)}{2}$$
$$\le\del'+\del''+m-\max\{r',r''\}\qquad\qquad\mbox{\proofend}$$

\begin{example}\label{ex1}
{\rm In the notation of Lemma \ref{l6}, if $C'_0$, $C''_0$ are
non-singular at $z$, then $\del_1=\del_2=0$, $r_1=r_2=1$, and the
number of nodes in a deformation does not exceed $m-1$, and this
number can be attained \cite{Ch}, Theorem 2.1. We can produce the
maximal number of nodes by means of suitable deformation patterns
(i.e., certain affine curves) as defined below in section
\ref{sec10}.}
\end{example}

\begin{remark}\label{rnew1}
In the notations of Lemma \ref{l6}, accept all the hypotheses but
assume that $C'_0,C''_0$ are not necessarily reduced. Furthermore,
let $C'_0$ (resp., $C''_0$) have $r'$ (resp., $r''$) reduced local
branches at $z$ of multiplicities $\rho'_1,...,\rho'_{r'}$ (resp.,
$\rho''_1,...,\rho''_{r''}$). Then the argument used in the proof
of Lemma \ref{l6} shows that, if $C_t$ is reduced in $U$, then
$$\check\chi(C_t\cap U)\le-\min|m'_1+...+m'_{r'}-m''_1-...-m''_{r''}|\
,$$ where integers $m'_1,...,m''_{r''}$ run over the range $1\le
m'_1\le \rho'_1$, ..., $1\le m''_{r''}\le \rho''_{r''}$.
\end{remark}

\subsection{Amoebas and tropicalizations of nodal curves passing through
generic points}\label{sec9} Let $\bx_1,...,\bx_r\in\R^2$ be
generic distinct points with rational coordinates, and let
$\bp_1,...,\bp_r\in(\K^*)^2$ be generic points satisfying
$\val(\bp_i)=\bx_i$, $i=1,...,r$, and having only rational
exponents of the parameter $t$.

Observe that the coefficients of a polynomial $f\in\K[x,y]$, which
defines a curve $C\in\Sig_\Del(nA_1)$, are Puiseux series with
rational exponents of $t$. A parameter change $t\mapsto t^m$ with
a suitable natural $m$ makes all these exponents integral, and the
convex piece-wise linear function $\nu_f:\Del\to\R$
integral-valued at integral points. We keep these assumptions
through out the rest of the paper.

Let $S_f:\ \Del=\Del_1\cup...\cup\Del_N$ be the subdivision
defined by $\nu_f$, $(C_1,...,C_N)$ the tropicalization of the
curve $C=\{f=0\}\in\Lam_K(\Del)$. The union of the divisors
$\Tor(\sig)\subset\Tor(\Del_k)$, where $\sig$ runs over all edges
of $\Del_k$, we shall denote by $\Tor(\partial\Del_k)$,
$k=1,...,N$. For any $i=1,...,N$, denote by $C_{ij}$,
$j=1,...,m_i$, the distinct irreducible components of the curve
$C_i\subset\Tor(\Del_i)$ and by $r_{ij}$, $j=1,...,m_i$, their
multiplicities. Denote by $s_{ij}$ the number of local branches of
$C_{ij}$ centered on $\partial\Del_i))$, $j=1,...,m_i$.

We intend to estimate $\check{\chi}(C^{(t)})$ from above and from
below and to compare these bounds.

Let $U$ be the union of small open balls $U_z$ in the three-fold
$Y$ (see the definition in section \ref{sec111}) centered at all
the points $z\in\bigcup_i(C_i\cap\Tor(\partial\Del_i))$. If $z\in
C_i\cap\Tor(\sig)$, where $\sig$ is an edge of $\Del_i$ lying on
$\partial\Del$, then $\check{\chi}(C^{(t)}\cap U_z)$ does not
exceed the number of local branches of $C_i$ at the points of
$C_i\cap\Tor(\sig)$. If $z\in\Tor(\sig)\cap C_i\cap C_k$, where
$\sig=\Del_i\cap\Del_k$ is a common edge, then
$\check{\chi}(C^{(t)}\cap U_z)\le 0$ by Remark \ref{rnew1}. Hence
\begin{equation}\check{\chi}(C^{(t)}\cap
U)\le|\partial\Del\cap\Z^2|\label{e31}\end{equation} with an
equality if and only if, for any edge
$\sig\subset\Del_i\cap\partial\Del$, the reduction of the curve
$C_i$ is non-singular along $\Tor(\sig)$ and meets $\Tor(\sig)$
transversally.

For the upper bound to $\check{\chi}(C^{(t)}$, we can assume that,
for any $i=1,...,N$, and $1\le j<j'\le m_i$, the components
$C_{ij}$ and $C_{ij'}$ do not glue up in $Y\backslash U$ when
$C^{(0)}$ deforms into $C^{(t)}$. Then the normalization of
$C^{(t)}\backslash U$ is the union of connected components, each
of them tending to some curve $C_{ij}\backslash U$. Furthermore,
the components which tend to a certain $C_{ij}\backslash U$ can be
naturally projected onto $C_{ij}\backslash U$, and this projection
is an $r_{ij}$-sheeted covering (possibly ramified at a finite
set). Hence
$$\check{\chi}(C^{(t)}\backslash
U)\le\sum_{i=1}^N\sum_{j=1}^{m_i}r_{ij}\check{\chi}(C_{ij}\backslash
U)=\sum_{i=1}^N\sum_{j=1}^{m_i}r_{ij}(\check{\chi}(C_{ij})-s_{ij})$$
$$\le 2\sum_{i=1}^Nm_i-\sum_{i=1}^N\sum_{j=1}^{m_i}s_{ij}$$ with
an equality only if all $C_{ij}$ are rational, and $r_{ij}=1$ as
far as $s_{ij}>2$. Next we notice that $s_{ij}\ge 2$ for any
$C_{ij}$, and $s_{ij}\ge 3$ for at least one of the components
$C_{ij}$ if $\Del_i$ has an odd number of edges, or $\Del_i$ has
an even number of edges, but not all pairs of opposite sides are
parallel. Hence (in the notation of Lemma \ref{l12})
\begin{equation}\check{\chi}(C^{(t)}\backslash
U)\le-N_3-\sum_{j\ge
2}(N_{2j+1}+N_{2j}-N'_{2j})\label{e30}\end{equation} with an
equality only if, for each triangle $\Del_i$, $C_i$ is irreducible
and satisfies $s_{ij}=3$; for each $\Del_i$ with an odd $\ge 5$
number of edges or with an even number of edges, but not all pairs
of opposite sides parallel, exactly one component $C_{ij}$
satisfies $s_{ij}=3$ and the others satisfy $s_{ij}=2$; and,
finally, $s_{ij}=2$ for all components $C_{ij}$ in the remaining
polygons $\Del_i$. Notice also that $s_{ij}=2$ means that $C_{ij}$
is defined by a binomial.

On the other hand,
$$\check{\chi}(C^{(t)})=2-2g(C^{(t)})=2-2(|\Int(\Del)\cap\Z^2|-n)$$
$$=2-2|\Int(\Del)\cap\Z^2|
+2(|\Del\cap\Z^2|-1-r)=2|\partial\Del\cap\Z^2|-2r$$
$$\ge 2|\partial\Del\cap\Z^2|-2\cdot\rk(S_f)=
2|\partial\Del\cap\Z^2|-2\cdot\rk_{\text{\rm exp}}(S_f)-2d(S_f)$$
with an equality only if $\rk(S_f)=\rk_{\text{\rm
exp}}(S_f)+d(S_f)=r$. Next, by (\ref{e32}) we have
$$\check{\chi}(C^{(t)})\ge 2|\partial\Del\cap\Z^2|-2|V(S_f)|+2+2\sum_{i=1}^N(|V(\Del_i)|-3)-2d(S_f)$$
$$=2|\partial\Del\cap\Z^2|-2|V(S_f)|+2-2|V(S_f)\cap\partial\Del|+4|E(S_f)|-6N-2d(S_f)\
,$$ where $E(S_f)$ denotes the set of edges of $S_f$. Since
$|V(S_f)|-|E(S_f)|+N=1$, and
$2|E(S_f)|=3N_3+4N_4+5N_5+...+|V(S_f)\cap\partial\Del|$, we
finally obtain $$\check{\chi}(C^{(t)})\ge
2(|\partial\Del\cap\Z^2|-|V(S_f)\cap\partial\Del|)+|V(S_f)\cap\partial\Del|-N_3+N_5+2N_6+...-2d(S_f)\
.$$ Combining this with (\ref{e31}) and (\ref{e30}), we obtain
$$(|\partial\Del\cap\Z^2|-|V(S_f)\cap\partial\Del|)+\sum_{|V(\Del_i)|=3}(m_i-1)$$
$$+\sum_{m\ge 2}((2m-3)N_{2m}-N'_{2m})+\sum_{m\ge
2}(2m-2)N_{2m+1}\le 2d(S_f)\ ,$$ which, in view of Lemma
\ref{l12}, yields that each integral point on $\partial\Del$ is a
vertex of $S_f$ and all the non-triangular $\Del_i$ are
parallelograms.

Altogether the equality conditions for the upper and lower bounds
to $\check{\chi}(C^{(t)})$ prove that the amoeba $A_f$ is nodal of
rank $r$. Furthermore, \begin{itemize}\item for each triangle
$\Del_i$, the curve $C_i$ is rational and meets
$\Tor(\partial\Del_i)$ at exactly three points, where it is
unibranch;
\item for each parallelogram $\Del_i$, the polynomial, defining
$C_i$, is of type \mbox{$x^ky^l(\alp x^a+\bet y^b)^p(\gam x^c+\del
y^d)^q$} with $(a,b)=(c,d)=1$ and $(a:b)\ne(c:d)$.
\end{itemize}
We shall describe these curves more precisely.

\begin{lemma}\label{l7} For any lattice triangle $\Del'\subset\R^2$, there exists a
polynomial with Newton polygon $\Del'$ and prescribed coefficients
at the vertices of $\Del'$, which defines a rational curve
$C\subset\Tor(\Del')$, meeting $\Tor(\partial\Del')$ at exactly
three points, where it is unibranch. Furthermore, the curves
defined by these polynomials are nodal, and nonsingular at the
intersection with $\Tor(\partial\Del')$. Moreover, the number of
such polynomials is finite and equal to $|\Del'|$. An additional
fixation of one or two intersection points of $C$ with
$\Tor(\partial\Del')$ divides the number of polynomials under
consideration by the product of the length\footnote{We define the
length $|\sig|$ of a segment $\sig$ with integral endpoints as
$|\sig\cap\Z^2|-1$.} of the corresponding edges of $\Del'$.
\end{lemma}

{\bf Proof}. By a suitable lattice preserving transformation, we
can turn $\Del'$ into a triangle with vertices $(p,0)$, $(q,0)$,
$(0,m)$, $0\le p<q\le m$. Assuming that the curve $C$ crosses
$\Tor(\partial\Del)$ at points corresponding to the values $0$,
$1$ and $\infty$ of a uniformizing parameter $\tet$, we
necessarily obtain that $C$ is given by $x=\alp\tet^m$,
$y=\bet\tet^p(\tet-1)^{q-p}$. If the restrictions of the defining
polynomial on the edges $[(p,0),(q,0)]$ and $[(p,0),(0,m)]$ are
$x^p(x+\eps_1)^{q-p}$ and $(y^{m/d}+\eps_2x^{p/d})^d$,
respectively, where $d=\gcd(m,p)$, $\eps_1^{q-p}=\eps_2^d=1$, then
$$\alp+\eps_1=0,\quad \bet^{m/d}(-1)^{m(q-p)/d}+\eps_2\alp^{p/d}=0\
,$$ which gives $m(q-p)=|\Del|$ solutions for $(\alp,\bet)$.
Additional fixation of intersection points with
$\Tor(\partial\Del')$ means fixation of $\eps_1$ or/and $\eps_2$
and the respective reduction of the number of solutions.
Prescribed coefficients of $x^p,x^q,y^m$ in the polynomial can be
achieved by an appropriate coordinate change.

It remains to show that the curve $x=\tet^m$,
$y=\tet^p(\tet-1)^{q-p}$ is nodal. Since $\dot x(\tet)\ne 0$ as
$\tet\ne 0$, the curve has no local singular branches. Assuming
\begin{equation}\tet^m=\tet_1^m=\tet_2^m,\quad\tet^p(\tet-1)^{q-p}=\tet_1^p(\tet_1-1)^{q-p}=\tet_2^p(\tet_2-1)^{q-p}\
,\label{e24}\end{equation} we successively obtain
\begin{equation}
\begin{cases}&\tet_1=\tet\eps_1,\quad\tet_2=\tet\eps_2,\quad\eps_1^m=\eps_2^m=1,\quad\eps_1\ne\eps_2\
,\\ &
\tet-1=\eps_3(\tet\eps_1-1)=\eps_4(\tet\eps_2-1),\quad\eps_3^{q-p}=\eps_1^p,\
\eps_4^{q-p}=\eps_2^p\ ,\\
&\tet=\frac{1-\eta_1^p}{1-\eta_1^q}=\frac{1-\eta_2^p}{1-\eta_2^q}\end{cases}\
,\label{e23}\end{equation} where
$$\eps_1=\eta_1^{q-p},\quad\eps_3=\eta_1^p,\quad\eps_2=\eta_2^{q-p},\quad\eps_4=\eta_2^p,\quad
\eta_1^p\ne1,\ \eta_1^{q-p}\ne1,\ \eta_2^p\ne 1,\ \eta_2^{q-p}\ne
1\ .$$ Then, plugging $\eta_1=\cos\omega_1+\sqrt{-1}\sin\omega_1$,
$\eta_2=\cos\omega_2+\sqrt{-1}\sin\omega_2$ into (\ref{e23}), we
get
$$\frac{\cos(p\omega_1/2)\cos(q\omega_2/2)}{\cos(q\omega_1/2)\cos(p\omega_2/2)}=
\cos\frac{(q-p)(\omega_2-\omega_1)}{2}+\sqrt{-1}\sin\frac{(q-p)(\omega_2-\omega_1)}{2}\
,$$ and finally, $$(q-p)(\omega_2-\omega_1)\in
2\pi\cdot\Z\quad\Longrightarrow\quad\eps_1=\eta_1^{q-p}=\eta_2^{q-p}=\eps_2$$
in contrast to (\ref{e23}), and we are done. \proofend

\begin{lemma}\label{l11} Given integers $a,b,c,d$ such that
$(a,b)=(c,d)=1$, $(a:b)\ne(c:d)$, and any non-zero
$\alp,\bet,\gam,\del$, the curve $(\alp x^a+\bet y^b)(\gam
x^c+\del y^d)=0$ has $|\Del'\cap\Z^2|-3$ nodes as its only
singularities in $(\C^*)^2$, where $\Del'$ is the lattice
parallelogram built on the vectors $(a,-b)$, $(c,-d)$.
\end{lemma}

{\bf Proof}. Straightforward. \proofend

\subsection{Irreducible curves and irreducible
amoebas} One can speak of non-Archimedean amoebas as corner loci
of all possible convex piece-wise linear functions, whose generic
gradients are integral vectors (cf. \cite{St}). In such a sense,
an amoeba is called reducible if it is the union of two proper
sub-amoebas. The above description of the tropicalization of nodal
curves has the following interesting consequence.

\begin{lemma}\label{l30}
In the notation of section \ref{sec9}, the amoeba of an
irreducible (resp., reducible) $r$-nodal curve passing through
$\bp_1,...,\bp_r$ is irreducible (resp., reducible).
\end{lemma}

{\bf Proof}. The amoeba of a reducible curve is the union of the
amoebas of the irreducible components, and thus, is reducible.

Let a nodal amoeba $A$ of rank $r$ be the union of distinct
amoebas $A'$, $A''$. The intersection points of $A'$ and $A''$ are
four-valent vertices of $A$, and they correspond to some
parallelograms in the dual subdivision
$\Del=\Del_1\cup...\cup\Del_N$. As we have shown in section
\ref{sec9}, in the deformation $C^{(t)}$, $t\in(\C,0)$, for any
parallelogram $\Del_i$, the distinct components of the curve $C_i$
do not glue up. Hence the amoebas $A'$, $A''$ are lifted up to
separate algebraic curves. \proofend

\subsection{Refinement of the tropicalization at an isolated singular point}\label{sec10}
In the notation and hypotheses of the preceding section, we shall
introduce a refinement of the tropicalization for each point
$z\in\Tor(\sig)$, where $\sig=\Del_k\cap\Del_l$ is a common edge,
and the curves $C_k$, $C_l$ meet $\Tor(\sig)$ at $z$ with
multiplicity $m\ge 2$.

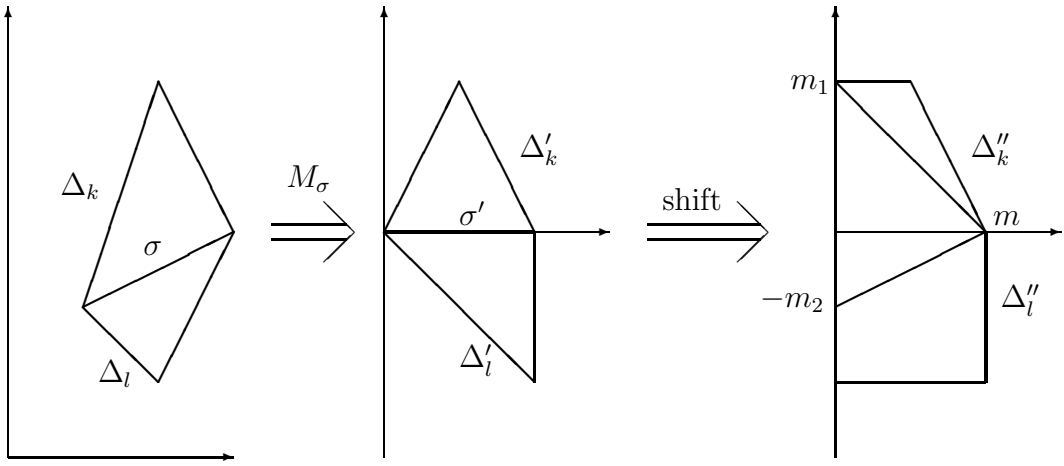
\begin{figure}
\setlength{\unitlength}{1cm}
\begin{picture}(13,6)(0,0)
\thinlines \put(0,0){\vector(1,0){3}}\put(0,0){\vector(0,1){6}}
\thicklines\put(1,2){\line(1,-1){1}}\put(2,1){\line(1,2){1}}
\put(1,2){\line(1,3){1}}\put(1,2){\line(2,1){2}}\put(2,5){\line(1,-2){1}}
\thinlines\put(3.5,2.9){\line(1,0){1}}\put(3.5,3.1){\line(1,0){1}}
\put(4.6,3.0){\line(-1,1){0.4}}\put(4.6,3.0){\line(-1,-1){0.4}}
\put(5,3){\vector(1,0){3}}\put(5,0){\vector(0,1){6}}
\thicklines\put(5,3){\line(1,0){2}}\put(5,3){\line(1,2){1}}
\put(5,3){\line(1,-1){2}}\put(6,5){\line(1,-2){1}}\put(7,1){\line(0,1){2}}
\thinlines\put(8.5,2.9){\line(1,0){1.5}}\put(8.5,3.1){\line(1,0){1.5}}
\put(10.1,3.0){\line(-1,1){0.4}}\put(10.1,3.0){\line(-1,-1){0.4}}
\put(11,3){\vector(1,0){3}}\put(11,0){\vector(0,1){6}}
\thicklines\put(11,5){\line(1,-1){2}}\put(11,5){\line(1,0){1}}
\put(12,5){\line(1,-2){1}}\put(11,2){\line(2,1){2}}\put(11,2){\line(0,-1){1}}
\put(11,1){\line(1,0){2}}\put(13,1){\line(0,1){2}}
\put(0.7,3.5){$\Del_k$}\put(1.2,1.0){$\Del_l$}\put(1.8,2.7){$\sig$}
\put(6.8,4){$\Del'_k$}\put(6,1.2){$\Del'_l$}\put(6,3.1){$\sig'$}
\put(12.8,4){$\Del''_k$}\put(13.2,2){$\Del''_l$}\put(13.1,3.1){$m$}
\put(10.4,4.9){$m_1$}\put(10,2){$-m_2$}\put(3.7,3.6){$M_\sig$}
\put(8.7,3.3){\text{\rm shift}}
\end{picture}
\caption{Refinement of the tropicalization, I}\label{f4}
\end{figure}

Let $\Del_k,\Del_l$ be triangles. Then $C_k$ and $C_l$ are
non-singular at $z$ and tangent to $\Tor(\sig)$ with multiplicity
$m$. To also cover the case of cuspidal curves, treated below in
section \ref{sec6}, we consider a more general situation. Namely,
assume that $C_k$ (resp., $C_l$) has at $z$ a semiquasihomogeneous
singularity topologically equivalent to $y^{m_1}+x^m=0$, $m_1\le
m$, (resp., $y^{m_2}+x^m=0$, $m_2\le m$). Perform the following
transformation of $f(x,y)$.
\begin{itemize}\item Let $M_\sig$ be an affine automorphism of
$\Z^2$ which takes $\Del$ in the right half-plane and takes $\sig$
into a segment $\sig'$ on the horizontal coordinate axis. This
corresponds to a monomial coordinate change $x=(x')^a(y')^b$,
$y=(x')^c(y')^d$ in $f(x,y)$ and further multiplication by a
monomial in $x',y'$. The truncation of the new polynomial
$f'(x',y')$ on the edge $\sig$ (i.e., the sum of the monomials of
$f'$ corresponding to the integral points in $\sig$) is a
polynomial in $x'$ over $\K$. Its tropicalization is a complex
polynomial in $x'$, the common truncation of the tropicalizations
$f'_k$ and $f'_l$ of $f'$ on the polygons $M_\sig(\Del_k)$,
$M_\sig(\Del_l)$. The point $z$ corresponds to a root $\xi\ne 0$
of $P_0(x')$. \item Without loss of generality, assume that
$\nu_f$ is zero along $\sig$ (just multiply $f(x,y)$ by a suitable
constant from $\K^*$). Then we perform the shift $x'=x''+\xi$,
$y'=y''$, and put $f''(x'',y'')=f'(x',y')$.
\end{itemize}
To understand the tropicalization of $f''$, apply the above
transformations to the polynomials
$$P_k(x,y)=\sum_{(i,j)\in\Del_k}c_{ij}^0t^{\lam_k(i,j)}x^iy^j,\quad
P_l(x,y)=\sum_{(i,j)\in\Del_l}c_{ij}^0t^{\lam_l(i,j)}x^iy^j\ ,$$
where
$$f_k(x,y)=\sum_{(i,j)\in\Del_k}c_{ij}^0x^iy^j,\quad
f_l(x,y)=\sum_{(i,j)\in\Del_l}c_{ij}^0x^iy^j$$ are the
tropicalizations of $f$ on $\Del_k,\Del_l$, and
$\lam_k=\nu_f\big|_{\Del_k}$, $\lam_l=\nu_f\big|_{\Del_l}$ are
linear functions. The Newton polygons $\Del''_k$, $\Del''_l$ of
the  resulting polynomials\footnote{From now on by ``polynomial"
we refer to a Laurent polynomial.} $P''_k$, $P''_l$ contain
segments $[(m,0),(0,m_1)]$, $[(m,0),(0,-m_2)]$, respectively (see
Figure \ref{f4}).

Clearly, $\lam''_k=\nu_{f''}\big|_{\Del''_k}$ and
$\lam''_l=\nu_{f''}\big|_{\Del''_l}$ are linear functions.
Furthermore, $\nu_{f''}(i,j)>\max\{\lam''_k(i,j),\lam''_l(i,j)\}$
for all points $(i,j)$ satisfying $0\le i<m$, $m_1i+mj<mm_1$,
$m_2i-mj<mm_2$ (i.e., inside the triangle $\Del_z$ with vertices
$(m,0)$, $(0,m_1)$, $(0,-m_2)$). This means, in particular, that
the coefficients $c''_{ij}(t)$ of $(x'')^i(y'')^j$ in $f''$
satisfy $c''_{m,0}(0)\ne 0$, $c''_{i,0}(0)=0$, $i<m$. Hence there
is a unique $\tau(t)\in\K$, $\tau(0)=0$, such that the polynomial
$\widetilde f(\wx,\wy)=f''(\wx+\tau(t),\wy)$ does not contain the
monomial $\wx^{m-1}$.

One can easily see that $\nu_{\widetilde
f}\big|_{\Del''_k}=\lam''_k$, $\nu_{\widetilde
f}\big|_{\Del''_l}=\lam''_l$, and
$\nu_{f''}(i,j)>\max\{\lam''_k(i,j),\lam''_l(i,j)\}$ as far as
$0\le i<m$, $m_1i+mj<mm_1$, $m_2i-mj<mm_2$. It follows that the
subdivision of the Newton polygon of $\widetilde f$ defined by the
function $\nu_{\widetilde f}$ contains a subdivision of the
triangle $\Del_z$, and moreover, this subdivision has no vertices
inside the segments $E_k^z=[(m,0),(0,m_1)]$,
$E_l^z=[(m,0),(0,-m_2)]$ and at the point $(m-1,0)$.

Finally, the fragment of the tropicalization of the polynomial
$\widetilde f$ and of the curve $\widetilde C=\{\widetilde f=0\}$
restricted to the triangle $\Del_z$, we call {\it the
$z$-refinement of the tropicalization of $f$ and of $C$} and
denote by ${\cal T}_z(f)$, ${\cal T}_z(C)$.

\begin{remark}\label{r4}
{\rm Notice that the truncation $\varphi^{(k)}_z(x,y)$ (resp.,
$\varphi^{(l)}_z(x,y)$) of the tropicalization of $\widetilde f$
on the segment $E_k^z$ (resp., $E_l^z$) is uniquely determined by
the polynomials $f_k$, $f_l$ and the point $z$. Any polynomial
with Newton polygon $\Del_z$, whose truncations to the edges
$E_k^z$, $E_l^z$ are just $\varphi_z^{(k)}$, $\varphi_z^{(l)}$,
and the coefficient of $x^{m-1}$ vanishes, will be called a
deformation pattern compatible with $f_k$, $f_l$ and $z$.}
\end{remark}

\medskip

In our situation, $C_k$, $C_l$ are non-singular at $z$, i.e.,
$m_1=m_2=1$, and the refinements of the tropicalization are
described in the following statements.

\begin{lemma}\label{l10}
For a given integer $m\ge 2$ and fixed $a,b,c\in\C^*$, the set of
polynomials $F(x,y)=ay^2+byg(x)+c$ with $g(x)=x^m+...$, $\deg
g=m$, defining plane rational curves, consists of $m$ disjoint
one-dimensional families. Each family has a unique representative
with the zero coefficient of $x^{m-1}$ in $g(x)$, and the rest of
the family can be obtained from the chosen representative by the
coordinate change $x\mapsto x+a$, $a\in\C$. Furthermore, all such
rational curves have $m-1$ nodes in $\C^2$ as their only
singularities.
\end{lemma}

{\bf Proof}. Without loss of generality, suppose that $a=c=1$,
$b=2$. Then the equations for singular points,
$F(x,y)=F_x(x,y)=F_y(x,y)=0$, reduce to the system
\begin{equation}g(x)^2=1,\quad \frac{dg}{dx}(x)=0\ .\label{e33}\end{equation} Solutions to this system are
the $x$-coordinates of singular points, and their multiplicities
in $g'(x)$ are Milnor numbers. Hence the total Milnor number does
not exceed $m-1$. On the other hand, the total Milnor number is at
least the $\del$-invariant. Thus, by our assumptions they
coincide, which is only possible in the case of $m-1$ nodes,
corresponding to $m-1$ distinct solutions to (\ref{e33}). The
latter condition on (\ref{e33}) holds if and only if $g(x)$ is (up
to a shift $x\mapsto x+a$) the Chebyshev polynomial
$\cos(m\cdot\arccos(2^{-(m-1)/m}x))$ or one of the $(m-1)$ other
polynomials of type $g(x\eps)$, $\eps^m=1$, or $-g(x\eps)$,
$\eps^m=-1$. \proofend

\begin{lemma}\label{l113} In the above notation and definitions,
let $z\in\Tor(\sig)\cap C_k\cap C_l$. If $C_k$ and $C_l$ are
non-singular at $z$, the intersection number of $C_k$ and
$\Tor(\sig)$ at $z$ is $(C_k\cdot\Tor(\sig))_z=m\ge 2$, and $z$
bears singularities with the total $\del$-invariant $m-1$, then
${\cal T}_z(f)$ consists of one polynomial with Newton polygon
$\Del_z$, which defines a rational curve with $m-1$ nodes in
$\Tor(\Del_z)$, and the singularities born of $z$ are $m-1$ nodes.
\end{lemma}

{\bf Proof}. We show that ${\cal T}_z(f)$ consists of the triangle
$\Del_z$ and, correspondingly, of one polynomial, which then is
described in Lemma \ref{l10}.

We use induction on $m$. Let $m=2$, and so ${\cal T}_z(f)$
contains more than one polynomial. Then $\Del_z$ is subdivided
into two triangles, \mbox{$\conv\{(0,0),(0,1),(2,0)\}$} and
\mbox{$\conv\{(0,0),(0,-1),(2,0)\}$}. The curves defined by the
polynomials with these Newton triangles are non-singular and cross
$\Tor([(0,0),(2,0)])$ transversally, since the coefficient of $x$
vanishes, but then no singular point appears in the deformation by
Lemma \ref{l6}. Let $m\ge 3$ and
$(i_1,0),...,(i_r,0)\in\Int(\Del_z)$ be the vertices of the
subdivision of $\Del_z$ associated with ${\cal T}_z(f)$, $r\ge 1$,
$0\le i_1<...<i_r\le m-2$. Singular points may appear only from
possible tangency points along $\Tor([(i_s,0),(i_{s+1},0)])$,
$s=1,...,r-1$, or along $\Tor([(i_r,0),(m,0)])$, or from a curve
with Newton triangle \mbox{$\conv\{(0,1),(0,-1),(i_1,0)\}$}, if
$i_1>0$. Since the curves with Newton triangles containing the
edge $[(i_r,0),(m,0)]$ must cross $\Tor([(i_r,0),(m,0)])$ at least
at two points due to the condition that the coefficient of
$\wx^{m-1}$ vanishes, the total $\del$-invariant of singular
points which may appear is at most
$$\max\{0,\ i_1-1\}+\sum_{s=1}^{r-1}(i_{s+1}-i_s-1)+(m-i_r-2)\le
m-2\ ,$$ which gives a contradiction. \proofend

\subsection{Refinement of the tropicalization along a non-isolated singularity}
\label{sec12} 
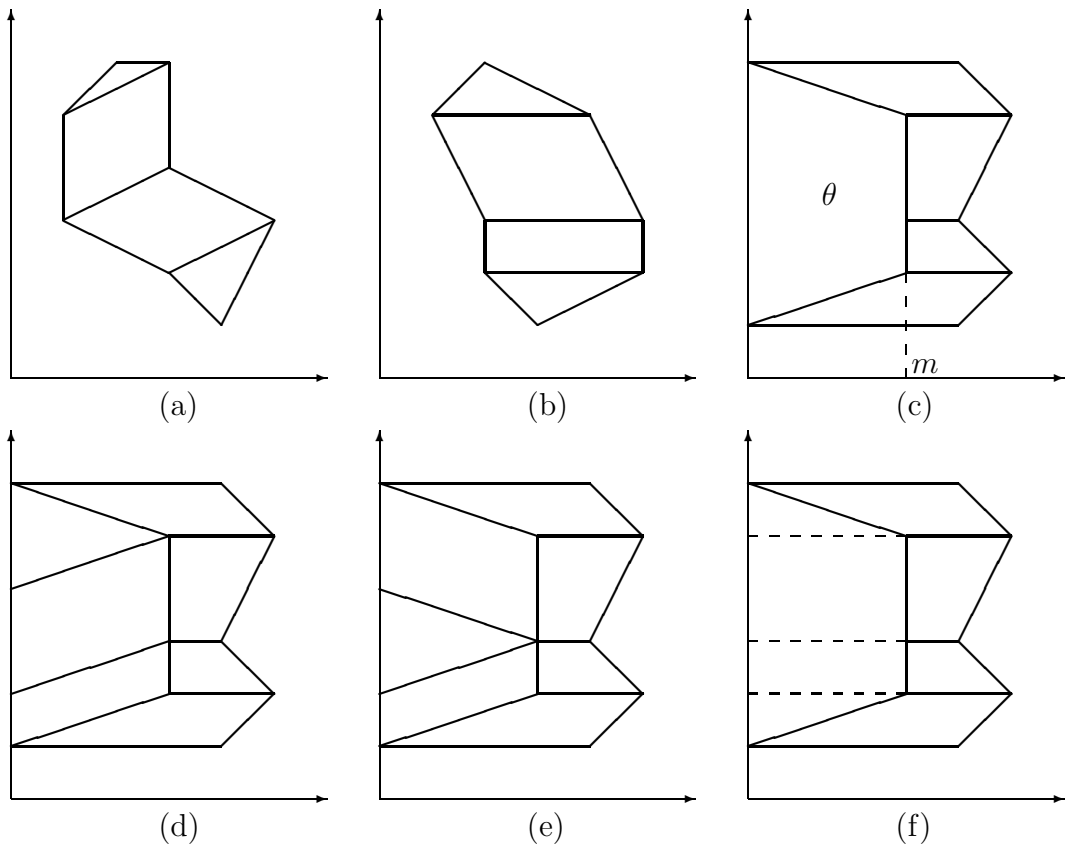
\begin{figure}
\setlength{\unitlength}{0.7cm}
\begin{picture}(14,16)(0,0)
\thinlines \put(0,9){\vector(1,0){6}}\put(0,9){\vector(0,1){7}}
\thicklines\put(1,12){\line(0,1){2}}\put(1,12){\line(2,1){2}}
\put(1,12){\line(2,-1){2}}\put(1,14){\line(1,1){1}}\put(1,14){\line(2,1){2}}
\put(2,15){\line(1,0){1}}\put(3,11){\line(1,-1){1}}
\put(3,11){\line(2,1){2}}\put(3,13){\line(0,1){2}}\put(3,13){\line(2,-1){2}}
\put(4,10){\line(1,2){1}}\put(2.8,8.3){\text{\rm (a)}} \thinlines
\put(7,9){\vector(1,0){6}}\put(7,9){\vector(0,1){7}}
\thicklines\put(8,14){\line(1,-2){1}}\put(8,14){\line(1,1){1}}
\put(9,11){\line(0,1){1}}\put(9,11){\line(1,0){3}}\put(8,14){\line(1,0){3}}
\put(9,11){\line(1,-1){1}}\put(9,12){\line(1,0){3}}
\put(9,15){\line(2,-1){2}}\put(10,10){\line(2,1){2}}\put(11,14){\line(1,-2){1}}
\put(12,11){\line(0,1){1}}\put(9.8,8.3){\text{\rm (b)}} \thinlines
\put(14,9){\vector(1,0){6}}\put(14,9){\vector(0,1){7}}
\thicklines\put(14,10){\line(3,1){3}}\put(14,10){\line(1,0){4}}
\put(18,10){\line(1,1){1}}\put(17,11){\line(1,0){2}}\put(17,11){\line(0,1){1}}
\put(17,12){\line(1,0){1}}\put(17,12){\line(0,1){2}}
\put(18,12){\line(1,-1){1}}\put(17,14){\line(1,0){2}}\put(18,12){\line(1,2){1}}
\put(14,15){\line(3,-1){3}}\put(14,15){\line(1,0){4}}\put(18,15){\line(1,-1){1}}\put(16.8,8.3){\text{\rm
(c)}}\thinlines\dashline{0.2}(17,9)(17,11)\put(17.1,9.1){$m$}\put(15.4,12.3){$\theta$}
\put(0,1){\vector(1,0){6}}\put(0,1){\vector(0,1){7}}
\thicklines\put(0,2){\line(3,1){3}}\put(0,2){\line(1,0){4}}
\put(4,2){\line(1,1){1}}\put(3,3){\line(1,0){2}}\put(3,3){\line(0,1){1}}
\put(3,4){\line(1,0){1}}\put(3,4){\line(0,1){2}}
\put(4,4){\line(1,-1){1}}\put(3,6){\line(1,0){2}}\put(4,4){\line(1,2){1}}
\put(0,7){\line(3,-1){3}}\put(0,7){\line(1,0){4}}\put(4,7){\line(1,-1){1}}
\put(0,3){\line(3,1){3}}\put(0,5){\line(3,1){3}}\put(2.8,0.3){\text{\rm
(d)}}\thinlines\put(7,1){\vector(1,0){6}}\put(7,1){\vector(0,1){7}}
\thicklines\put(7,2){\line(3,1){3}}\put(7,2){\line(1,0){4}}
\put(11,2){\line(1,1){1}}\put(10,3){\line(1,0){2}}\put(10,3){\line(0,1){1}}
\put(10,4){\line(1,0){1}}\put(10,4){\line(0,1){2}}
\put(11,4){\line(1,-1){1}}\put(10,6){\line(1,0){2}}\put(11,4){\line(1,2){1}}
\put(7,7){\line(3,-1){3}}\put(7,7){\line(1,0){4}}\put(11,7){\line(1,-1){1}}
\put(7,3){\line(3,1){3}}\put(7,5){\line(3,-1){3}}\put(9.8,0.3){\text{\rm
(e)}}\thinlines\put(14,1){\vector(1,0){6}}\put(14,1){\vector(0,1){7}}
\thicklines\put(14,2){\line(3,1){3}}\put(14,2){\line(1,0){4}}
\put(18,2){\line(1,1){1}}\put(17,3){\line(1,0){2}}\put(17,3){\line(0,1){1}}
\put(17,4){\line(1,0){1}}\put(17,4){\line(0,1){2}}
\put(18,4){\line(1,-1){1}}\put(17,6){\line(1,0){2}}\put(18,4){\line(1,2){1}}
\put(14,7){\line(3,-1){3}}\put(14,7){\line(1,0){4}}\put(18,7){\line(1,-1){1}}\put(16.8,0.3){\text{\rm
(f)}}\thinlines\dashline{0.2}(14,3)(17,3)\dashline{0.2}(14,4)(17,4)\dashline{0.2}(14,6)(17,6)
\end{picture}
\caption{Refinement of the tropicalization, II}\label{f5}
\end{figure}

Assume that $z\in\Tor(\sig)$,
$\sig=\Del_k\cap\Del_l$ is a common edge,
$(C_k\cdot\Tor(\sig))_z=(C_l\cdot\Tor(\sig))_z=m\ge 2$, and at
least one of the $\Del_k,\Del_l$ is a parallelogram. Then the pair
$\Del_k,\Del_l$ extends up to a chain (after a renumbering)
$\Del_1,...,\Del_p$, $p\ge 3$, where $\Del_1,\Del_p$ are
triangles, $\Del_2,...,\Del_{p-1}$ are parallelograms,
$\sig_1=\Del_1\cap\Del_2,...,\sig_{p-1}=\Del_{p-1}\cap\Del_p$ are
common edges which are parallel to each other (see Figure
\ref{f5}(a)). We shall associate a refinement of the
tropicalization with the union $Z$ of the multiple components of
the curves $C_2,...,C_{p-1}$ which cross the lines
$\Tor(\sig_1),...,\Tor(\sig_{p-1})$.

Multiplying $f(x,y)$ be a suitable constant from $\K^*$, we can
achieve the constancy of $\nu_f$ along the edges
$\sig_1,...,\sig_{p-1}$. Then we apply $M\in\Aff(\Z^2)$, which
puts $\Del$ into the right half-plane and makes
$\sig_1,...,\sig_{p-1}$ horizontal (Figure \ref{f5}(b)). The
corresponding monomial coordinate change transforms $f(x,y)$ into
a polynomial $f'(x',y')$. For the latter polynomial, the
truncations of the edges $\sig_1,...,\sig_{p-1}$ of the
tropicalizations to $\Del_1,...,\Del_p$ contain a factor
$(x-\xi)^m$ with some $\xi\in\C^*$. Then we introduce the
polynomial $f''(x'',y'')=f'(x''+\xi,y'')$ and consider its
tropicalization. Again, for better understanding of ${\cal
T}(f'')$, we apply the above coordinate changes to the polynomials
$$P_k(x,y)=\sum_{(i,j)\in\Del_k}c^0t^{\lam_k(i,j)}x^iy^j,\quad
k=1,...,p\ ,$$ where
$$f_k(x,y)=\sum_{(i,j)\in\Del_k}c^0x^iy^j,\quad\lam_k=\nu_f\Big|_{\Del_k},\quad
k=1,...,p\ .$$ The Newton polygons $\Del''_k$ of the resulting
polynomials $P''_k(x'',y'')$, $i=1,...,p$, appear as shown in
Figure \ref{f5}(c) and bound a trapezoid $\theta$ with vertices
$(0,a-1)$, $(0,b+1)$, $(m,a)$, $(m,b)$. Then, in particular,
$\nu_{f''}\big|_{\Del''_k}=\lam''_k$, $k=1,...,p$, are linear
functions, and
\begin{equation}\nu_{f''}(i,j)>\max_{1\le k\le p}\lam''_k(i,j),\quad 0\le
i<m,\ a\le j\le b\ .\label{e52}\end{equation}

Consider now subdivisions of $\theta$ into parallelograms and one
triangle with edges parallel to the edges of $\theta$ (see, for
example, Figure \ref{f5}(d,e)). Exactly one of them can be induced
by a convex piece-wise linear function, defined as $\nu_{f''}$ on
$\Del''_1\cup...\cup\Del''_p$ and extended to $\theta$. Here we
suppose that the function $\nu_f$ is generic among the convex
piece-wise linear functions determining the same subdivision $S_f$
of $\Del$ (that means it had generic rational slopes before we
have multiplied it by a large natural number), and then its graph
necessarily has a break along the edges of the triangle. Let
$(m,d)$ be a vertex of the triangle on the chosen subdivision of
$\theta$. We perform one more shift now. Namely, in view of
(\ref{e52}) and the fact that the linear functions $\lam''_k$,
$k=1,...,p$, are constant in the horizontal direction, we conclude
that there exists a unique $\tau(t)\in\K$, $\tau(0)=0$, such that
the polynomial $\widetilde f(\wx,\wy)=f''(\wx+\tau(t),\wy)$ has no
monomial $\wx^{m-1}\wy^d$ (next to the vertex of the triangle).

We claim that the function $\nu_{\widetilde f}$ defines the
subdivision of $\theta$ into one triangle and $p-1$ parallelograms
as described above. Furthermore, the tropicalizations of
$\widetilde f$ on the parallelograms inside $\theta$ are products
of binomials, and the tropicalization of $\widetilde f$ on the
triangle inside $\theta$ (which we denote $\Del_Z$) is
$\wy^dP(\wx,\wy)$, where $P$ is a polynomial from Lemma \ref{l7}
with the vanishing coefficient of $\wx^{m-1}\wy$.

Indeed, deformation of the tropicalization of the curve
$\widetilde C=\{\widetilde f=0\}$ on the polygons $\Del'\subset
\theta$, $\Del'\in P(S_{\widetilde f})$, describes the deformation
of the tropicalization $C^{(0)}$ of the original curve $C=\{f=0\}$
in a neighborhood of $Z$. The argument from section \ref{sec9}
implies, first, that the truncations of the tropicalization of
$\widetilde f$ of the vertical edges of
$\Del''_2,...,\Del''_{p-1}$, lying on $\partial\theta$ (see Figure
\ref{f5}(c)), define irreducible components of the restriction of
${\cal T}(\widetilde C)$ to the polygons subdividing $\theta$. In
particular, each vertical edge of $\Del''_2,...,\Del''_{p-1}$,
lying on $\partial\theta$, is joined to a segment, lying on the
vertical coordinate axis, by a sequence of parallelograms.
Furthermore, the argument from section \ref{sec9} yields that all
other components of ${\cal T}(\widetilde C)$ in $\Tor(\Del')$,
$\Del'\subset\theta$, must be rational, and each of them crosses
$\bigcup_\sig\Tor(\sig)$ at two points, where $\sig$ runs over all
non-vertical edges of the subdivision $S_{\widetilde f}$ in
$\theta$. All this leaves only the possibility proclaimed above
for ${\cal T}(\widetilde C)\big|_\theta$. For instance, the
subdivision cannot be as shown by dashes in Figure \ref{f5}(f),
since the tropicalization of $\widetilde f$ of the horizontal
dashed segment $\sig$, lying on the level $d$, cannot be a power
of a binomial (notice that the monomial $\wx^{m-1}\wy^d$ is
absent). Thus the component of ${\cal T}(\widetilde C)$
corresponding to $\sig$ is irreducible and crosses $\Tor(\sig)$ at
least at two points.

\begin{remark}\label{r5}
{\rm Notice that the truncations $\varphi_Z^{(1)}$ and
$\varphi_Z^{(p)}$ of the tropicalization of $\widetilde f$ on the
non-vertical edges of the triangle $\Del_Z$ are uniquely
determined by the polynomials $f_1$ and $f_p$, respectively. Any
polynomial with Newton triangle $\conv\{(0,0),(m,1),(0,2)\}$
($\Del_Z$ shifted down), whose truncations to the non-vertical
edges coincide with $\varphi_Z^{(1)}$, $\varphi_Z^{(p)}$ (up to
multiplication by a suitable monomial), and the coefficient of
$x^{m-1}y$ vanishes, will be called a deformation pattern for the
set $Z$.}
\end{remark}

\subsection{Restoring a nodal curve out of tropical
data}\label{sec11} Denote by ${\cal Q}_\Del(nA_1)$ the set of
quadruples $(A,S,F,R)$, where \begin{itemize}\item $A\in{\cal
A}(\Del)$ is a nodal amoeba of rank $r$, $S:\
\Del=\Del_1,...,\Del_N$ is a subdivision of $\Del$ dual to $A$,
and $F$, $R$ are collections of the following polynomials in
$\C[x,y]$ which together are defined up to multiplication by the
same non-zero (complex) constant;
\item $F=(f_1,...,f_N)$, $f_i$ is a polynomial with
Newton polygon $\Del_i$, $i=1,...,N$, such that, if $\Del_i$ is a
triangle, then $f_i$ defines a rational curve in $\Tor(\Del_i)$ as
described in Lemma \ref{l7}, if $\Del_i$ is a parallelogram, then
$f_i$ defines a curve in $\Tor(\Del_i)$ as described in Lemma
\ref{l11}, and, for any common edge $\sig=\Del_i\cap\Del_j$, the
truncations $f_i^\sig$ and $f_j^\sig$ coincide; \item $R$ is a
collection of deformation patterns compatible with $F$ as defined
in Remarks \ref{r4} and \ref{r5}.
\end{itemize}

We are given the points $\bx_1,...,\bx_r\in\R^2$ and
$\bp_1,...,\bp_r\in(\K^*)^2$ such that $\val(\bp_i)=\bx_i$,
$i=1,...,r$, and we intend to find \begin{itemize}\item how many
elements $(A,S,F,R)\in{\cal Q}_\Del(nA_1)$ correspond to a nodal
amoeba $A\in{\cal A}(\Del)$ of rank $r$ passing through
$\bx_1,...,\bx_r$, and \item how many polynomials $f\in\K[x,y]$
(determined up to multiplication by a non-zero $\K$-constant) with
Newton polygon $\Del$, which define curves $C\in\Sig_\Del(nA_1)$
passing through $\bp_1,...,\bp_r$, arise from a tropicalization
$(A,S,F,R)\in{\cal Q}_\Del(nA_1)$.
\end{itemize}

\medskip

{\it Step 1}. Let $A\in{\cal A}(\Del)$ be a nodal amoeba of rank
$r$ passing through the given points $\bx_1,...,\bx_r\in\R^2$.

Observe, first, that $A$ uniquely determines a dual subdivision
$S$ of $\Del$. Indeed, the unbounded components of $\R^2\backslash
A$ are in a natural one-to-one correspondence with
$\partial\Del\cap\Z^2$. The bounded edges of $A$ in the boundary
of the above components define germs of the edges of $S$ starting
at $\partial\Del\cap\Z^2$. There is a pair of non-parallel
neighboring germs which start at distinct points of
$\partial\Del\cap\Z^2$, and their extension uniquely determines a
triangle or a parallelogram in the subdivision $S$. Then we remove
this polygon out of $\Del$ and continue the process.

Second, $A$ determines (uniquely up to a constant shift) a convex
piece-wise linear function $\nu:\Del\to\R$ whose graph projects
onto the subdivision $S$. More precisely, the points
$\bx_1,...,\bx_r$ lie on $r$ distinct edges of $A$ which
correspond to some $r$ edges of $S$. If $\sig_i\in E(S)$
corresponds to a point $\bx_i$, and $\bi'_i,\bi''_i$ are the
endpoints of $\sig_i$, $1\le i\le r$, then we have linear
conditions on $\nu(\bi'_i)$ and $\nu(\bi''_i)$:
\begin{equation}\nu(\bi'_i)-\nu(\bi''_i)=(\bi''_i-\bi'_i)\bx_i,\quad
i=1,...,r\ .\label{e144}\end{equation} Since $\bx_1,...,\bx_r$ are
generic, system (\ref{e144}) is independent. Furthermore,
parallelograms $\Del_j\in P(S)$, $j=1,...,N_4$, corresponding to
the $4$-valent vertices of $A$, impose the following linear
conditions on the values of $\nu$ at the vertices
$\bi_j^{(1)},\bi_j^{(2)},\bi_j^{(3)},\bi_j^{(4)}$ of $\Del_j$
(listed, say, clockwise):
\begin{equation}\nu_f(\bi_j^{(1)})+\nu_f(\bi_j^{(3)})=\nu_f(\bi_j^{(2)})+\nu_f(\bi_j^{(4)}),\quad
j=1,...,N_4\ .\label{e145}\end{equation} Lemma \ref{l12} yields
that the united system (\ref{e144}), (\ref{e145}) is independent,
and, since it contains $|V(S)|-1$ equations, it determines the
values of $\nu$ at the vertices of $S$ uniquely up to a constant
shift.

\medskip

{\it Step 2}. We are looking for polynomials of the form
\begin{equation}f(x,y)=\sum_{(i,j)\in\Del}\widetilde c_{ij}(t)t^{\nu(i,j)}x^iy^j,\quad
\widetilde c_{ij}(0)=c_{ij},\ (i,j)\in\Del\
,\label{e63}\end{equation} where
$$f_k(x,y)=\sum_{(i,j)\in\Del_i}c_{ij}x^iy^j,\quad
k=1,...,N\ .$$ We claim that the condition
\begin{equation}f(\bp_1)=...=f(\bp_r)=0\label{e58}\end{equation} uniquely determines the
coefficients of $f_1,...,f_N$ at the vertices of $S$, and the
truncations of $f_1,...,f_N$ on the edges $\sig_1,...,\sig_r$,
corresponding to $\bx_1,...,\bx_r$, up to multiplication by the
same non-zero constant.

Indeed, let $$\bx_i=(-\alp_i,-\bet_i),\quad
\bp_i=(\xi,\eta),\quad\xi=\xi^0_it^\alp_i+\
\text{h.o.t.},\quad\eta=\eta^0_it^\bet_i+\ \text{h.o.t.},\quad
\xi^0_i,\eta^0_i\in\C^*\ ,$$ and let the endpoints of the edge
$\sig_i$ be $\bi'_i=(i_1,j_1)$, $\bi''_i=(i_2,j_2)$. The
conditions $f(\bp_i)=0$, $i=1,...,r$, then transform into the
following equations:
$$f(\bp_i)=t^{\nu_f(i_1,j_1)+i_1\alp_i+j_1\bet_i}\left(g_i(\xi^0_i,\eta^0_i)+O(t)\right)=0$$
\begin{equation}\Longrightarrow\quad g_i(\xi^0_i,\eta^0_i)=0\
,\label{e146}\end{equation} where a quasihomogeneous polynomial
$g_i(x,y)=c_{\bi'_i}x^{i_1}y^{j_1}+...+c_{\bi''_i}x^{i_2}y^{j_2}$
is the tropicalization of $f^{\sig_i}$. Since $g_i$ is the product
of a monomial and a power of an irreducible binomial, (\ref{e146})
determines it uniquely up to a constant factor. On the other hand,
the coefficients $b_1,b_2,b_3,b_4$ of the polynomial $f_j$, having
Newton parallelogram $\Del_j$ with respective clockwise ordered
vertices $\bi_j^{(1)},\bi_j^{(2)},\bi_j^{(3)},\bi_j^{(4)}$,
satisfy $b_1b_3=b_2b_4$. We see that all the conditions imposed on
the coefficients of $f_1,...,f_N$ at $V(S)$ are just a
multiplicative form of a system like (\ref{e144}), (\ref{e145}),
and hence the claim follows.

\medskip

{\it Step 3}. We next compute how many ways are there to restore
$f_1,...,f_N$, if one fixes the coefficients of the polynomials
$f_1,...,f_N$ at $V(S)$ and their truncations to the edges
$\sig_1,...,\sig_r$.

To formulate the answer, consider the amoeba $A$ and introduce the
set $E^*(A)$ of extended edges of $A$, i.e., the maximal straight
line intervals contained in $A$. That is, the edges of $A$, which
lie on the same straight line and are connected by four-valent
vertices, we join into one extended edge. We then claim that there
are $W(A)\prod_{\sig^*\in
E^*(A)}|\sig^*|^{-1}\prod_{i=1}^r|\sig_i|^{-1}$ collections of
polynomials $f_1,...,f_N$ compatible with the given data, where
$|\sig^*|$ is understood as the length of one of the edges of $S$
dual to $\sig^*\in E^*(A)$.

Let $G_0$ be the union of all edges in $E^*(A)$, passing through
$\bx_1,...,\bx_r$. It is a union of trees, whose vertices
different from the four-valent vertices of $A$ have valency at
most two. Notice that the position of the remaining edges in
$E^*(A)$ is prescribed by $G_0$. This yields that the graph $G_0$
has a bivalent vertex $v_0$, and there is an extended edge
$\eps_1$ in $A\backslash G_0$ starting at $v_0$. Notice that
$\eps_1$ cannot join two vertices of $G_0$ of valency two. Indeed,
the position of two non-adjacent bivalent vertices of $G_0$ is
uniquely determined by the position of some four points among
$\bx_1,...,\bx_r$, and thus, due to the generality of the latter
points, the straight line through the given vertices is not
orthogonal to any of the segments joining integral points in
$\Del$. Put $G_1=G_0\cup\eps_1$. Next, by a similar reason, there
is a bivalent vertex $v_1$ of $G_1$, and the extended edge
$\eps_2$ of $A\backslash G_1$, starting at $v_1$, which does not
end up at another bivalent vertex of $G_1$. Proceeding in the same
manner, we reconstruct the whole amoeba $A$.

By Lemma \ref{l7}, for a triangle $\Del_k$, $1\le k\le N$ and
given truncations to two edges $\sig',\sig''$ of $\Del_k$, an
admissible polynomial $f_k$ (i.e., defining a rational nodal curve
with precisely three unibranch intersection points with
$\Tor(\partial\Del_k)$) can be restored in
$|\Del_k|(|\sig'|\cdot|\sig''|)^{-1}$ ways. Following the
construction of the graphs $G_0,G_1,...$, we obtain that the
collection of polynomials $f_k$ with Newton triangles and
compatible with the given data, can be restored in
$W(A)\prod_{\sig^*\in
E^*(A)}|\sig^*|^{-1}\prod_{i=1}^r|\sig_i|^{-1}$ ways, and the
polynomials with Newton parallelograms are then determined
uniquely.

\medskip

{\it Step 4}. By Lemma \ref{l10}, any collection $f_1,...,f_N$ can
be completed with any of $\prod_{\sig\in E^*(A)}|\sig|$
collections of deformation patterns; hence we can find
$W(A)\prod_{i=1}^r|\sig_i|^{-1}$ elements $(A,S,F,R)\in{\cal
Q}_\Del(nA_1)$ compatible with the given nodal amoeba $A$ and the
points $\bx_1,...,\bx_r\in\R^2$, $\bp_1,...,\bp_r\in(\K^*)^2$.

We complete the proof of Theorem \ref{t5} with the following
statement, which will be proven after the main patchworking
Theorem \ref{t1} in section \ref{sec300}

\begin{lemma}\label{l18} In the above notation, given the points $\bx_1,...,\bx_r\in\R^2$,
$\bp_1,...,\bp_r\in(\K^*)^2$, and a compatible $(A,S,F,R)\in{\cal
Q}_\Del(nA_1)$, the polynomials $f\in\K[x,y]$ with Newton polygon
$\Del$, which tropicalize into $(A,S,F,R)$, define exactly
$\prod_{i=1}^r|\sig_i|$ curves in $\Lam_\K(\Del)$ having $n$ nodes
and passing through $\bp_1,...,\bp_r$.
\end{lemma}

\section{Counting curves with one cusp}\label{sec6}

\subsection{Formulation of the result}
Let $\Del$ be a non-degenerate convex lattice polygon having at
least one interior integral point. We are interested to find the
degree of the variety $\Sig_\Del(A_2)$ of curves
$C\in\Lam_\K(\Del)$ having an ordinary cusp as their only
singularity, and we intend to express this degree as the number of
certain non-Archimedean amoebas.

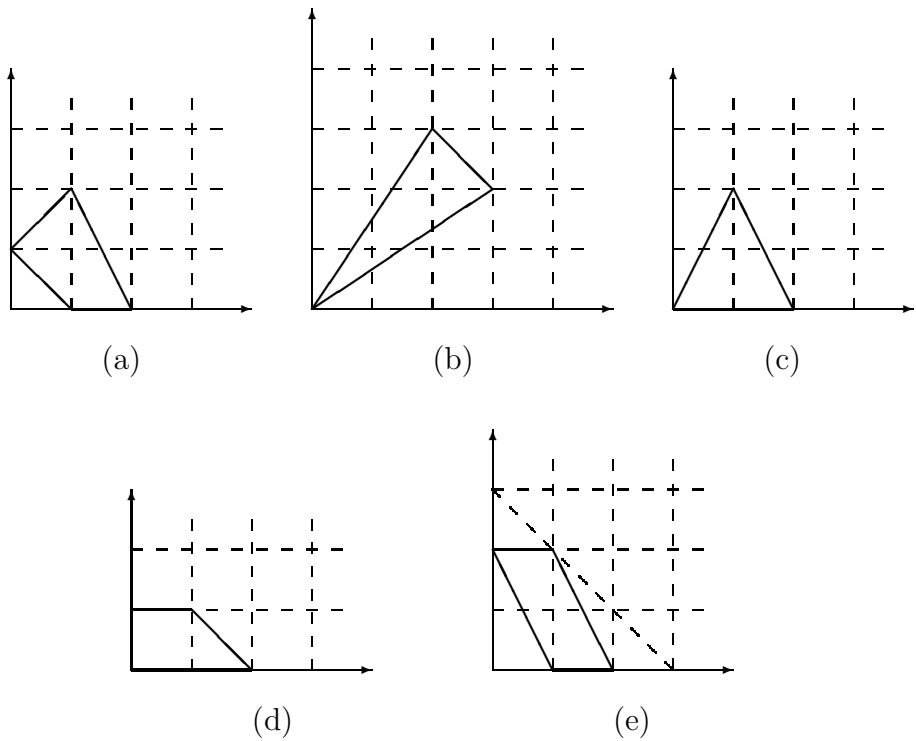
\begin{figure}
\setlength{\unitlength}{0.8cm}
\begin{picture}(17,12)(0,0)
\thinlines\put(12,7){\vector(0,1){4}}\put(9,1){\vector(0,1){4}}\put(1,7){\vector(0,1){4}}
\put(3,1){\vector(0,1){3}}\put(6,7){\vector(0,1){5}}
\put(1,7){\vector(1,0){4}}\put(9,1){\vector(1,0){4}}\put(12,7){\vector(1,0){4}}
\put(6,7){\vector(1,0){5}}\put(3,1){\vector(1,0){4}}
\dashline{0.2}(12,8)(15.5,8)\dashline{0.2}(12,9)(15.5,9)
\dashline{0.2}(12,10)(15.5,10)\dashline{0.2}(1,8)(4.5,8)
\dashline{0.2}(1,9)(4.5,9)\dashline{0.2}(1,10)(4.5,10)
\dashline{0.2}(3,2)(6.5,2)\dashline{0.2}(3,3)(6.5,3)
\dashline{0.2}(6,8)(10.5,8)\dashline{0.2}(6,9)(10.5,9)
\dashline{0.2}(6,10)(10.5,10)\dashline{0.2}(6,11)(10.5,11)
\dashline{0.2}(13,7)(13,10.5)\dashline{0.2}(14,7)(14,10.5)
\dashline{0.2}(15,7)(15,10.5)\dashline{0.2}(4,1)(4,3.5)
\dashline{0.2}(5,1)(5,3.5)\dashline{0.2}(6,1)(6,3.5)
\dashline{0.2}(2,7)(2,10.5)\dashline{0.2}(3,7)(3,10.5)\dashline{0.2}(4,7)(4,10.5)
\dashline{0.2}(7,7)(7,11.5)\dashline{0.2}(8,7)(8,11.5)\dashline{0.2}(9,7)(9,11.5)
\dashline{0.2}(10,7)(10,11.5)\dashline{0.2}(9,2)(12.5,2)\dashline{0.2}(9,3)(12.5,3)
\dashline{0.2}(9,4)(12.5,4)\dashline{0.2}(10,1)(10,4.5)\dashline{0.2}(11,1)(11,4.5)
\dashline{0.2}(12,1)(12,4.5)\dashline{0.2}(12,1)(9,4)\thicklines
\put(12,7){\line(1,2){1}} \put(14,7){\line(-1,2){1}}
\put(12,7){\line(1,0){2}} \put(3,1){\line(0,1){1}}
\put(3,1){\line(1,0){2}}
\put(3,2){\line(1,0){1}}\put(4,2){\line(1,-1){1}}\put(1,8){\line(1,-1){1}}
\put(1,8){\line(1,1){1}}\put(2,7){\line(1,0){1}}\put(3,7){\line(-1,2){1}}
\put(6,7){\line(2,3){2}}\put(6,7){\line(3,2){3}}\put(8,10){\line(1,-1){1}}
\put(10,1){\line(-1,2){1}}\put(10,1){\line(1,0){1}}\put(11,1){\line(-1,2){1}}
\put(9,3){\line(1,0){1}}\put(13.5,6){{\rm (c)}} \put(5,0){{\rm
(d)}}\put(11,0){{\rm (e)}}\put(2.5,6){{\rm (a)}} \put(8,6){{\rm
(b)}}
\end{picture}
\caption{Newton polygons of auxiliary curves with nodes and cusps
}\label{f3}
\end{figure}

An amoeba $A\in{\cal A}(\Del)$ is called {\it 1-cuspidal} if its
dual subdivision $S$ of $\Del$ satisfies one of the following
conditions:
\begin{enumerate} \item[(i)] $S$ contains a quadrangle,
$\Aff(\Z^2)$-equivalent to that shown in Figure \ref{f3}(a), and
the rest of $S$ consists of triangles of area $1/2$;
\item[(ii)] $S$ contains a triangle, $\Aff(\Z^2)$-equivalent to that
shown in Figure \ref{f3}(b), and the rest of $S$ consists of
triangles of area $1/2$;
\item[(iii)] $S$ contains an edge of length $2$ common for a triangle,
$\Aff(\Z^2)$-equivalent to that shown in Figure \ref{f3}(c), and
for a triangle of area $1$, and the rest of $S$ consists of
triangles of area $1/2$;
\item[(iv)] $S$ contains an edge of length $2$ common for a quadrangle,
$\Aff(\Z^2)$-equivalent to that shown in Figure \ref{f3}(d), and
for a triangle of area $1$, and the rest of $S$ consists of
triangles of area $1/2$; \item[(v)] $S$ contains an edge of length
$3$, common for two triangles of area $3/2$, and the rest of $S$
consists of triangles of area $1/2$.
\end{enumerate}

Observe that a $1$-cuspidal amoeba has rank $r=|\Z^2\cap\Del|-3$
and determines the dual subdivision uniquely.

Let $\bx_1,...,\bx_r$ be generic points in $\R^2$, $A$ a
1-cuspidal amoeba passing through these points. We shall introduce
the weight $W(A,\bx_1,...,\bx_r)$.

Assume that $A$ has only $3$-valent vertices, i.e., $S$ contains
only triangles. Then put $W(A,\bp_1,...,\bp_r)$ equal to $5$, $6$,
or $6$ in accordance with cases (ii), (iii), or (v) in the
definition of 1-cuspidal amoebas.

Assume that $A$ contains a quadrangle $\Del'$. The vertices of $S$
and the $r$ edges of $S$ dual to the edges of $A$, which contain
the fixed points, form a graph $\Gamma$. Observe that $\Gamma$ has
no cycles, since, otherwise, as shown in Step 2 of section
\ref{sec11}, one would have a dependent sequence of relations of
type (\ref{e144}), which is impossible in view of the generic
choice of $\bp_1,...,\bp_r$. Thus, $|V(S)|=|\Del\cap\Z^2|-1=r+2$
yields that $\Gamma$ consists of two disjoint trees (a tree may be
one point). Furthermore, the vertices of $\Del'$ cannot all belong
to one component of $\Gamma$, and, for the case of $\Del'$ shown
in Figure \ref{f3}(d), it cannot be that the two upper vertices
belong to one component of $\Gamma$ and the lower vertices belong
to the other. Take the vectors joining the vertices of $\Del'$,
belonging to the same component of $\Gamma$, and take one vector
$v$ joining two vertices from distinct components of $\Gamma$, and
denote by $w(A,\bx_1,...,\bx_r)$ the minimal positive coefficient
of $v$ in the possible linear combinations with integral
coefficient of all the vectors taken. Now put
$W(A,\bx_1,...,\bx_r)$ equal to $w(A,\bx_1,...,\bx_r)$ or
$3w(A,\bx_1,...,\bx_r)$ in the cases (i) and (iv), respectively.

\begin{theorem}\label{t6} In the above notation, $$\deg\Sig_\Del(A_2)=\sum
W(A,\bx_1,...,\bx_n)\ ,$$ where $\bx_1,...,\bx_r\in\R^2$ is a
collection of generic distinct points, and $A$ ranges over all
1-cuspidal amoebas in ${\cal A}(\Del)$, passing through
$\bx_1,...,\bx_r$.
\end{theorem}

\subsection{Auxiliary curves with nodes and cusps}
We start by describing the nodal complex curves which will be used
in the proof of Theorem \ref{t6}.

\begin{lemma}\label{l17}
Up to the action of the group $\Aff(\Z^2)$ of affine automorphisms
of $\Z^2$,
\begin{itemize}
\item the polygons in Figure \ref{f3}(a,e) are the only lattice quadrangles with one
interior integral point and all edges of length $1$;
\item the polygon in Figure \ref{f3}(b) is the only lattice
triangle with two interior integral points and all edges of length
$1$;
\item the polygon in Figure \ref{f3}(c) is the only lattice
triangle with one interior integral point, one edge of length $2$
and the others of length $1$;
\item the polygon in Figure \ref{f3}(d) is the only lattice
quadrangle without interior integral points, with one edge of
length $2$ and the others of length $1$.
\end{itemize}
There is no lattice pentagon which contains vertices as its only
integral points.
\end{lemma}

This is an elementary geometric fact, and we omit the proof.

\begin{lemma}\label{l16} Denote by $\Del_i$, $i=1,...,5$, the polygons shown in Figure
\ref{f3}(a-e), respectively. \begin{enumerate}\item[(i)] A curve
in $\Tor(\Del_i)$, defined by a polynomial with Newton polygon
$\Del_i$, cannot have a singularity more complicated than an
ordinary cusp if $i=1,2$, and has at most one node if $i=3,4$.
\item[(ii)] Given the coefficients at the vertices of $\Del_2$, there exist exactly
five polynomials with Newton triangle $\Del_2$ defining a curve
with a cusp. Furthermore, such curves have no other singular
points.
\item[(iii)] Given the coefficients at the vertices of $\Del_3$,
there exist exactly two polynomials with Newton triangle $\Del_3$
defining a curve with a node, which lies on $\Tor([(0,0),(2,0)])$.
\item[(iv)] A polynomial with Newton polygon $\Del_1$ and
coefficients $c_{10},c_{01},c_{12},c_{20}\in\C^*$ at the vertices
of $\Del_1$ defines a curve in $\Tor(\Del_1)$, having a cusp, if
and only if
\begin{equation}c_{10}^3c_{1,2}=c_{01}^2c_{20}^2\ .\label{e66}\end{equation} Moreover, for fixed
$c_{10},c_{01},c_{12},c_{20}$, such a polynomial is unique and the
corresponding curve has a cusp as its only singularity.
\item[(v)] A polynomial with Newton
polygon $\Del_4$ and coefficients
$c_{00},c_{01},c_{1,1},c_{2,0}\in\C^*$ at the vertices of $\Del_4$
defines a curve in $\Tor(\Del_4)$, having a node on
$\Tor([(0,0),(2,0)])$, if and only if
\begin{equation}c_{00}c_{1,1}^2=c_{01}^2c_{20}\ .\label{e67}\end{equation} Moreover, for fixed
$c_{00},c_{01},c_{1,1},c_{2,0}$, such a polynomial is unique and
the corresponding curve has a node as its only singularity.
\item[(vi)] A curve in $\Tor(\Del_5)$ defined by a polynomial
with Newton polygon $\Del_5$ cannot have cusps.\end{enumerate}
\end{lemma}

{\bf Proof}. The statements come from a direct computation, and we
explain only (vi). Indeed, a polynomial with Newton polygon
$\Del_5$ defines a plane cubic which admits two tangent lines
intersecting at some point on a curve, which is impossible for a
cuspidal cubic by Pl\"ucker formulas. \proofend

\subsection{Amoebas and tropicalizations of 1-cuspidal curves}
The dimension of the stratum of curves with a cusp in $\Lam(\Del)$
is $r=|\Z^2\cap\Del|-3$. Pick $r$ distinct generic points
$\bx_1,...,\bx_r\in\R^2$ with rational coordinates, and points
$\bp_1,...,\bp_r\in(\K^*)^2$ such that $\val(\bp_i)=\bx_i$ and the
exponents of $t$ in the coordinates of $\bp_i$ are rational,
$i=1,...,r$. Then there are finitely many $1$-cuspidal curves in
$\Lam_\K(\Del)$ passing through $\bp_1,...,\bp_r$, and their
amoebas pass through $\bx_1,...,\bx_r$.

Observe that the coefficients of a polynomial $f\in\K[x,y]$,
defining a cuspidal curve in $\Lam_\K(\Del)$, passing through
$\bp_1,...,\bp_r$, are Puiseux series with rational exponents of
$t$. A parameter change $t\mapsto t^M$ with a suitable natural $M$
makes all these exponents integral, and the convex piece-wise
linear function $\nu_f:\Del\to\R$ integral-valued at integral
points. Through the rest of the proof we keep these assumptions.

We claim that the images of $1$-cuspidal curves in
$\Lam_\K(\Del)$, passing through $\bp_1,...,\bp_r$, are
$1$-cuspidal amoebas, passing through $\bx_1,...,\bx_r$.

Indeed, an amoeba $A_f$ passing through $\bx_1,...,\bx_r$ must
satisfy
$$|\Z^2\cap\Del|-1\ge\rk(A_f)\ge|\Z^2\cap\Del|-3\ .$$ 

If $\rk(A_f)=|\Z^2\cap\Del|-1$, then the subdivision $S_f$ of
$\Del$ consists of triangles of area $1/2$ which all are
$\Aff(\Z^2)$-equivalent to $\conv\{(0,0),(1,0),(0,1)\}$, and thus,
the tropicalization ${\cal T}(f)=(f_1,...,f_N)$ defines curves
$C_i\subset\Tor(\Del_i)$ which are non-singular and cross
$\Tor(\partial\Del_i)$ transversally. Thereby
$\{f=0\}\subset\Tor_\K(\Del)$ is non-singular.

If $\rk(A_f)=|\Z^2\cap\Del|-2$, then the subdivision $S_f$ of
$\Del$ contains either a parallelogram of area $1$, or a triangle
with edges of length $1$ and one interior integral point, or an
edge of length $2$, common for two triangles of area $1$. Then, by
Lemmas \ref{l7} and \ref{l13}(i), $f$ defines a curve with at most
node as singularity.

If, $\rk(A_f)=r=|\Z^2\cap\Del|-3$, then $A_f$ is either
$1$-cuspidal, or nodal. The subdivision $S_f$ for a nodal amoeba
or rank $r$, which is not $1$-cuspidal, has (besides triangles of
area $1/2$) either a parallelogram of area $2$ with edges of
length $1$ (cf. Figure \ref{f3}(e)) or two fragments of the
following three kinds: triangle with edges of length $1$ and one
interior integral point, parallelogram of area $1$, two triangles
of area $1$ with a common edge of length $2$. By Lemmas \ref{l7},
\ref{l13}(i), \ref{l16}(vi), a polynomial $f$, defining a curve in
$\Lam_\K(\Del)$ with a cusp, cannot have tropicalization
associated to a subdivision dual to a nodal but non-$1$-cuspidal
amoeba.

Furthermore, the polynomials in ${\cal T}(f)$ corresponding to the
polygon in $S_f$, which is $\Aff(\Z^2)$-equivalent to that shown
in Figure \ref{f3}(a,b,c,d), must be as described in Lemma
\ref{l16}(i-v), and the polynomials in ${\cal T}(f)$ corresponding
to the triangles without interior integral points must be as
described in Lemma \ref{l7}. If $S_f$ contains a polygon
$\Aff(\Z^2)$-equivalent to one of these shown in Figure
\ref{f3}(c,d), there is an edge $\sig=\Del_k\cap\Del$ and a point
$z\in\Tor(\sig)\cap C_k\cap C_l$ such that
$(C_k\cdot\Tor(sig))_z=(C_l\cdot\Tor(\sig))_z\ge 2$, and thus, we
can construct a $z$-refinement of the tropicalization of $f$ as
explained in section \ref{sec10}. Possible refinements of the
tropicalization of $f$ are described in the following statements.

\begin{figure}
\setlength{\unitlength}{1cm}
\begin{picture}(13,8)(0,8)
\thicklines\put(1,9){\vector(0,1){6}}\put(8,9){\vector(0,1){7}}
\put(1,12){\vector(1,0){5}}\put(8,12){\vector(1,0){4}}
\put(1,11){\line(3,1){3}} \put(1,13){\line(3,-1){3}}
\put(8,11){\line(2,1){2}} \put(8,14){\line(1,-1){2}}\thinlines
\dashline{0.2}(1,10)(5.5,10)\dashline{0.2}(1,11)(5.5,11)
\dashline{0.2}(1,13)(5.5,13)\dashline{0.2}(1,14)(5.5,14)
\dashline{0.2}(8,10)(11.5,10)\dashline{0.2}(8,11)(11.5,11)
\dashline{0.2}(8,13)(11.5,13)\dashline{0.2}(8,14)(11.5,14)
\dashline{0.2}(8,15)(11.5,15)\dashline{0.2}(2,9.5)(2,14.5)
\dashline{0.2}(3,9.5)(3,14.5)\dashline{0.2}(4,9.5)(4,14.5)
\dashline{0.2}(5,9.5)(5,14.5)\dashline{0.2}(9,9.5)(9,15.5)
\dashline{0.2}(10,9.5)(10,15.5)\dashline{0.2}(11,9.5)(11,15.5)
\put(3.5,8){{\rm (a)}} \put(10,8){{\rm (b)}}
\end{picture}
\caption{Newton polygons for deformation patterns}\label{f2}
\end{figure}
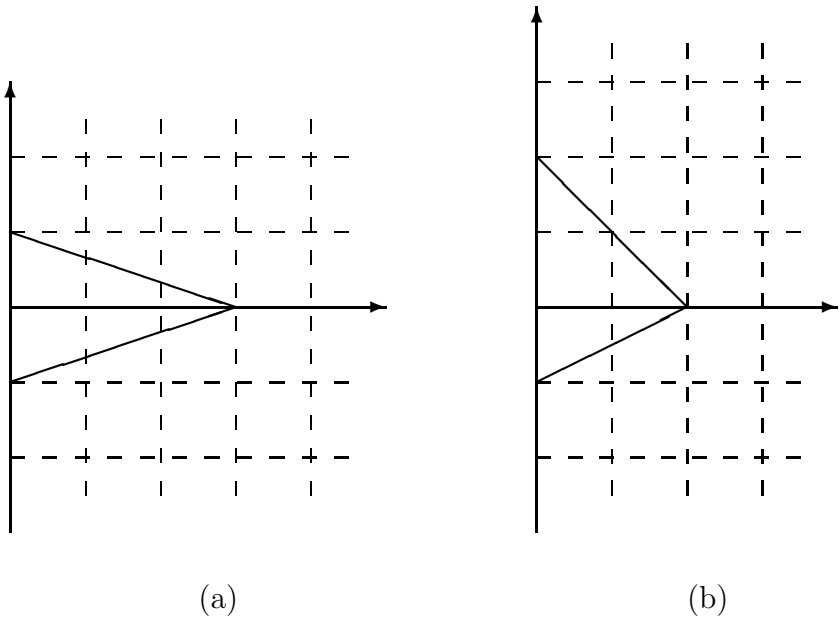

\begin{lemma}\label{l13} In the above notation and definitions,
let $z\in\Tor(\sig)\cap C_k\cap C_l$.
\begin{enumerate}
\item[(i)] If $C_k$ and $C_l$ are non-singular at $z$,
$(C_k\cdot\Tor(\sig))_z=3$, and $z$ bears a cusp $A_2$, then
${\cal T}_z(f)$ consist of one polynomial with Newton polygon
$\Del_z$ as shown in Figure \ref{f2}(a), which defines an elliptic
curve with one cusp in $\Tor(\Del_z)$.
\item[(ii)] If $C_k$ has a node at $z$, $C_l$ is non-singular at $z$,
$(C_k\cdot\Tor(\sig))_z=2$, and $z$ bears a cusp $A_2$, then
${\cal T}_z(f)$ consist of one polynomial with Newton polygon
$\Del_z$ as shown in Figure \ref{f2}(b), which defines a rational
curve with one cusp in $\Tor(\Del_z)$.
\end{enumerate}
\end{lemma}

{\bf Proof}. We show that ${\cal T}_z(f)$ in each case consists of
one polynomial, and describe these polynomials in Lemma \ref{l14}
below.

The proof proceeds in the same way as the proof of Lemma
\ref{l113}. As an example, we treat the situation (ii).

Besides the vertices of $\Del_z$ only the points $(0,0)$ and
$(0,1)$ may serve as vertices of the subdivision $S_{f^{(z)}}$.
Notice that the possible intersections of the curves defined by
${\cal T}_z(f)$ with $\Tor([(0,0),(2,0)])$ are transversal, since
the coefficient of $x$ is zero by construction. Hence $z$ may bear
the only singular points, coming from singularities in $(\C^*)^2$
curves defined by polynomials with Newton polygon
\mbox{$\conv\{(0,0),(2,0),(0,2)\}$} or
\mbox{$\conv\{(0,-1),(0,1),(2,0)\}$}, which are at most nodes.
\proofend

\begin{lemma}\label{l14}
\begin{enumerate}
\item[(i)] There exist exactly two polynomials with Newton polygon $\Del'$, as shown in Figure
\ref{f2}(a), which have prescribed coefficients at the vertices of
$\Del'$, zero coefficient of $x^2$, and define curves in
$\Tor(\Del')$ with a cusp as its only singularity;
\item[(ii)] there exist exactly three polynomials with Newton polygon $\Del''$, as shown in Figure
\ref{f2}(b), which have prescribed coefficients at the vertices of
$\Del''$, zero coefficient of $x^2$, and define curves in
$\Tor(\Del'')$ with a cusp as its only singularity.
\end{enumerate}
\end{lemma}

{\bf Proof}. The statement results from a direct computation. As
an example, we consider the second case.

After a suitable coordinate change, we reduce the question to the
study of polynomials
$$F(x,y)=y^3+yx^2+1+ay+by^2,\quad a,b\in\C\ .$$
The system $F(x,y)=F_x(x,y)=F_y(x,y)=0$ reduces in $\C^2$ to the
system
$$x=0,\quad 3y^2+2by+a=0,\quad y^3+by^2+ay+1=0\ ,$$ which
must have a solution of multiplicity $2$, that is
$y^3+by^2+ay+1=(y+\alp)^3$, $\alp^3=1$, and the statement follows.
\proofend

\subsection{Restoring a cuspidal curve out of tropical
data} We proceed along the argument of section \ref{sec11}.

Denote by ${\cal Q}_\Del(A_2)$ the set of quadruples $(A,S,F,R)$,
where
\begin{itemize}\item $A\in{\cal A}(\Del)$ is a 1-cuspidal amoeba,
$S:\ \Del=\Del_1,...,\Del_N$ is a subdivision of $\Del$ dual to
$A$, and $F$, $R$ are collections of the following polynomials in
$\C[x,y]$ which together are defined up to multiplication by the
same non-zero (complex) constant;
\item $F=(f_1,...,f_N)$, $f_i$ is a polynomial with
Newton polygon $\Del_i$, $i=1,...,N$, such that, if $\Del_i$ is a
triangle without interior integral points, then $f_i$ defines a
rational curve in $\Tor(\Del_i)$ as described in Lemma \ref{l7},
if $\Del_i$ is a triangle $\Aff(\Z^2)$-equivalent to that shown in
Figure \ref{f3}(b,c), then $f_i$ defines a curve in $\Tor(\Del_i)$
as described in Lemma \ref{l16}(ii,iii), if $\Del_i$ is a
quadrangle $\Aff(\Z^2)$-equivalent to that shown in Figure
\ref{f3}(a,d), then $f_i$ defines a curve in $\Tor(\Del_i)$ as
described in Lemma \ref{l16}(iv,v), and, at last, for any common
edge $\sig=\Del_i\cap\Del_j$, the truncations $f_i^\sig$ and
$f_j^\sig$ coincide; \item $R$ is a collection of deformation
patterns compatible with $F$ as defined in Remarks \ref{r4} and
\ref{r5}.
\end{itemize}

Let $(A,S,F,R)\in{\cal Q}_\Del(A_2)$, and $\bx_1,...,\bx_r\in A$.
These points lie on $r$ distinct edges of $A$ which correspond to
some $r$ edges of $S$, and they impose conditions (\ref{e144}) on
the values of the function $\nu:\Del\to\R$. Since
$\bx_1,...,\bx_r$ are generic, system (\ref{e144}) is independent.
A quadrangle which may appear in $S$ imposes one linear condition
on the values of $\nu$ at its vertices, which, for the case of the
shape shown in Figure \ref{f3}(a), reads (up to $\Aff(\Z^2)$
action)
$$3\nu(1,0)+\nu(1,2)=2\nu(0,1)+2\nu(2,0)\ ,$$ and, for the case of the
shape shown in Figure \ref{f3}(d), reads
$$\nu(0,0)+2\nu(1,1)=2\nu(0,1)+\nu(2,0)\ .$$ The latter condition together with
(\ref{e144}) form a system of $|V(S)|-1$ independent equations,
which determines the values of $\nu$ at the vertices of $S$
uniquely up to a shift.

We look for a polynomial $f\in\K[x,y]$ defining a curve in
$\Sig_\Del(A_2)$ in the form (\ref{e63}). Similarly to Step 2 of
section \ref{sec11}, the conditions $f(\bp_1)=...=f(\bp_r)=0$
transform into a system of equations (\ref{e146}). In the case of
a triangular subdivision $S$, the latter system determines the
coefficients of $f_1,...,f_N$ at $V(S)$ as well as the truncations
of $f_1,...,f_N$ on the edges $\sig_1,...,\sig_r$ uniquely up to
proportionality. If $S$ contains a quadrangle, then we supply
system of equations (\ref{e146}) with equation (\ref{e66}) or
(\ref{e67}), and the system obtained produces
$w(A,\bx_1,...,\bx_r)$ collections of the coefficients of
$f_1,...,f_N$ at $V(S)$ and truncations to $\sig_1,...\sig_r$ (up
to proportionality). Finally, taking into account Lemmas \ref{l16}
and \ref{l14}, we decide that given an 1-cuspidal amoeba
$A\in{\cal A}(\Del)$, $A\supset\{\bx_1,...,\bx_r\}$, there are
$W(A,\bx_1,...,\bx_r)\prod_{i=1}^r|\sig_i|^{-1}$ quadruples
$(A,S,F,R)\in{\cal Q}_\Del(A_2)$ which may serve as
tropicalizations of polynomials $f\in\K[x,y]$ with Newton polygon
$\Del$, defining curves $C\in\Sig_\Del(A_2)$ passing through
$\bp_1,...,\bp_r$.

The proof of Theorem \ref{t6} is completed with

\begin{lemma}\label{l20} In the above notation, given the points $\bx_1,...,\bx_r\in\R^2$,
$\bp_1,...,\bp_r\in(\K^*)^2$, and a compatible $(A,S,F,R)\in{\cal
Q}_\Del(A_2)$, the polynomials $f\in\K[x,y]$ with Newton polygon
$\Del$, which tropicalize into $(A,S,F,R)$ define exactly
$\prod_{i=1}^r|\sig_i|$ curves in $\Lam_\K(\Del)$ having one cusp
and passing through $\bp_1,...,\bp_r$.
\end{lemma}

The proof of Lemma \ref{l20} is completely similar to the proof of
Lemma \ref{l18} (section \ref{sec300}) and we omit it.

\section{Patchworking singular algebraic curves}\label{sec8}

\subsection{Initial data for patchworking}\label{sec1} Let $\Del\subset\R^2$ be a
non-degenerate convex lattice polygon, $S:\
\Del=\Del_1\cup...\cup\Del_N$ its subdivision into convex lattice
polygons, defined by a convex piece-wise linear function
$\nu:\Del\to\R$ such that $\nu(\Z^2)\subset\Z$.

Let $a_{ij}\in\C$, $(i,j)\in\Del\cap\Z^2$, be such that $a_{ij}\ne
0$ for each vertex $(i,j)$ of all the polygons
$\Del_1,...,\Del_N$. Then we define polynomials
$$f_k(x,y)=\sum_{(i,j)\in\Del_k\cap\Z^2}a_{ij}x^iy^j,\quad
k=1,...,N\ ,$$ and curves $C_k=\{f_k=0\}\subset\Tor(\Del_k)$,
$k=1,...,N$, on which we impose the following conditions.
\begin{itemize} \item[(A)] For any $k=1,...,N$, each multiple
component of $C_k$ (if it exists) is defined by a binomial; it
crosses any other component of $C_k$ transversally, only at
non-singular points, and not on $\Tor(\partial\Del_k)$.
\item[(B)] For any edge $\sig\subset\partial\Del$, $\sig\subset\Del_k$, $1\le k\le N$,
the curve $C_k$ is non-singular along $\Tor(\sig)$ and crosses
$\Tor(\sig)$ transversally. \item[(C)] If $\sig$ is an edge of
$\Del_k$, $1\le k\le N$, and $z\in\Tor(\sig)\cap C_k$ is an
isolated singular point of $C_k$, then the germ $(C_k,z)$ is
topologically equivalent to $(y'')^{m(k,z)}+(x'')^m=0$, in local
coordinates $x'',y''$ with $y''$-axis coinciding with
$\Tor(\sig)$.
\end{itemize}

Now we introduce additional polynomials which will play the role
of deformation patterns, as defined in Remarks \ref{r4} and
\ref{r5}.

Consider all the triples $(k,\sig,z)$, where $1\le k\le N$,
$\sig\not\subset\partial\Del$ is an edge of $\Del_k$,
$z\in\Tor(\sig)\cap C_k$ and $(C_k\cdot\Tor(\sig))_z=m\ge 2$. Then
introduce the equivalence of triples: (i)
$(k,\sig,z)\sim(l,\sig,z)$ if $\sig=\Del_k\cap\Del_l$, and (ii)
$(k,\sig,z)\sim(k,\sig',z')$ if $\sig,\sig'$ are parallel sides of
$\Del_k$ and $z,z'$ belong to the same multiple component of
$C_k$. The transitive extension of this equivalence distributes
the triples into disjoint classes. We denote the set of
equivalence classes by $\Pi$. In fact, a pair of points $z,z'$
from equivalent triples $(k,\sig,z)$, $(l,\sig',z')$ determines an
element of $\Pi$ uniquely, and we write simply $(z,z')\in\Pi$.

To any element of $\Pi$ we assign a deformation pattern. Namely,
in any class there are exactly two triples $(k,\sig,z)$,
$(l,\sig',z')$ with coinciding or parallel edges $\sig$, $\sig'$,
and isolated singular (or non-singular) points $z$, $z'$ of the
curves $C_k$, $C_l$, respectively. In some local coordinates in
neighborhoods of $z$ and $z'$ as required in the above property
(C), the curves $C_k$ and $C_l$ are defined by
$$\sum_{i\cdot m(k,z)+jm\ge
m\cdot m(k,z)}\alp_{ij}x^iy^j=0,\quad\sum_{i\cdot m(l,z')+jm\ge
m\cdot m(l,z')}\bet_{ij}x^iy^j=0\ ,$$ respectively, with
$\alp_{m0}=\bet_{m0}$, and non-degenerate homogeneous polynomials
$$\varphi^{(k)}_z(x,y)=\sum_{i\cdot m(k,z)+jm=m\cdot m(k,z)}\alp_{ij}x^iy^j,\quad
\varphi^{(l)}_{z'}(x,y)=\sum_{i\cdot m(l,z')+jm=m\cdot
m(l,z')}\bet_{ij}x^iy^j\ .$$ A deformation pattern attached to the
chosen class of triples is a curve
$C_{z,z'}\subset\Tor(\Del_{z,z'})$,
$\Del_{z,z'}=\conv\{(m,0),(0,m(k,z)),(0,-m(l,z'))\}$, defined by a
polynomial $F_{z,z'}(x,y)$ with Newton triangle $\Del_{z,z'}$ and
truncations $\varphi^{(k)}_z(x,y)$, $\varphi^{(l)}_{z'}(x,y^{-1})$
on the edges $[(m,0),(0,m(k,z)]$, $[(m,0),(0,-m(l,z'))]$,
respectively.

\subsection{Transversality}\label{sec3} Transversality of equisingular strata provides
sufficient conditions for the patchworking (cf. \cite{Sh2,Sh3}).

Let $\ks$ be a topological or (contact) analytic equivalence of
isolated planar curve singular points. We intend to define the
$\ks$-transversality for triples $(\Del_k,\Del_k^-,C_k)$, $1\le
k\le N$, where $\Del_k^-$ is a connected (or empty) union of some
edges of $\Del_k$, and for deformation patterns.

Pick a triple $(\Del_k,\Del_k^-,C_k)$, $1\le k\le N$.

Denote by $\Sing^{\text{\rm iso}}(C_k)$ the set of isolated
singular points of $C_k$. If $z\in\Sing^{\text{\rm
iso}}(C_k)\cap(\C^*)^2$, denote by $M^\ks(C_k,z)$ the germ at
$C_k$ of the $\ks$-equisingular stratum of $(C_k,z)$ in
$\Lam(\Del_k)$. The (projective) Zariski tangent space to
$M^\ks(C_k,z)$ at $C_k$ is formed by the curves
$\{g=0\}\in\Lam(\Del_k)$, with $g\in
I^\ks(C_k,z)\subset\ko_{\Tor(\Del_k),z}$, where $I^\ks(C_k,z)$ is
the equisingular ideal or the Tjurina ideal (see \cite{DH,W}),
according to whether $\ks$ is the topological or analytic
equivalence.

Let $z\in C_k\cap\Tor(\sig)$ be a non-singular or singular
isolated point of $C_k$, where $\sig$ is an edge of $\Del_k$, and
let $x'',y''$ be local coordinates in a neighborhood of $z$ in
$\Tor(\Del_k)$ as introduced in section \ref{sec1}. The ideals
$$I^{\sqh}_0(C_k,z)=\Big\{g\in\ko_{\Tor(\Del_k),z}\ :\ g=\sum_{i\cdot m(k,z)+jm\ge m\cdot m(k,z)}
b_{ij}(x'')^i(y'')^j\Big\}\ ,$$
$$I^{\sqh}(C_k,z)=I^{\sqh}_0(C_k,z)+\langle \frac{\partial
f''}{\partial x''}\rangle\ $$ naturally define the linear
subsystems $M^{\sqh}_0(C_k,z)$, $M^{\sqh}(C_k,z)$ in
$\Lam(\Del_k)$, respectively.

Let $C^{\text{\rm red}}_k$ be the reduction of the curve $C_k$.
Let $z\in\Sing^{\text{\rm iso}}(C^{\text{\rm
red}}_k)\backslash\Sing^{\text{\rm iso}}(C_k)$. Then $z$ is an
intersection point of two distinct components $\{g'=0\},\{g''=0\}$
of $C_k$ having multiplicities $m',m''$, respectively, with
$m'+m''>2$. Denote by $M^{eg}(C_k,z)$ the closure of the germ at
$C_k$ of the family of curves $C\in\Lam(\Del_k)$, having $m'm''$
nodes in a neighborhood of $z$.

\begin{lemma}\label{l21}
(i) The (projective) Zariski tangent space to $M^{eg}(C_k,z)$ at
$C_k$ is formed by the curves $\{g=0\}$, $g\in\Lam(\Del_k)$, with
$g\in I^{eg}(C_k,z):=
\langle(g')^{m'},(g'')^{m''}\rangle\subset\ko_{\Tor(\Del_k),z}$.

(ii) Let $\Del_k$ be a parallelogram with a pair of non-parallel
edges $\sig_1,\sig_2$, the curve $C_k$ given by $\{f_k=0\}$, where
$f_k$ is a product of a monomial and binomials. Then the germ
$M^{eg}(C_k)=\bigcap_zM^{eg}(C_k,z)$, where $z$ runs over all
intersection points of distinct components of $C_k$, is smooth of
codimension $\Area(\Del_k)$ in $\Lam(\Del_k)$, and intersects
transversally with the space of curves, defined by polynomials $f$
with Newton polygon $\Del_k$ such that $f^{\sig_i}=f_k^{\sig_i}$,
$i=1,2$.
\end{lemma}

{\bf Proof}. (i) In a neighborhood of $z$, the curves $C\in
M^{eg}(C_k,z)$ are unions of $m'+m''$ discs (counting
multiplicities), and are represented by equations
$((g')^{m'}+g'_1)((g'')^{m''}+g''_1)=0$ with $||g'_1||,||g''_1||$
sufficiently small; thus, the claim follows.

(ii) Observing that
$|\{f_k^{\sig_1}=0\}\cap\{f_k^{\sig_2}=0\}\cap(\C^*)^2|=\Area(\Del_k)|\sig_1|^{-1}|\sig_2|^{-1}$,
we derive the required statement, when showing that $f_k$ is the
only polynomial with Newton polygon $\Del_k$, the fixed
truncations on $\sig_1,\sig_2$, and belonging to the ideal
$\langle f_k^{\sig_1},\ f_k^{\sig_2}\rangle_w\subset{\cal
O}_{\C^2,w}$, for any point
$w\in\{f_k^{\sig_1}=0\}\cap\{f_k^{\sig_2}=0\}\cap(\C^*)^2$. The
latter claim immediately follows from B\'ezout's theorem.
\proofend

\begin{definition}\label{d1} {\rm In the above notation, let $\Del_k^+$ be
the union of the edges $\sig$ of $\Del_k$ such that
$\sig\not\subset\Del_k^-$. The triad $(\Del_k,\Del_k^-,C_k)$ is
called $\ks$-transversal, if all the germs
$$\begin{cases}M^{\ks}(C_k,z),\quad &z\in\Sing^{\text{\rm iso}}(C_k)\cap(\C^*)^2\ ,\\ M^{eg}(C_k,z),\quad &z\in
\Sing^{\text{\rm iso}}(C^{\text{\rm
red}}_k)\backslash\Sing^{\text{\rm iso}}(C_k)\ ,\\
M^{\sqh}_0(C_k,z),\quad &z\in C_k\cap\Tor(\Del_k^-)\
\text{is not a non-isolated singular point}\ ,\\
M^{\sqh}(C_k,z),\quad &z\in C_k\cap\Tor(\Del_k^+)\ \text{is not a
non-isolated singular point}
\end{cases}$$
are smooth of expected dimension and intersect transversally in
$\Lam(\Del_k)$.}
\end{definition}

\begin{definition}\label{d2}
{\rm A deformation pattern $C_{z,z'}\subset\Del_{z,z'}$, is called
$\ks$-transversal, if the triad
$(\Del_{z,z'},\Del_{z,z'}^-,C_{z,z'})$ is $\ks$-transversal, where
$\Del_{z,z'}^-$ is the union of the non-vertical edges of
$\Del_{z,z'}$.}
\end{definition}

\begin{lemma}\label{l1} In the above notation,
the triad $(\Del_k,\Del_k^-,C_k)$ is $\ks$-transversal if
\begin{equation}
H^1(\Tor(\Del_k),\ko_{\Tor(\Del_k)}(C_k)\otimes\kj_{Z_k})=0\
,\label{e1}
\end{equation}
where $\kj_{Z_k}\subset\ko_{\Tor(\Del_k)}$ is the ideal sheaf of
the zero-dimensional scheme $Z_k\subset\Tor(\Del_k)$, defined at
the points $z\in C_k$ mentioned in Definition \ref{d1} by the
ideals $I^{\ks}(C_k,z)$, $I^{eg}(C_k,z)$, $I^{sqh}_0(C_k,z)$,
$I^{sqh}(C_k,z)$, respectively.
\end{lemma}

The statement immediately follows from the cohomology
interpretation of transversality.

Following \cite{Sh2,Sh3}, we provide an explicit numerical
criterion for the $h^1$-vanishing (\ref{e1}). To formulate it, we
use topological invariants $b(C,\xi)$, $\widetilde b(C,\xi)$
defined for a curve $C$ and its local branch $\xi$, and the
Tjurina number $\tau(C,z)$, equal to the codimension of the
Tjurina ideal in the local ring of an ambient surface. The
complete definition of the invariants $b$ and $\widetilde b$ can
be found in \cite{Sh2}, section 4.1, or in \cite{Sh3}, section 4,
Definition 1. We only recall it for a few cases. If $C$ has a
node, then $b(C,\xi)=0$ for both branches; if $C$ has a cusp, then
$b(C,\xi)=1$; if $C$ is locally given by $\{x^{pr}+y^{qr}=0\}$,
$(p,q)=1$, then $\widetilde b(C,\xi)=p+q-1$.

\begin{lemma}\label{l2} (i) The $\ks$-transversality of a triad $(\Del_k,\Del_k^-,C_k)$
persists if $\Del_k^-$ contains at most two edges, and one removes
from $\Del_k^-$ edges of length $1$.

(ii) If $C_k$ is irreducible, then the triad
$(\Del_k,\Del_k^-,C_k)$ is transversal with respect to the
topological equivalence of singular points, provided,
$${\sum}'b(C_k,\xi)+{\sum}''\widetilde b(C_k,{\cal Q})+
{\sum}'''((C_k\cdot\Tor(\sig))_z-\eps)<
\sum_{\sig\subset\partial\Del_k}(C_k\cdot\Tor(\sig))\ ,$$ where
$\sum'$ ranges over all local branches $\xi$ of $C_k$, centered at
the points $z\in\Sing(C_k)\cap(\C^*)^2$, $\sum''$ ranges over all
local branches ${\cal Q}$ of $C_k$, centered at the points
$z\in\Sing(C_k)\cap\Tor(\partial\Del_k)$, and $\sum'''$ ranges
over all non-singular points $z$ of $C_k$ on
$\Tor(\partial\Del_k)$ with $\eps=0$ if $\sig\subset\Del_k^-$ and
$\eps=1$ otherwise.

(iii) If $C_k$ is irreducible, then the triad
$(\Del_k,\Del_k^-,C_k)$ is transversal with respect to the
analytic equivalence of singular points, provided,
$${\sum}'(\tau(C_k,z)-1)+{\sum}''\widetilde b(C_k,{\cal Q})+
{\sum}'''((C_k\cdot\Tor(\sig))_z-\eps)<
\sum_{\sig\subset\partial\Del_k}(C_k\cdot\Tor(\sig))\ ,$$ where
$\sum'$ ranges over all the points $z\in\Sing(C_k)\cap(\C^*)^2$,
and $\sum''$, $\sum'''$ are as above.

(iv) If $C_k$ is reduced and reducible, then the triad
$(\Del_k,\Del_k^-,C_k)$ is transversal with respect to the
topological equivalence of singular points, provided, for any
irreducible component $C$ of $C_k$,
\begin{equation}{\sum}'b(C_k,\xi)+{\sum}''\widetilde b(C_k,{\cal Q})+
{\sum}'''((C\cdot\Tor(\sig))_z-\eps)<
\sum_{\sig\subset\partial\Del_k}(C\cdot\Tor(\sig))\
,\label{e1n}\end{equation} where $\sum'$ ranges over all local
branches $\xi$ of $C$, centered at the points
$z\in\Sing(C_k)\cap(\C^*)^2$, $\sum''$ ranges over all local
branches ${\cal Q}$ of $C$, centered at the points
$z\in\Sing(C_k)\cap\Tor(\partial\Del_k)$, and $\sum'''$ ranges
over all non-singular points $z$ of $C_k$ on $C\cap\Tor(\sig)$,
$\sig\subset\partial\Del_k$, with $\eps=0$ if
$\sig\subset\Del_k^-$ and $\eps=1$ otherwise.

(v) If $C_k$ is non-reduced, then the triad
$(\Del_k,\Del_k^-,C_k)$ is transversal with respect to the
topological equivalence of singular points, provided that any
component of $C_k$, which is not defined by a binomial, satisfies
(\ref{e1n}), and any component of $C_k$ defined by a binomial
crosses $\Tor(\Del_k^-)$ at most in one point and crosses the
reduced union of all the other component of $C_k$ transversally in
only non-singular points.
\end{lemma}

\begin{lemma}\label{l3}
In the notation of section \ref{sec1} and Definition \ref{d2},
\begin{enumerate}\item[(i)] an irreducible deformation pattern $C_{z,z'}$ is transversal with
respect to the topological equivalence of singular points if
$$\sum_{w\in\Sing(C_{z,z'})\cap\C^2}b(C_{z,z'},w)<\#(\Z\cap(-m(l,z'),m(k,z)))+\eps_0\ ,$$
and is transversal with respect to the analytic equivalence of
singular points if
$$\sum_{w\in\Sing(C_{z,z'})\cap\C^2}(\tau(C_{z,z'},w)-1)<\#(\Z\cap(-m(l,z'),m(k,z)))+\eps_0\
,$$ where $\eps_0$ is the number of edges of length $1$ in
$\Del_{z,z'}^-$;
\item[(ii)] a reducible deformation pattern $C_{z,z'}$ is transversal
with respect to the topological equivalence of singular points if,
for any irreducible component $C$ of $C_{z,z'}$,
$$\sum_{\xi}b(C_{z,z'},\xi)<(C\cdot\Tor(\partial\Del_{z,z'}))-(C\cdot\Tor(\Del_{z,z'}^-))+\eps_0(C)\ ,$$
where $\xi$ ranges on all local branches of $C$ centered at
$\Sing(C_{z,z'})\cap\C^2$, and $\eps_0(C)$ is the number of edges
of length $1$ in $\Del_{z,z'}^-$.\end{enumerate}
\end{lemma}

{\bf Proof}. Both Lemmas \ref{l2} and \ref{l3} are slightly
modified particular cases of Theorem 4.1 in \cite{Sh2}, and we
shall explain only the modifications.

To obtain $\eps_0$, or, more generally, to remove edges of length
$1$ from $\Del_k^-$, we notice that the $\ks$-transversality for
the triad $(\Del_k,\widetilde\Del_k^-,C_k)$ means that in the
space ${\cal P}(\Del_k)$ of polynomials with Newton polygon
$\Del_k$, the corresponding $\ks$-equisingular stratum is smooth
and transversally intersects with the (affine) subspace of
polynomials having fixed coefficients at the integral points in
$\widetilde\Del_k^-$. The action of $(\C^*)^3$ on ${\cal
P}(\Del_k)$ defined as $(\lam_0,\lam_1,\lam_2)\cdot
F(x,y)=\lam_0F(\lam_1x,\lam_2y)$, arbitrarily varies the
coefficients, corresponding to integral points in the edges of
length $1$ in $\Del_k^-\backslash\widetilde\Del_k^-$, whereas the
coefficients at the integral points in $\widetilde\Del_k^-$ stay
fixed. Since the considered $\ks$-equisingular stratum is
invariant with respect to this action, we conclude that it
transversally intersects the subspace of polynomials having fixed
coefficients at the integral points in $\Del_k^-$.

In case (v) of Lemma \ref{l2}, one derives (\ref{e1}) by
successively eliminating the components of $C_k$, defined by
binomials, using the Horace method \cite{Hi}, and then applying
statement (iii) to the rest of the curve. Indeed, if $C'_k$ is the
union of the components of $C_k$, which are not defined by
binomials, then, for any component $C\subset C'_k$, and any local
branch $\xi$ of $C$, participating in the left-hand side of
(\ref{e1n}), we have respectively $b(C'_k,\xi)=b(C_k,\xi)$, or
$\widetilde b(C'_k,\xi)=\widetilde b(C_k,\xi)$. In turn, the
Horace method step is performed as follows. Let $C\subset C_k$ be
a component of $C_k$ defined by a binomial $f$, and $C''_k$ be the
curve on $\Tor(\Del_k)$ defined by the polynomial $f''_k=f_k/f$.
We have an exact sequence of ideal sheaves on $\Tor(\Del_k)$
$$0\to\kj_{Z_k:C}(C''_k)\to\kj_{Z_k}(C_k)\to\kj_{Z_k\cap C}(C_k)\to 0\ .$$ Here the ideal
sheaf $\kj_{Z_k:C}(C''_k)$ can be considered on the surface
$\Tor(\Del''_k)$, where $\Del''_k$ is the Newton polygon of the
polynomial $f''_k$. It is easy to see that $\Del''_k$ can be
obtained from $\Del_k$, when shrinking two opposite parallel edges
of $\Del_k$ (the edges parallel to the Newton segment of the
binomial $f$) by a unit length segment. The scheme $Z_k:C$ is the
zero-dimensional scheme, which can be defined for the curve
$C''_k$ on $\Tor(\Del''_k)$ along the formulation of Lemma
\ref{l1}, in the same way as the scheme $Z_k$ for $C_k$. In turn
the sheaf $\kj_{Z_k\cap C}(C_k)$ can be considered on the
non-singular rational curve $C$, and one can easily check that, by
construction and hypotheses of Lemma \ref{l2}, $$\deg Z_k\cap C\le
\deg\ko_C(C_k)+1\quad\Longrightarrow\quad\deg\kj_{Z_k\cap
C}(C_k)\ge-1>-2=2g(C)-2\ .$$ Hence $h^1(C,\kj_{Z_k\cap
C}(C_k))=0$, and thus, (\ref{e1}) reduces to
$h^1(\Tor(\Del''_k),\kj_{Z_k:C}(C''_k))=0$, which is the induction
assumption. \proofend

\subsection{Patchworking theorem}\label{sec4} Let us be given the
data introduced in section \ref{sec1}, i.e., subdivision
\mbox{$S:\ \Del=\Del_1\cup...\cup\Del_N$}, induced by a function
$\nu:\Del\to\R$, amoeba $A$, polynomials $f_1,...,f_N$, and
deformation patterns defined by polynomials $f_{z,\wz}$. Let
${\cal G}$ be the set of orientations of the amoeba $A$ (as a
graph), which have no oriented cycles and obey the following
requirements. For $\Gam\in{\cal G}$, denote by $\Del_k^-(\Gam)$
the union of those edges of $\Del_k$ which correspond to arcs of
$A$, which are $\Gamma$-oriented inside $\Del_k$. We assume that
$\Del_k^-(\Gam)$ is connected for any $k=1,...,N$, and any two
arcs of $A$, having a common vertex and lying on a straight line,
are cooriented. Denote by $\Arc(\Gam)$ the set of ordered pairs
$(k,l)$, where $\Del_k$, $\Del_l$ have a common edge, and the
corresponding arc of $A$ is $\Gam$-oriented from $\Del_k$ to
$\Del_l$.

\begin{theorem}\label{t1} Under the assumptions of sections
\ref{sec1}, suppose that all the given deformation patterns are
$\ks$-transversal, and there is $\Gam\in{\cal G}$ such that every
triad $(\Del_k,\Del_k^-(\Gamma),C_k)$ is $\ks$-transversal,
$k=1,...,N$. Then there exists a polynomial $f\in\K[x,y]$ with
Newton polygon $\Del$, whose refined tropicalization consists of
the given data, $\nu$, $S$, $f_1,...,f_N$, and the given
deformation patterns, and which defines a family of reduced curves
$C^{(t)}\subset\Tor(\Del)$, $t\ne 0$, such that there is an
$\ks$-equivalent 1-to-1 correspondence between $\Sing(C^{(t)})$
and the disjoint union of \begin{itemize}\item the sets
$\Sing^{\text{\rm iso}}(C_k)\cap(\C^*)^2$, $k=1,...,N$, \item the
sets\ $\Sing(C_{z,\wz})\cap\C^2$, $\{z,\wz\}\in\Pi$, \item the set
of $\sum_{k=1}^N\sum_z\dim \ko_{\C^2,z}/I^{eg}(C_k,z)$ nodes,
where $z$ runs over $\Sing(C^{\text{\rm
red}}_k)\backslash\Sing^{\text{\rm iso}}(C_k)$,
$k=1,...,N$.\end{itemize}

Furthermore, take any set $B\subset V(S)$ such that, for each
$k=1,...,N$, either $|B\cap\Del_k|\le 3$, or
$B\cap\Del_k\subset\Del_k^-$. Then one has a family of polynomials
as above, which can be described by the relations
\begin{equation}f(x,y)=\sum_{(i,j)\in\Del}(a_{ij}+c_{ij})x^iy^jt^{\nu(i,j)}\ ,\label{e2n}\end{equation}
\begin{equation}\begin{cases}&c_{ij}=c_{ij}(t)\in\K,\quad c_{ij}(0)=0,\quad(i,j)\in\Del\
,\\
&c_{ij}(t)=\Phi^B_{ij}(\{c_{kl}(t),\ (k,l)\in B\},t),\quad
(i,j)\in\Del\cap\Z^2\backslash B\
,\end{cases}\label{e3n}\end{equation} with certain complex
analytic functions $\Phi^B_{ij}$, $(i,j)\in\Del\cap\Z^2\backslash
B$.
\end{theorem}

The proof is a routine adaptation of the proofs of similar
patchworking theorems in \cite{Sh2,Sh3}. We however present it
here for a curious reader.

{\bf Proof of Theorem \ref{t1}}.

{\it Step 1}. Fix $k=1,...,N$. The space
$H^0(\Tor(\Del_k),\ko_{\Tor(\Del_k)}(C_k))$ can naturally be
identified with the space of polynomials ${\cal P}(\Del_k)$. We
shall split this linear space into subspaces and choose specific
bases in them.

Introduce $${\cal P}(\Del_k^-)=\Span\{x^iy^j,\
(i,j)\in\Del_k^-\},\quad\overline{\cal
P}(\Del_k^-)=\Span\{x^iy^j,\ (i,j)\in\Del_k\backslash\Del_k^-\}\
.$$ The $\ks$-transversality of the triad $(\Del_k,\Del_k^-,C_k)$
yields, by Lemma \ref{l1}, the surjectivity of the map
$\pr_k:{\cal P}(\Del_k)\to H^0(Z'_k,\ko_{Z'_k})$, where $Z'_k$ is
the part of the zero-dimensional scheme $Z_k$, supported at points
on $\Tor(\partial\Del_k)$, and $Z_k$ is introduced in Lemma
\ref{l1}. Since $f_k\in\Ker(\pr_k)$ and $\Del_k^-$ is connected,
$$\Ker(\pr_k)=\Span\{f_k\}\oplus(\Ker(\pr_k)\cap\overline{\cal
P}(\Del_k^-))\ .$$

The variety germs $M^{\ks}(C_k,z)$ and $M^{eg}(C_k,z)$, listed in
Definition \ref{d1} and related to singular points $z$ of $C_k$ in
$(\C^*)^2$, naturally lift to germs at $f_k$ of varieties in
${\cal P}(\Del_k)$. The latter germs we extend up to the germs
$\widetilde M^{\ks}(C_k,z)$ and $\widetilde M^{eg}(C_k,z)$ at
$f_k$ of the corresponding equisingular varieties in the space
${\cal P}(\Del)$. Due to the $\ks$-transversality condition, the
intersection $M_k$ of all the germs $\widetilde M^{\ks}(C_k,z)$,
$\widetilde M^{eg}(C_k,z)$, corresponding to the singular points
of $C_k$ in $(\C^*)^2$, is smooth of expected codimension
$n_k=\sum_z\dim\ko_{\C^2,z}/I^{\ks}(C_k,z)+\sum_z\dim\ko_{\C^2,z}/I^{eg}(C_k,z)$
in ${\cal P}(\Del)$. Furthermore, $M_k$ intersects transversally
with $\Ker(\pr_k)\cap\overline{\cal P}(\Del_k^-)$ in ${\cal
P}(\Del)$. That is, in a neighborhood of $f_k$, this germ is given
by a system of analytic equations
\begin{equation}\Phi^{(k)}_1(F)=...=\Phi^{(k)}_{n_k}(F)=0\ ,\label{e48}\end{equation} where $F$ stands
for a variable polynomial in ${\cal P}(\Del)$, and there is a set
$B_k$ of $n_k$ linearly independent elements of
$\Ker(\pi_k)\cap\overline{\cal P}(\Del_k^-)$ such that
\begin{equation}
\det\left(\frac{\partial\{\Phi^{(k)}_i(F),\
i=1,...,n_k\}}{\partial B_k}\right)\bigg|_{F=F_k}\ne 0\
.\label{e3}
\end{equation} Moreover, the elements of $B_k$ can be chosen as
the disjoint union of sets $B_{k,z}$, $z\in\Sing(C^{\text{\rm
red}}_k)$ such that $B_{k,z}$ projects to a basis of the space
$\ko_{\C^2,z}/I^{\ks}(C_k,z)$ or $\ko_{\C^2,z}/I^{eg}(C_k,z)$,
according as $z\in\Sing^{\text{\rm iso}}(C_k)\cap(\C^*)^2$ or
$z\in\Sing(C^{\text{\rm red}}_k)\backslash\Sing^{\text{\rm
iso}}(C_k)$, and, in addition $B_{k,z}$ projects to zero in any
other such space corresponding to any point $z'\ne
z\in\Sing(C^{\text{\rm red}}_k)\cap(\C^*)^2$.

\medskip

{\it Step 2}. Under the hypotheses of Step 1, we choose a basis of
$H^0(Z'_k,\ko_{Z'_k})$ reflecting its splitting
\begin{equation}\bigoplus_{z\in C_k\cap\Tor(\Del_k^-)\backslash C_k^{nr}}
\ko_{\Tor(\Del_k),z}/I^{\sqh}_0(C_k,z)\oplus\bigoplus_{\renewcommand{\arraystretch}{0.6}
\begin{array}{c}
\scriptstyle{z\in C_k\cap\Tor(\Del_k^+)\backslash C_k^{nr}}\\
\scriptstyle{(C_k\cdot\Tor(\Del_k^+))_z\ge 2}
\end{array}}\ko_{\Tor(\Del_k),z}/I^{\sqh}(C_k,z)\
,\label{e4}\end{equation} where $C_k^{nr}$ stands for the union of
the multiple components of $C_k$.

For a point $z$, occurring in the above splitting, we have local
coordinates $x'',y''$, introduced in condition (C), section
\ref{sec1}, in which $I^{\sqh}_0(C_k,z)$ is generated by monomials
lying on or above the segment $[(0,m(k,z)),(m,0)]$. Then we can
choose the monomial basis
\begin{equation}(x'')^i(y'')^j,\quad m(k,z)i+mj<m\cdot m(k,z)\ ,\label{e68}\end{equation} for
$\ko_{\Tor(\Del_k),z}/I^{\sqh}_0(C_k,z)$, and the monomial basis
\begin{equation}(x'')^i(y'')^j,\quad m(k,z)i+mj<m\cdot m(k,z),
\quad (i,j)\ne(m-1,0)\ ,\label{e5}\end{equation} for
$\ko_{\Tor(\Del_k),z}/I^{\sqh}(C_k,z)$.

Take a point $z\in C_k\cap\Tor(\Del_k^+)\backslash C_k^{nr}$ such
that $(C_k\cdot\Tor(\Del_k^+))_z\ge 2$. In the local coordinates
$x'',y''$, the ideal $I^{\sqh}(C_k,z)$ is generated by monomials
lying on or above segment $[(0,m(k,z)),(m,0)]$ and by $\partial
f''_k/\partial x''$. We lift the monomial basis (\ref{e5}) to
polynomials $\pi_{ij}^{(k,z)}\in{\cal P}(\Del_k)$, which can be
chosen obeying the following restrictions:
\begin{itemize}
\item $\pi_{ij}^{(k,z)}$ vanishes in all summands of (\ref{e4})
corresponding to points different from $z$,
\item $\pi_{ij}^{(k,z)}$ belongs to ${\cal
P}(\Del_k\backslash\Del_k^-)$,
\item $\pi_{ij}^{(k,z)}$ with $j>0$ does not contain the
monomials $x^py^q$, $(p,q)\in\sig$, where $z\in
C_k\cap\Tor(\sig)$.
\end{itemize}

Similarly, if $z\in C_k\cap\Tor(\Del_k^-)$ such that
$(C_k\cdot\Tor(\Del_k^-))_z\ge 2$, then, in the local coordinates
$x'',y''$, the ideal $I^{\sqh}_0(C_k,z)$ is generated by monomials
lying on or above segment $[(0,m(k,z)),(m,0)]$, and we lift the
monomial basis (\ref{e68}) to polynomials
$\pi_{ij}^{(k,z)}\in{\cal P}(\Del_k)$, which can be chosen obeying
the following restrictions:
\begin{itemize}
\item $\pi_{ij}^{(k,z)}$ vanishes in all summands of (\ref{e4})
corresponding to points $\ne z$,
\item $\pi_{ij}^{(k,z)}$, $j>0$, belongs to ${\cal
P}(\Del_k\backslash\Del_k^-)$.
\end{itemize}

\medskip

{\it Step 3}. Let $\{(k,\sig,z),(l,\widetilde\sig,\wz)\}\in\Pi$,
$C_{z,\wz}\subset\Tor(\Del_{z,\wz})$ the given deformation
pattern, defined by a polynomial $f_{z,\wz}$ with Newton polygon
$\Del_{z,\wz}$. We can identify
$H^0(\Tor(\Del_{z,\wz}),\ko_{\Tor(\Del_{z,\wz})}(C_{z,\wz}))$ with
${\cal P}(\Del_{z,\wz})$. Assume that ${\cal P}(\Del_{z,\wz})$ is
embedded into some finite-dimensional linear space ${\cal V}$ of
polynomials. The $\ks$-transversality of the deformation pattern
$C_{z,\wz}$ means that the germ at $f_{z,\wz}$ of the
$\ks$-equisingular stratum in ${\cal V}$, corresponding to
$\Sing(C_{z,\wz})\cap\C^2$, is smooth of expected dimension (which
we denote by $n_{z,\wz}$), and is the intersection of smooth
analytic transverse hypersurfaces
$$\Phi^{(z,\wz)}_1(F)=0,\quad...\quad ,\ \Phi^{(z,\wz)}_{n_{z,\wz}}(F)=0,\quad
F\in{\cal V}\ ,$$ and furthermore, there is the set $B_{z,\wz}$ of
$n_{z,\wz}$ coefficients of monomials
\mbox{$(i,j)\in\Del_{z,\wz}\backslash\Del_{z,\wz}^-$}, such that
$$\det\left(\frac{\partial\{\Phi^{(z,\wz)}_i(F),\
i=1,...,n_{z,\wz}\}}{\partial\{B_{z,\wz}\}}\right)\bigg|_{F=f_{z,\wz}}\ne
0\ .$$

\medskip

{\it Step 4}. We intend to write a formula for the desired
polynomial $f\in\K[x,y]$ with unknown coefficients, which then
will be found as a solution to certain system of equations.

For $k=1,...,N$, the restriction $\nu\big|_{\Del_k}$ is a linear
function $\lam_k(i,j)=\alp_ki+\bet_kj+\gam_k$. Introduce a
$\C$-linear map
$$T_k:\C[x,y,x^{-1},y^{-1}]\to\C[x,y,x^{-1},y^{-1}],\quad
T_k(x^iy^j)=x^iy^jt^{\lam_k(i,j)}\ .$$

Put
$$f(x,y)=\sum_{(i,j)\in\Del}a_{ij}x^iy^jt^{\nu(i,j)}+\sum_{k=1}^N\sum_{h\in B_k}
c_hT_k(h(x,y))$$
$$+\sum_{\renewcommand{\arraystretch}{0.6}
\begin{array}{c}
\scriptstyle{\{(k,\sig,z),(l,\widetilde\sig,\wz)\}\in\Pi}\\
\scriptstyle{\sig\subset\Del_k^+}
\end{array}}\sum_{\renewcommand{\arraystretch}{0.6}
\begin{array}{c}
\scriptstyle{m(k,z)i+mj<m\cdot m(k,z)}\\
\scriptstyle{(i,j)\ne(m-1,0)}
\end{array}}tc_{ij}^{(k,z)}T_k(\pi_{ij}^{(k,z)}(x,y))$$
\begin{equation}+\sum_{\renewcommand{\arraystretch}{0.6}
\begin{array}{c}
\scriptstyle{\{(k,\sig,z),(l,\widetilde\sig,\wz)\}\in\Pi}\\
\scriptstyle{\sig\subset\Del_k^-}
\end{array}}\sum_{\renewcommand{\arraystretch}{0.6}
\begin{array}{c}
\scriptstyle{m(k,z)i+mj<m\cdot m(k,z)}\\
\scriptstyle{j>0}
\end{array}}tc_{ij}^{(k,z)}T_k(\pi_{ij}^{(k,z)}(x,y))\
,\label{e18}\end{equation} where all the coefficients
$c_h=c_h(t)$, $c^{(k,z)}_{ij}=c^{(k,z)}_{ij}(t)$ are elements of
$\K$.

\medskip

{\it Step 5}. Pick $k=1,...,N$, and consider the polynomial
$$\hat f_k(x,y):=T_k^{-1}(f(x,y))=\sum_{(i,j)\in\Del_k}a_{ij}x^iy^j+\sum_{h\in B_k}c_hh(x,y)$$ $$+
\sum_{(i,j)\in\Del\backslash\Del_k}a_{ij}x^iy^jt^{\nu(i,j)-\lam_k(i,j)}
+\sum_{l\ne k}\sum_{h\in B_l}c_hT_k^{-1}T_l(h(x,y))$$
$$+\sum_{\renewcommand{\arraystretch}{0.6}
\begin{array}{c}
\scriptstyle{\{(p,\sig,z),(q,\widetilde\sig,\wz)\}\in\Pi}\\
\scriptstyle{\sig\subset\Del_p^+}
\end{array}}\sum_{\renewcommand{\arraystretch}{0.6}
\begin{array}{c}
\scriptstyle{m(p,z)i+mj<m\cdot m(p,z)}\\
\scriptstyle{(i,j)\ne(m-1,0)}
\end{array}}tc_{ij}^{(p,z)}T_k^{-1}T_p(\pi_{ij}^{(p,z)}(x,y))$$
\begin{equation}+\sum_{\renewcommand{\arraystretch}{0.6}
\begin{array}{c}
\scriptstyle{\{(p,\sig,z),(q,\widetilde\sig,\wz)\}\in\Pi}\\
\scriptstyle{\sig\subset\Del_p^-}
\end{array}}\sum_{\renewcommand{\arraystretch}{0.6}
\begin{array}{c}
\scriptstyle{m(p,z)i+mj<m\cdot m(p,z)}\\
\scriptstyle{j>0}
\end{array}}tc_{ij}^{(p,z)}T_k^{-1}T_p(\pi_{ij}^{(p,z)}(x,y))\
.\label{e69}\end{equation} This is a deformation of $f_k(x,y)$ in
${\cal P}(\Del)$.

Our first requirement about $\hat f_k(x,y)$ is that $\hat f_k\in
M_k$, where $M_k\subset{\cal P}(\Del)$ is the variety germ
introduced in Step 1. Consequently, by (\ref{e3}), this can be
expressed in the form
\begin{equation}c_h=L^k_h\left(\{c_{h'}\ :\ h'\in B_u,\
(u,k)\in\Arc(\Gam)\}\right) +O(t),\quad h\in B_k\
,\label{e9}\end{equation} where $L^k_h$ are linear functions with
constant coefficients in their variables $c_{h'}$.

\medskip

{\it Step 6}. Let $Z=\{(k,\sig,z),(l,\widetilde\sig,\wz)\}\in\Pi$,
$(C_k\cdot\Tor(\sig))_z=m\ge 2$. If $\widetilde\sig=\sig$,
$\wz=z$, then $\sig=\Del_k\cap\Del_l$, $z\not\in C_l^{nr}\cup
C_l^{nr}$. If $\widetilde\sig\ne\sig$, then the polygons $\Del_k$,
$\Del_l$ are joined by a well ordered sequence of polygons such
that any two neighboring polygons in the sequence have a common
edge, and all these common edges, among them $\sig$ and
$\widetilde\sig$, are parallel. We intend to specify conditions
under which there appear $m-1$ nodes of the curve $C_t$,
corresponding to the pair $Z$. We shall treat only the case
$\widetilde\sig\ne\sig$, since this includes the situation
$\widetilde\sig=\sig$ as a particular case.

Let $\Del_1,...,\Del_s$, $s\ge 1$, be the sequence of polygons
joining $\Del_k$ and $\Del_l$ (cf. Figure \ref{f5}(a)).
Geometrically, this means that the points $z$ and $\widetilde z$
are joined in $\bigcup_{i=1}^N\Tor(\Del_i)$ by a sequence of
non-singular rational components $C'_1\subset C_1$, ...,
$C'_s\subset C_s$, such that each component $C'_i$ appears in
$C_i$ with multiplicity $m$. Without loss of generality assume
that $\nu$ is constant on $\sig,\widetilde\sig$ and on all
parallel to them edges of $\Del_1,...,\Del_s$. Perform the
coordinate change $(x,y)\mapsto(x'',y'')$ as described in section
\ref{sec10}. The Newton polygons
$\Del''_k,\Del''_l,\Del''_1,...,\Del''_s$ will be located as shown
in Figure \ref{f5}(c). Denote the ordinates of the horizontal
edges of $\Del''_1,...,\Del''_s$ by $p_0>p_1>...>p_s$ (see Figure
\ref{f5}(f)). Then the vertices of $\theta(z,\wz)$ are
respectively $(0,p_s-m(l,\wz))$, $(0,p_0+m(k,z))$, $(m,p_0)$,
$(m,p_s)$. Divide $\theta$ by horizontal lines as shown in Figure
\ref{f5}(f). Denote by $\Theta$ the $\C$-linear map which takes
any polynomial $g(x,y)$ to $g''(x'',y'')$ along the above
coordinate change, and then projects $g''(x'',y'')$ to the space
\mbox{${\cal P}(\theta_0)$}, \mbox{$
\theta_0:=\Int(\theta)\cup[(0,p_s-m(l,\wz)+1),(0,p_0+m(k,z)-1)]$}.
Observe, that $\Theta$ induces the isomorphisms
$$\begin{cases}&\Span\left(\{\pi^{(k,z)}_{ij}\ :\ m(k,z)i+mj<m\cdot m(k,z),\ (i,j)\ne(m-1,0)\}\cup\left\{\partial
f''_k/\partial x''\right\}\right)\\
&\qquad\simeq\Span\{(x'')^\alp(y'')^\bet\ :\ \bet\ge p_0,\
m(k,z)\alp+m(\bet-p_0)<m\cdot m(k,z)\}\ ,\\
&\Span\{\pi^{(l,\wz)}_{ij}\ :\ j>0,\ m(l,\wz)i+mj<m\cdot m(l,\wz)\}\\
&\qquad\simeq\Span\{(x'')^\alp(y'')^\bet\ :\ \bet<p_s,\
m(l,\wz)\alp+m(p_s-\bet)<m\cdot m(l,\wz)\}\ ,\\
&\Span(\bigcup_{z\in
C'_u}B_{u,z})\simeq\Span\{(x'')^\alp(y'')^\bet\ :\ 0\le\alp<m,\
p_u\le\bet<p_{u-1}\}\ ,\\ &\qquad u=1,...,s\ .\end{cases}$$ The
latter isomorphism statement comes from the fact that the
monomials $(x'')^\alp(y'')^\bet$, $0\le\alp<m$,
$p_u\le\bet<p_{u-1}$, project to a basis of the space
$\bigoplus_{z\in C'_u}\ko_{\C^2,z}/I^{eg}(C_u,z)$, whereas
$\bigcup_{z\in C'_u}B_{u,z}$ by construction of Step 1 projects to
a basis of the aforementioned space $\bigoplus_{z\in C'_u}
\ko_{\C^2,z}/I^{eg}(C_u,z)$.

Put
$\widetilde\nu(i,j)=\max\{\lam''_k(i,j),\lam''_1(i,j),...,\lam''_s(i,j),\lam''_l(i,j)\}$.
Any coefficient $A_{ij}$, $(i,j)\in\theta_0$, of $f''(x'',y'')$,
satisfies $\val(A_{ij})<-\widetilde\nu(i,j)$. Furthermore, we have
$$\Theta(f''(x'',y''))=t\left(\sum_{(i,j)\in\theta_0}(L_{1,ij}+O(t))t^{\widetilde\nu(i,j)}+
L_2\cdot\frac{\partial f''}{\partial x''}\right)\ ,$$ where
$L_{1,ij}$ is a linear function with constant coefficients
depending on the parameters
\begin{equation}\begin{cases}&c^{k,z}_{\alp\bet},\quad\text{where}\
m(k,l)\alp+m\bet<m\cdot m(k,z),\ (\alp,\bet)\ne(m-1,0)\ ,\\
&c^{l,\wz}_{\alp\bet},\quad\text{where}\ m(l,\wz)\alp+m\bet<m\cdot m(l,\wz)\ ,\\
&c_h,\quad\text{where}\ h\in\bigcup_{z\in C'_1}B_{1,z}\cup...\cup
\bigcup_{z\in C'_s}B_{s,z}\ ,\end{cases}\label{e73}\end{equation}
and
\begin{equation}\left\{c_h,\quad\text{where}\
h\in\bigcup_{(u,k)\in\Arc(\Gam)}B_u\cup\bigcup_{(u,l)\in\Arc(\Gam)}B_u\cup\bigcup_{\alp=1}^s
\bigcup_{(u,\alp)\in\Arc(\Gam)}B_u\right\}\
,\label{e71}\end{equation} and $L_2$ is a linear function with
constant coefficients depending on the parameters
\begin{equation}\begin{cases}&\{c_h,\ h\in
B_k\},\quad\{c_h,\ h\in\bigcup_{(s,k)\in\Arc(\Gam)}B_s\},\\
&\{c^{(k,z)}_{ij}),\
m(k,z)i+mj<m\cdot m(k,z),\ (i,j)\ne(m-1,0)\},\\
&\{c^{(k,w)}_{ij},\ w\in C_k\cap\Tor(\partial\Del_k),\ w\ne z\},\\
&\{c^{(u,w)}_{ij},\
w\in\bigcup_{(u,k)\in\Arc(\Gam)}(C_u\cap\Tor(\partial\Del_u))\}\end{cases}\
. \label{e72}\end{equation}

Using the isomorphism induced by $\Theta$, we conclude that there
exist
\begin{equation}\begin{cases}&c^{(k,z)}_{ij},\quad m(k,z)i+mj<m\cdot m(k,z),\ (i,j)\ne(m-1,0)\
,\\ &c_h,\quad h\in\bigcup_{z\in C'_1}B_{1,z}\cup...\cup
\bigcup_{z\in C'_s}B_{s,z}\ ,\\
&c^{(l,\wz)}_{ij},\quad j>0,\ m(l,\wz)i+mj<m\cdot m(l,\wz)\
,\end{cases}\label{e78}\end{equation} and $\tau(t)\in\K$,
$\val(\tau)\le 0$ such that
\begin{equation}\Theta(f''(x''+t\tau(t),y'')=\sum_{(i,j)\in\theta_0}d^{z,\wz}_{ij}(t)t^{\nu_{z,\wz}(i,j)}(x'')^i(y'')^j\
,\label{e79}\end{equation} where $\val(d^{z,\wz}_{ij})\le 0$ and
the restriction of the tropicalization of $f''(x''+\tau(t),y'')$
on $\theta_0$ determines the subdivision of $\theta_0$ induced by
the function $\nu_{z,\wz}$ (see the definition in section
\ref{sec1} illustrated in Figure \ref{f5}(d,e)). Formally,
(\ref{e79}) reduces to a system of equations for the variables
(\ref{e78})
\begin{equation}\begin{cases}&c^{(k,z)}_{ij}=L^{(k,z)}_{ij}+O(t)\ ,\\ &\qquad m(k,z)i+mj<m\cdot m(k,z),\ (i,j)\ne(m-1,0)\
,\\ &c^{(l,\wz)}_{ij}=L^{(l,\wz)}_{ij}+O(t),\quad j>0,\
m(l,\wz)i+mj<m\cdot m(l,\wz)\
,\end{cases}\label{e80}\end{equation}
\begin{equation}c_h=L_h+O(t),\quad h\in\bigcup_{z\in C'_1}B_{1,z}\cup...\cup
\bigcup_{z\in C'_s}B_{s,z}\ ,\label{e90}\end{equation} where
$L^{(k,z)}_{ij},L_h,L^{(l,\wz)}_{ij}$ are linear function with
constant coefficients, whose variables are $d^{z,\wz}_{ij}$,
$(i,j)\in\theta_0$, and additionally
\begin{itemize}\item for $L^{(k,z)}_{ij}$,
\begin{equation}\begin{cases}&c_h,\quad\text{where}\quad
h\in\bigcup_{(s,k)\in\Arc(\Gam)}B_s\ ,\\
&c^{(k,w)}_{pq},\quad\text{where}\quad w\in\Tor(\Del_k^-),\
(C_k\cdot\Tor(\Del_k^-))_w\ge 2\
,\end{cases}\label{e81}\end{equation}
\item for $L_h$, $h\in B_u$, $1\le u\le s$,
\begin{equation}c_{pq},\quad\text{where}\quad (p,q)\in \Del_u^-\
,\label{e82}\end{equation} and \begin{equation}
\begin{cases}&c_h,\quad\text{where}\quad h\in B_k\ ,\\
&c_h,\quad\text{where}\quad h\in\bigcup_{(v,k)\in\Arc(\Gam)}B_v\ ,\\
&c^{(k,z)}_{pq},\quad\text{where}\quad
m(k,z)p+mq<m\cdot m(k,z)\ ,\\ &\qquad\qquad (p,q)\ne(m-1,0)\ ,\\
&c^{(k,w)}_{pq},\quad\text{where}\quad w\in C_k\cap\Tor(\partial\Del_k),\ w\ne z\ ,\\
&c^{(s,w)}_{pq},\quad\text{where}\quad
w\in\bigcup_{(v,k)\in\Arc(\Gam)}(C_v\cap\Tor(\partial\Del_v))\end{cases}
\label{e83}\end{equation}\item for $L^{(l,\wz)}_{ij}$,
\begin{equation}\begin{cases}&b_h(t),\quad\text{where}\quad
h\in\bigcup_{(v,l)\in\Arc(\Gam)}B_v\ ,\\
&c^{(v,w)}_{pq},\quad\text{where}\quad (v,l)\in\Arc(\Gam),\ w\in
C_l\cap\Tor(\Del_l^-) \ ,\end{cases}\label{e84}\end{equation} and
(\ref{e83}).
\end{itemize}

Our demands of $d^{z,\wz}_{ij}(t)$, $(i,j)\in\theta_0$ are as
follows. Let $\Del_Z$ be the only triangle in the subdivision of
$\theta_0$, and let $(m,p_u)$, $(0,p_u+m(k,z))$,
$(0,p_u-m(l,\wz))$ be its vertices. We write
$d^{z,\wz}_{ij}(t)=d^{z,\wz}_{ij}(0)+e^{z,\wz}_{ij}(t)$, where
$e^{z,\wz}_{ij}(0)=0$, and suppose that the tropicalization of
$f''(x''+\tau(t),y'')$ on any parallelogram of the subdivision of
$\theta_0$ is the product of a monomial and some irreducible
binomials, and that the tropicalization of $f''(x''+\tau(t),y'')$
on $\Del_Z$ is $(y'')^{p_u}P_{z,\wz}(x'',y'')$, i.e., the given
deformation pattern multiplied by a monomial. All this, clearly,
determines $d^{z,\wz}_{ij}(0)$, $(i,j)\in\theta_0$, uniquely.

Next we impose conditions on $e^{z,\wz}_{ij}(t)$,
$(i,j)\in\theta_0$. Namely, pick $v=1,...,s$ and consider the
parallelogram $\theta^{z,\wz}_u$ from the subdivision of $\theta$,
whose vertices are $$(m,p_{v-1}),\ (m,p_v),\ (0,p_{v-1}+m(k,z)),\
(0,p_v+m(k,z))$$ if $v\le u$ (i.e., the parallelogram lies above
the triangle $\Del_Z$, Figure \ref{f5}(e)), or
$$(m,p_{v-1}),\ (m,p_v),\ (0,p_{v-1}-m(l,\wz)),\ (0,p_v-m(l,\wz))$$ if
$v>u$ (i.e., the parallelogram lies below the triangle $\Del_Z$,
Figure \ref{f5}(e)). Let
$\nu_{z,\wz}(i,j)\big|_{\theta^{z,\wz}_v}=\alp_vi+\bet_vj+\gam_v$.
Then
\begin{equation}t^{-\gam_v}f''(x''t^{-\alp_v}+t\tau(t),y''t^{-\bet_v})=\sum_{(i,j)\in\theta^{z,\wz}_v}
d^{z,\wz}_{ij}(t)(x'')^i(y'')^j+O(t)\ .\label{e75}\end{equation}
The tropicalization of the latter polynomial on $\theta^{z,\wz}_v$
defines a curve $C^{z,\wz}_v\subset\Tor(\theta^{z,\wz}_v)$ which
consists of components defined by binomials. We demand that the
curve, defined by the polynomial (\ref{e75}), belongs to the
variety $M^{eg}(C^{z,\wz}_v)$, introduced in Lemma \ref{l21}(ii).
In turn, Lemma \ref{l21}(ii) yields that this can be expressed by
a system of equations
\begin{equation}e^{z,\wz}_{ij}=L^{z,\wz}_{ij}+O(t),\quad
(i,j)\in\theta^{z,\wz}_{v,0}\ ,\label{e76}\end{equation} where
$\theta^{z,\wz}_{v,0}$ is obtained from $\theta^{z,\wz}_v$ by
removing its upper and right edges, if $v\le u$, or removing the
lower and right edges, if $v>u$ (see Figure \ref{f5}(e)), and
$L^{z,\wz}_{ij}$ are linear function with constant coefficients
and variables
$e^{z,\wz}_{i'j'}:=d^{z,\wz}_{i'j'}-d^{z,\wz}_{i'j'}(0)$ as
$(i',j')$ ranges over
$\theta_0\cap(\theta^{z,\wz}_v\backslash\theta^{z,\wz}_{v,0})$ and
$c_h$, $h\in\bigcup_{(q,v)\in\Arc(\Gam)}B_q$.

At last, let $\nu_{z,\wz}(i,j)\big|_{\Del_Z}=\alp i+\bet j+\gam$.
Then the polynomial
$$t^{-\gam}f''(x''t^{-\alp}+\tau(t),y''t^{-\bet})=\sum_{(i,j)\in\Del_Z}
d^{z,\wz}_{ij}(t)(x'')^i(y'')^j+O(t)$$ represents a deformation of
$(y'')^{p_u}P_{z,\wz}(x'',y'')$, which we want to be
$\ks$-equisingular with respect to the singularities of
$\{P_{z,\wz}=0\}$ in $\C^2$. As pointed in Step 4, this can be
expressed by a system of equations
\begin{equation}e^{z,\wz}_{ij}=L^{z,\wz}_{ij}+O(t),\quad (i,j)\in
B_{z,\wz}+p_u\ ,\label{e77}\end{equation} where $L^{z,\wz}_{ij}$
are linear functions with constant coefficients and variables
$e^{z,\wz}_{i'j'}$, $(i',j')\in\theta_0\backslash(B_{z,\wz}+p_u)$,
and $c_h$, $h\in\bigcup_{(q,u)\in\Arc(\Gam)}B_q$ .

The variables in systems (\ref{e76}), $v=1,...,s$, and (\ref{e77})
are naturally ordered so that, for $t=0$, each variable depends
linearly only on the preceding variables; hence by the implicit
function theorem this bunch of equations can be resolved with
respect to $e^{z,\wz}_{ij}$, $(i,i)\in\theta_0$. We then plug the
solution obtained to system (\ref{e80}), (\ref{e90}), noticing
that, in this substitution, the variables $c^{(q,w)}_{ij}$,
mentioned in (\ref{e84}), enter the terms $O(t)$ for all
$q,w,i,j$.

\medskip

{\it Step 7}. Before we join all the equation obtained in the
preceding steps into one system, we should like to notice that
some equations may be dependent, and hence must be removed from
the system, since we finally intent to apply the implicit function
theorem. Namely, the system of equations (\ref{e80}), (\ref{e90}),
obtained in Step 6, is, in fact, included in the system
(\ref{e80}), (\ref{e9}). Indeed, in our setting, (\ref{e9}) takes
form
\begin{equation}c_h=L^q_h\left(\{c_{h'}\ :\ h'\in B_v,\
(v,q)\in\Arc(\Gam)\}\right) +O(t),\quad h\in B_q,\ q=1,...,s\
.\label{e91}\end{equation} By the implicit function theorem we can
resolve system (\ref{e80}), (\ref{e90}) with respect to the
variables in the left-hand side, then we substitute the
expressions for $c^{k,z}_{ij}$, $c^{l,\wz}_{ij}$ into (\ref{e90}).
The right-hand sides of the resulting system (\ref{e90}) depend on
the same bunch of variables as in (\ref{e91}), and, by our
construction, system (\ref{e91}) implies the property that the
distinct multiple components of any of the curves $C_1,...,C_s$ do
not glue up with each other and with any other component in a
neighborhood of $\bigcup_{i=1}^s(\Sing(C^{\text{\rm
red}}_i)\backslash\Sing^{\text{\rm iso}}(C_i)$ along the
deformation defined by $f(x,y)$. In turn system (\ref{e90}) simply
expresses the latter property for some of the multiple components
of $C_1,...,C_s$. Hence the claim follows, and we get rid of all
equations (\ref{e90}), including instead equations (\ref{e91}) in
the final system.

\medskip

{\it Step 8}. All the conditions imposed on the required
polynomial $f(x,y)$, we have expressed as systems of equations
\begin{itemize}\item (\ref{e9}) for all $k=1,...,N$, \item
(\ref{e80}) for all pairs
$\{(k,\sig,z),(l,\widetilde\sig,\wz)\}\in\Pi$.
\end{itemize}
The orientation $\Gam$ induces an ordering of the variables in the
above united system such that, for $t=0$, each variable is
expressed only via strongly preceding variables, and hence the
system can be solved by the implicit function theorem.

Geometric meaning of the imposed conditions is that $f(x,y)$
induces a $\ks$-equisingular one-parametric deformation for each
point $z\in\Sing^{\text{\rm iso}}(C_k)\cap(\C^*)^2$, for all
$k=1,...,N$, and for each point $z\in\Sing(C_{z,\wz})$,
$\{(k,\sig,z),(l,\widetilde\sig,\wz)\}\in\Pi$. Furthermore, each
point $z\in\Sing(C^{\text{\rm red}}_k)\backslash\Sing^{\text{\rm
iso}}(C_k)$, $1\le k\le N$, bears $\dim\ko_{\C^2,z}/I^{eg}(C_k,z)$
nodes, because the curves $C^{(t)}\subset\Tor(\Del)$ have no
multiple components (the curve $C^{(t)}$ crosses
$\Tor(\partial\Del)$ with multiplicity $1$ at each point by
assumptions of section \ref{sec1}). At last, notice that $C^{(t)}$
has no other singular points, for example, a point
$z\in\Tor(\sig)\cap C_k\cap C_l$, $\sig=\Del_k\cap\Del_l$, with
$(C_k\cdot\Tor(\sig))_z=1$ bears no singular points in view of
Lemma \ref{l6}.

\medskip

{\it Step 9}. We complete the last task, explaining that the
polynomials $f(x,y)\in\K[x,y]$, constructed in the preceding
steps, contain a family described by (\ref{e2n}), (\ref{e3n}).

To obtain the required family, it is sufficient to show that the
polynomials $h\in\bigcup_{k=1}^NB_k$ and all the polynomials
$\pi_{ij}^{(k,z)}$, introduced in Step 2, can be chosen so that
they do not contain monomials $x^\alp y^\bet$, $(\alp,\bet)\in B$.
Indeed, if the latter holds, then the solution to the system
considered in Step 8 depend on the coefficients $c_\omega$,
$\omega\in B$, as free parameters.

Pick $k=1,...,N$. If $B\cap\Del_k\subset\Del^-_k$, then by
construction, the polynomials $h\in B_k$ and the polynomials
$\pi_{ij}^{(k,z)}$, $j>0$, do not contain the monomials $x^\alp
y^\bet$, $(\alp,\bet)\in B$. Assume that $|B\cap\Del_k|\le 3$.
Then by transformations $g(x,y)\in\C[x,y]\mapsto ag(bx,cy)$,
$a,b,c\in\C^*$, we can freely vary the coefficients of $x^\alp
y^\bet$, $(\alp,\bet)\in B\cap\Del_k$, in any polynomial $g(x,y)$
with Newton polygon $\Del_k$. On the other hand, all the strata
$M^\ks(C_k,z)$, $M^{eg}(C_k,z)$, and $M^\sqh(C_k,z)$, which appear
in Definition \ref{d1}, are invariant with respect to the above
transformations (close to the identity). Hence (cf. the proof of
Lemmas \ref{l2} and \ref{l3}) the polynomials $h\in B_k$ and all
$\pi_{ij}^{(k,z)}$ can be chosen free of the monomials $x^\alp
y^\bet$, $(\alp,\bet)\in B$. \proofend

\subsection{Proof of Lemma \ref{l18}}\label{sec300}
Let $(A,S,F,R)\in{\cal Q}(nA_1)$. Take any vector
$\zeta\in\R^2\backslash\{0\}$ which in not parallel to any of the
edges of $S$, and orient the arcs of $A$ so that they form acute
angles with the chosen vector. This defines an orientation $\Gam$
of $A$, meeting the requirements of Theorem \ref{t1}. Furthermore,
the deformation patterns $R$ and all the triads
$(\Del_k,\Del_k^-,C_k)$ are transversal according to Lemmas
\ref{l2} and \ref{l3}. For example, if $\Del_k$, $1\le k\le N$, is
a triangle, the inequality of Lemma \ref{l2}(i), serving as the
transversality criterion, holds true, since the $b$-invariant
vanishes for nodes, $C_k$ is non-singular along
$\Tor(\partial\Del_k)$, and $\eps=1$ for all edges in
$\Del_k^+\ne\emptyset$. Thus, Theorem \ref{t1} applies, and the
set $B$ can be chosen as follows. For any parallelogram $\Del_k$
the set $\Del_k^+$ is the union of two neighboring edges. Then we
take $V(S)$ and remove all the interior vertices of $\Del_k^+$ for
all parallelograms $\Del_k$.

Notice that $|B|=r+1$, that is, formulas (\ref{e2n}), (\ref{e3n})
describe all the polynomials $f\in\K[x,y]$ with Newton polygon
$\Del$, defining $n$-nodal curves in $\Lam_K(\Del)$ and
tropicalizing into $(A,S,F,R)$.

We separate equations for $c_{ij}$, $(i,j)\in V(S)\backslash B$,
from (\ref{e3n}). Namely, let
$\bi_{1,k},\bi_{2,k},\bi_{3,k},\bi_{4,k}$ be the vertices of a
parallelogram $\Del_k$, listed clockwise, and
$\bi_{1,k}\not\in\Del_k^-$. Then an equation for $M^{eg}(C_k)$
(see section \ref{sec3}), involving the coefficients at the
vertices of $\Del_k$, reads
$$(a_{\bi_{1,k}}+c_{\bi_{1,k}})
(a_{\bi_{3,k}}+c_{\bi_{3,k}})-(a_{\bi_{2,k}}+c_{\bi_{2,k}})(a_{\bi_{4,k}}+c_{\bi_{4,k}})=O(t)\
,$$ where $O(t)$ includes the terms containing $t$ to a positive
power. Since
\mbox{$a_{\bi_{1,k}}a_{\bi{3,k}}-a_{\bi_{2,k}}a_{\bi_{4,k}}=0$},
we can rewrite system (\ref{e3n}) in the form
\begin{equation}c_{ij}=\Phi_{ij}(\{c_{\bi}\ :\ \bi\in
V(S)\}),\quad (i,j)\in\Del\backslash V(S)\
,\label{e95}\end{equation}\begin{equation}\begin{cases}&c_{\bi_{1,k}}
a_{\bi_{3,k}}+c_{\bi_{3,k}}a_{\bi_{1,k}}-c_{\bi_{2,k}}a_{\bi_{4,k}}-c_{\bi_{4,k}}a_{\bi_{2,k}}=O(t)\
,\\ &\qquad 1\le k\le N,\ |V(\Del_k)|=4\
.\end{cases}\label{e96}\end{equation}

Consider now equations $f(\bp_1)=...=f(\bp_r)=0$. Let
$\bx_s=\val(\bp_s)$ correspond to an edge $\sig_s$ of $S$. Without
loss of generality, suppose that $\sig_s$ lies on the horizontal
coordinate axis, $\nu\big|_{\sig_s}=0$, $\nu(i,j)>0$ as
$(i,j)\not\in\sig_s$. Then
$\bp_s=(\xi_s^0+\xi_s^1t,\eta_s^0+\eta_s^1t)$, where
$\xi_s^0,\eta_s^0\in\C^*$, $\xi_s^1,\eta_s^1\in\K$,
$\val(\xi_s^1),\val(\eta_s^1)\le 0$.

Assume that $|\sig_s|=1$, i.e., $\sig_s=[\bi'_s,\bi''_s]$,
$\bi'_s=(i,0),\bi''_s=(i+1,0)$. Then the equation $f(\bp_s)=0$ in
the form (\ref{e146}) reads
$$(a_{\bi'_s}+c_{\bi'_s})+(a_{\bi''_s}+c_{\bi''_s})(\xi_s^0+\xi_s^1t)=O(t)\
,$$ that, in view of $a_{\bi'_s}+a_{\bi''_s}\xi_s^0=0$, transforms
into $$c_{\bi'_s}+c_{\bi''_s}\xi_s^0=O(t)\ .$$

Assume that $|\sig_s|=m\ge 2$, i.e, without loss of much
generality, $\sig_s=[\bi'_s,\bi''_s]$,
$\bi'_s=(0,0),\bi''_s=(m,0)$. We consider only the case when
$\sig_s$ is a common edge of two triangles $\Del_k,\Del_l$ (cf.
Figure \ref{f4}), since the situation when $\sig_s$ is an edge of
a parallelogram can be treated in the same way, but requires a
more complicated notation. Let $z=\Tor(\sig_s)\cap C_k\cap C_l$.
We have
$$f(x,y)=\sum_{i=0}^m(a_{i,0}+c_{i,0})x^i+O(t)=\sum_{i=0}^mc_{i,0}x^i+a_{m,0}(x+\xi_s^0)^m+O(t)\ .$$ The coordinate
change $x=x'+\xi_s^0$ takes $f(x,y)$ into
$$f'(x',y)=\sum_{i=0}^{m-1}c'_{i,0}(x')^i+(a_{m,0}+c_{m,0})(x')^m+O(t)\
,$$ where
$$c'_{i,0}=\sum_{j=i}^mc_j(\xi_s^0)^{j-i}\left(\begin{matrix}j\\
i\end{matrix}\right),\quad i=0,...,m-1\ .$$ Furthermore, there is
$$\tau=-\frac{c'_{m-1,0}}{ma_{m,0}}+O(t)+\text{h.o.t.}\ \in\K\ ,$$
where ``h.o.t." contains all monomials in $c_\bi$, $\bi\in\Del$,
of degree $\ge 2$, such that the polynomial
$f''(x'',y):=f'(x'+\tau,y)$ has zero coefficient of $x^{m-1}$.
Then
$$f''(x'',y)=\sum_{i=1}^{m-2}c''_{i,0}(x'')^i+(a_{m,0}+c''_{m,0})(x'')^m+y(a''_{01}+c''_{01})t^p+
y^{-1}(a''_{0,-1}+c''_{0,-1})t^q+...\ .$$ where we omit monomials
$(x'')^iy^j$ with $(i,j)\not\in\Del_z$,
$\Del_z=\conv\{(m,0),(0,1),(0,-1)\}$, and have
$$c''_{i,0}=c'_{i,0}+O(t)+\text{h.o.t.},\ i=0,...,m-2,\quad c''_{m,0}(0)=c''_{01}(0)=c''_{0,-1}(0)=0\ ,$$
whereas $a''_{01},a''_{0,-1}\in\C^*$, and $p,q$ are distinct
positive integers, about which we assume $p<q$. By Lemma
\ref{l113}, the tropicalization of $f''(x'',y)$ determines a
subdivision, containing the triangle $\Del_z$, and the
corresponding deformation pattern. In particular,
$c''_{i,0}=O(t)$, $i=0,...,m-2$. Another consequence is that,
plugging the coordinates $x''=t\xi_s^1-\tau$,
$y=\eta_s^0+\eta_s^1t$ of $\bp_s$ into $f''(x'',y)$, we obtain
that the minimal powers of $t$ appear from monomials $x^m$ and
$y$, and they must compensate each other, since the coordinates of
$\bp_s$ annihilate $f''$, that is
$$\eta_s^0a''_{01}t^p+a_{m,0}(\xi_s^1t-\tau)^m+\text{h.o.t.}=0\
.$$ The latter equation leads to
$$\tau=\xi_s^1t-\left(-\frac{\eta_s^0a''_{01}}{a_{m,0}}\right)^{1/m}t^{p/m}+\text{h.o.t.}\
.$$ Combining this with the above formulas for $\tau$,
$c''_{i,0}$, $c'_{i,0}$, it is not difficult to derive the
relation
$$c_{00}-\frac{a_{00}}{a_{m,0}}c_{m,0}=c_{00}+(-1)^{m+1}(\xi_s^0)^mc_{m,0}=(-1)^mma_{m,0}(\xi_s^0)^{m-1}\xi_s^1t$$
\begin{equation}+
(-1)^{m-1}m(\xi_s^0)^{m-1}(-\xi_s^0a''_{01}a_{m,0}^{m-1})^{1/m}t^{p/m}
+\Phi_s\ ,\label{e98}\end{equation} where $\Phi_s$ is some
analytic function of parameters $a_{ij}$, $c_{ij}$,
$(i,j)\in\Del$, $\xi_s^0$, $\xi_s^1$, $\eta_s^0$, $\eta_s^1$, $t$,
whose terms contain $t$ to a positive power, or variables $c_{ij}$
to the total power $\ge 2$. We point out that formula (\ref{e98})
gives $m$ distinct equations.

Thus, we finally obtain $\prod_{s=1}^r|\sig_s|$ distinct systems
of equations for the coefficients of $f(x,y)$. Each system
consists of equations (\ref{e95}), (\ref{e96}) and (\ref{e98}),
the latter one being in the form
\begin{equation}c_{\bi'_s}a_{\bi''_s}-c_{\bi''_s}a_{\bi'_s}=O(t)+\text{h.o.t.},
\quad s=1,...,r\ .\label{e99}\end{equation} We now put
$c_{i_0j_0}=0$ for some $(i_0,j_0)\in B$, and apply the implicit
function theorem in order to conclude that the system has a unique
solution. The conditions of the implicit function theorem are
fulfilled, since, for example, the independence of the linearized
system (\ref{e96}), (\ref{e99}) for $t=0$ is equivalent to that of
system (\ref{e144}), (\ref{e145}) treated in section \ref{sec11},
Step 1.

\section{Counting real nodal curves}\label{sec301}

The complex conjugation naturally acts in $\K$, which allows us to
speak on real $\K$-curves, i.e., defined over the subfield $\K_\R$
of Puiseux series with real coefficients. If, for example, the
given points $\bp_1,...,\bp_r$ belong to $(\K_\R^*)^2$, using the
formulas from Lemmas \ref{l7}, \ref{l11}, \ref{l10}, one can count
how many real tropicalizations $(A,S,F,R)\in{\cal Q}_\Del(nA_1)$
correspond to a nodal amoeba $A$. Then, taking the real solutions
of equations (\ref{e98}), we can decide how many real nodal curves
correspond to a given amoeba, and thereby confirm the formulae
suggested by Mikhalkin in \cite{M2}. Here we focus on a related
problem of computing the Welschinger number
$\chi_\Del(\bp_1,...,\bp_r)$, introduced in a general symplectic
setting in \cite{Wel}. In our situation it is the number of real
nodal irreducible curves passing through the given real points and
counted with the sign $(-1)^{n^{\text{\rm sol}}}$, where, for a
given real nodal curve, $n^{\text{\rm sol}}$ is the number of its
real solitary nodes (i.e., locally given by $x^2+y^2=0$). The
importance of this number comes from Welschinger's theorem
\cite{Wel} that, for rational nodal curves,
$\chi_\Del(\bp_1,...,\bp_r)$ does not depend on the choice of the
fixed points. This means, when calculated for a special
configuration of real $r=|\partial\Del\cap\Z^2|-1$ points,
$|\chi_\Del(\bp_1,...,\bp_r)|$ provides a lower bound for the
number of real rational curves passing through an {\it arbitrary}
collection of $r$ real generic points in $\Tor(\Del)$.

As a consequence of the results of the preceding sections we state

\begin{proposition}\label{p1} In the notation of section
\ref{sec302}, given generic points $\bx_1,...,\bx_r\in\Q^2$ and
$\bp_1,...,\bp_r\in(\K_\R^*)^2$ such that $\val(\bp_i)=\bx_i$,
$i=1,...,r$, and an irreducible nodal amoeba $A$ of rank $r$,
passing through $\bp_1,...,\bp_r$, the following holds:
\begin{enumerate}\item[(i)] if the dual subdivision $S$ contains
an edge of even length, the contribution to
$\chi_\Del(\bp_1,...,\bp_r)$ of real $n$-nodal curves, passing
through $\bp_1,...,\bp_r$ and projecting onto $A$, is zero;
\item[(ii)] if the dual subdivision $S$ has only edges of odd
length, there exists a unique real irreducible $n$-nodal curve,
passing through $\bp_1,...,\bp_r$ and projecting onto $A$, and its
contribution to $\chi_\Del(\bp_1,...,\bp_r)$ is $(-1)^s$, where
$s$ is the total number of interior integral points in the
triangles of $S$.
\end{enumerate}
\end{proposition}

{\bf Proof}. Recall that by Lemma \ref{l30}, to count irreducible
nodal curves, we have to consider only irreducible nodal amoebas.

Let $S$ contain an edge $\sig$ of even length $m$. Given a real
tropicalization $f_1,...,f_N$, by the formulas of Lemma \ref{l10},
we can associate with the edge $\sig$ either zero, or two real
deformation patterns, which in turn are independent on how many
real solutions equations (\ref{e98}) have. If the real deformation
patterns do exist, their explicit formulas can be extracted from
the computation in the proof of Lemma \ref{l10}. Namely, one real
deformation pattern corresponds to the Chebyshev polynomial
$P(x)=\cos(m\cdot\arccos(2^{-(m-1)/m}x))$, and this deformation
pattern has $m-1$ real solitary nodes by \cite{Sh1}, Proposition
2.5. The other real deformation pattern corresponds to the
polynomial $-P(x\sqrt{-1})$, and it has one non-solitary node
besides $m-2$ imaginary nodes. Thus, claim (i) follows, since an
exchange of the above deformation patterns changes the parity of
the number of solitary real nodes.

If $S$ contains only edges of odd length, then the formulas in the
proof of Lemmas \ref{l7}, \ref{l10} and equations (\ref{e98}) give
a unique real choice for an $n$-nodal curve through
$\bp_1,...,\bp_r$, projecting onto $A$. It is a simple exercise to
check that the real tropicalizations to triangles and real
deformation patterns associated with edges of odd length have only
imaginary or real solitary nodes, whereas the real
tropicalizations to parallelograms do not bear solitary real
nodes. Thus, statement (ii) follows. \proofend

{\it Address}: School of Mathematical Sciences, Tel Aviv
University, Ramat Aviv, 69978 Tel Aviv, Israel.

{\it E-mail}: shustin@post.tau.ac.il
\end{document}